# The monodromy groups of Schwarzian equations on closed Riemann surfaces

By Daniel Gallo, Michael Kapovich, and Albert Marden

*To the memory of Lars V. Ahlfors*

## Abstract

Let $\theta : \pi_1(R) \rightarrow \mathrm{PSL}(2, \mathbb{C})$ be a homomorphism of the fundamental group of an oriented, closed surface $R$ of genus exceeding one. We will establish the following theorem.

THEOREM. *Necessary and sufficient for $\theta$ to be the monodromy representation associated with a complex projective stucture on $R$, either unbranched or with a single branch point of order 2, is that $\theta(\pi_1(R))$ be nonelementary. A branch point is required if and only if the representation $\theta$ does not lift to* $\mathrm{SL}(2, \mathbb{C})$.

## Contents

## 1. Introduction and background

1.1. *Introduction.* The goal of this paper is to present a complete, self-contained proof of the following result:

THEOREM 1.1.1. *Let $R$ be an oriented closed surface[1] of genus exceeding one, and*

$$\theta : \pi_1(R; O) \to \Gamma \subset \mathrm{PSL}(2, \mathbb{C})$$

*a homomorphism of its fundamental group onto a nonelementary group $\Gamma$ of Möbius transformations. Then*:

(i) *$\theta$ is induced by a complex projective structure for some complex structure on $R$ if and only if $\theta$ lifts to a homomorphism*

$$\theta^* : \pi_1(R; O) \to \mathrm{SL}(2, \mathbb{C}).$$

(ii) *$\theta$ is induced by a branched complex projective structure with a single branch point of order two for some complex structure on $R$ if and only if $\theta$ does not lift to a homomorphism into $\mathrm{SL}(2, \mathbb{C})$.*

The terms will be explained in §§1.2–1.4.

Theorem 1.1.1 characterizes the class of groups arising as monodromy groups of Schwarzian differential equations or equivalently, of the projectivized monodromy groups for the associated linear second-order differential equations. Poincaré himself explicitly raised the question by noting (for punctured spheres) second-order equations depend on the same number of parameters as their monodromy groups (the position of the singularities–the conformal structure–is allowed to change) and from this observation boldly concluded, "On peut *en général* trouver une équation du 2d ordre, sans points à apparence singulière qui admette un groupe donné" [P, p. 218]. In our own time, the question was raised in [Gu3] and [He1]; in fact Gunning conjectured Part (i) of our theorem and Tan [Tan] conjectured Part (ii).

Schwarzian equations themselves have long been an important tool in the study of Riemann surfaces and their uniformization. Their relation with algebraic geometry was established by Gunning in [Gu1]: For a fixed complex structure on $R$, the linear monodromy representations of the complex projective structures correspond to flat *maximally unstable* rank 2 holomorphic vector bundles over $R$. A similar relation for branched structures was later studied by Mandelbaum; see e.g. [Man 2, 3] (see also §11).

---

[1] In this paper, all surfaces are assumed to be connected. A *closed* surface is one which is compact, without boundary.

In §11, we will present an analogue Theorem 11.3.3 of our main theorem in the context of holomorphic vector bundles over Riemann surfaces. Namely, let $S$ be an oriented closed surface of genus exceeding one and $\rho : \pi_1(S) \to \mathrm{SL}(2, \mathbb{C})$ a nonelementary representation. Then $\rho$ is the monodromy of a holomorphic flat connection on a *maximally unstable holomorphic vector bundle* of rank two over a Riemann surface $R$, where $R$ is diffeomorphic to $S$ via an orientation preserving diffeomorphism $R \to S$.

Besides the fuchsian groups of uniformization, the class of monodromy groups includes the discrete, isomorphic groups of quasifuchsian deformations (Bers slices which model Teichmüller spaces and their boundaries), and discrete groups such as Schottky groups which are covered by fuchsian surface groups. See [Mas2] for a wide array of possibilities.

Theorem 1.1.1 further implies that the image in $\mathrm{PSL}(2, \mathbb{C})$ of "almost" every homomorphism of the fundamental group has a geometric structure. This is quite astonishing, especially so as the image groups are often not discrete and not even finitely presentable.

R. Rubinsztein [R] observed that if $G_0 \subset G = \pi_1(R)$ is *any* index two subgroup, the restriction of $\theta$ to $G_0$ can be lifted from $\mathrm{PSL}(2, \mathbb{C})$ to $\mathrm{SL}(2, \mathbb{C})$ in $2^{2g}$ ways. Consequently by Theorem 1.1.1, a homomorphism whose restriction to an index two subgroup is nonelementary is always associated with a complex projective structure for some complex structure on the corresponding two sheeted cover. One such index two subgroup is constructed in §8.6.

Special cases of Theorem 1.1.1(i) were proved in [He1] and the case of homomorphisms into $\mathrm{PSL}(2, \mathbb{R})$ was investigated in [Ga], [Go], [Por] and [Tan].

Proofs of Theorem 1.1.1(i) have been announced before. Gallo's research announcement [Ga1] proposed an innovative strategy for a proof, but the promised details have not been published or confirmed. Gallo's strategy had been developed in consultation with W. Goldman and W. P. Thurston, and was particularly inspired by Thurston's approach to the deformation of fuchsian groups by bending. Goldman's paper [Go1] is an exemplar of this strategy applied in the interesting special case where $\theta$ is an isomorphism onto a fuchsian group; it deals with the problem of determining *all* complex projective structures with the prescribed monodromy. This question is discussed further in §12.

The recent paper [Ka] proposed a proof confirming Theorem 1.1.1(i). Although the argument presented is incomplete (Lemmas 1 is incorrect and a condition is omitted in Lemma 4; they are corrected in the present paper, and some details are missing in the proofs of Propositions 1 and 2), the paper contains new ideas and directly motivated a fresh examination of the whole issue.

The present work was begun by Marden with the goal of settling the validity of the claims. In a general sense, Gallo's and Kapovich's strategy is followed, although the details, especially in Part B, are quite different from those suggested in [Ga1] or [Ka]. In the latter phase of the investigation, a collaboration with Kapovich began. Almost immediately this produced a breakthrough in understanding the connection between a certain construction invariant and the lifting obstruction (§§9–10). Instead of using the difficult continuity arguments proposed in [Ka], we use *branched* structures. Motivated by Tan's work [Tan] on real branched structures, we found a technique for constructing branched projective structures complementing that developed earlier for joining pants. This approach exhibits clearly the connection. It also clarifies the role of the 2nd Stiefel-Whitney class and degree of instability of holomorphic bundles which is discussed by Kapovich in §11. In fact, one of our discoveries is that it is *easier* to prove Theorem 1.1.1 simultaneously for branched and unbranched structures than to establish the unbranched case by itself.

Part C of our paper brings together additional results that fill out the picture presented by our main theorem. These are developed in the context of holomorphic bundles over Riemann surfaces. For example, in some respects Theorem 1.1.1 is more clearly seen in the context of a more general existence theorem for *branched* complex projective structures with a prescribed branching divisor and monodromy representation. This refinement, Theorem 11.2.4, is expressed in terms of the 2nd Stiefel-Whitney class. In addition, we present the full proof of the divergence theorem briefly outlined in [Ka]. This Theorem 11.4.1 deals with sequences of monodromy homomorphisms $\theta_n : \pi_1(R) \to \mathrm{PSL}(2, \mathbb{C})$ associated with divergent Schwarzian equations on a *fixed* Riemann surface. Such a sequence of homomorphisms *cannot* converge algebraically to a homomorphism, either nonelementary *or* elementary. In terminology of Teichmüller theory, the extension of a Bers slice to the full representation variety is properly embedded. In §12 we list and briefly discuss a number of open problems arising from our work.

We three authors decided to join together to pool the fruits of a decade of our individual and collaborative research relating to the main result. By doing so we have arrived at a rather larger understanding of the fundamental existence problem for the monodromy of projective structures.

Our topic falls under the ancient and revered subject heading of linear ordinary differential equations on Riemann surfaces, a subject introduced by Poincaré. The problem we consider fits comfortably with those associated with "the Riemann-Hilbert Problem" (Hilbert's 21st problem) for first-order fuchsian systems and $n$-th order fuchsian equations. Yet our approach is quite different than that associated with this theory [A-B], [I-K-Sh-Y], [Sib], [Y], [He2]. For one thing, our approach is special to second order equations. Then we work primarily with *projectivized* monodromy in $\mathrm{PSL}(2, \mathbb{C})$. This turns the problem into one largely involving the geometry of surfaces and Möbius groups.

Another difference is that here we are mainly dealing with equations without singularities. Finally we do not prescribe the complex structure in advance, rather it is determined as part of the solution: the number of parameters in the equations matches the number in the representations. The need to introduce a branch point to handle part (ii) of our theorem is however reminiscent of the need for "apparent singularities" in that theory.

Except for a particular case, we have left aside the general existence problem for surfaces with punctures and branch points. However, we believe that the foundation laid here will stimulate (further) exploration of these and other important aspects of the subject, including a characterization of the nonuniqueness, that are not now well understood.

*Acknowledgments.* Marden would like to thank the Mathematics Institute of the University of Warwick, the Forschungsinstitut für Mathematik at ETH, Zürich, and the Mathematical Sciences Research Institute in Berkeley, for the privilege of participating in their programs while his research was carried out. In addition he thanks David Epstein, Dennis Hejhal, Yasutaka Sibuya, and Kurt Strebel for helpful discussions. David in particular provided insightful suggestions for some of the proofs.

This research additionally received support from the NSF grants DMS-9306140 and DMS-96-26633 (Kapovich) and DMS-9022140 at MSRI (Kapovich and Marden).

All us authors thank Silvio Levy for providing invaluable editorial and LaTeX assistance and the referee for many helpful comments and suggestions.

### 1.2. *Möbius transformations.*

Möbius transformations correspond to elements of $\mathrm{PSL}(2,\mathbb{C})$ according to

$$\alpha(z) = \frac{az+b}{cz+d} \qquad \longleftrightarrow \qquad \pm \begin{pmatrix} a & b \\ c & d \end{pmatrix} \quad \text{with } ad - bc = 1.$$

They extend from their action on the extended plane $\mathbb{C} \cup \infty$ to upper half-three-space or, via stereographic projection, from the 2-sphere $\mathbb{S}^2$ to the 3-ball. The extensions form the group of orientation-preserving isometries of hyperbolic three-space, which we denote by $\mathbb{H}^3$ (in either the ball model or the upper half-space model) with $\partial \mathbb{H}^3$ denoting the "sphere at infinity," that is, the extended plane or $\mathbb{S}^2$, depending on the model. Throughout our paper, we will identify the extended plane with $\mathbb{S}^2$.

We recall the standard classification:

- A transformation $\alpha$ is *parabolic* if it has exactly one fixed point on $\partial \mathbb{H}^3$, or, equivalently, if it is not the identity and its trace satisfies $\mathrm{tr}^2 \alpha = (a+d)^2 = 4$. Parabolic transformations are those conjugate to $z \mapsto z+1$.

- An *elliptic* transformation has two fixed points in $\partial\mathbb{H}^3$ and also fixes pointwise its *axis of rotation*, that is, the hyperbolic line in $\mathbb{H}^3$ joining the fixed points. Its trace satisfies $0 \leq \mathrm{tr}^2\alpha < 4$, and it is conjugate to an element of the form $z \mapsto e^{2i\theta}z$, for $0 < \theta < \pi$.

- A *loxodromic* transformation $\alpha$ likewise has two fixed points in $\partial\mathbb{H}^3$, one repulsive and the other attractive; it preserves the line in $\mathbb{H}^3$ between them which is called the *axis*. The trace of $\alpha$ satisfies $\mathrm{tr}^2\alpha \notin [0,4]$, and $\alpha$ is conjugate to $z \mapsto \lambda^2 z$, where $\lambda$ satisfies $|\lambda| > 1$ and $\mathrm{tr}^2\alpha = (\lambda + \lambda^{-1})^2$. The transformation $\alpha$ acts on its axis by by moving each point hyperbolic distance $2\log|\lambda|$ toward the attractive fixed point.

The identity is not part of this classification.

A group $\Gamma$ is *elementary* if there is a single point on $\partial\mathbb{H}^3$, or a pair of points on $\partial\mathbb{H}^3$, or a single point in $\mathbb{H}^3$, which is invariant under all elements of $\Gamma$.

The generic group $\Gamma$ with two or more generators is nonelementary, and is likely to be nondiscrete as well. For example, any two loxodromic transformations $\alpha$ and $\beta$ without a common fixed point generate a nonelementary group $\Gamma = \langle\alpha,\beta\rangle$. The group $\Gamma$ is the homomorphic image, in many ways, of any surface group of genus $g \geq 2$.

The most important class of groups ruled out by the condition that $\Gamma$ be nonelementary are groups of rotations of the two-sphere and groups conjugate to them (unitary groups). We recall that a group, discrete or not, that is composed solely of elliptic transformations is conjugate to a group of rotations of the 2-sphere.

In anticipation of our later work, we also recall the definition of a two-generator *classical Schottky group* $G = \langle\alpha,\beta\rangle$. There are four mutually disjoint circles with mutually disjoint interiors, arranged as two pairs $(c_1, c_1')$ and $(c_2, c_2')$. The generator $\alpha$ sends the exterior of $c_1$ onto the interior of $c_1'$, and $\beta$ does the same for $(c_2, c_2')$. The common exterior of all four circles serves as a fundamental region for its action on its regular set $\Omega$.

Let $\pi : \Omega \to S := \Omega/G$ denote the natural projection. The surface $S$ has genus two, and $\pi(c_1)$ and $\pi(c_2)$ are disjoint, nondividing simple loops on $S$. If $d \subset S$ is a simple loop with an $\alpha$-invariant lift $d^* \subset \Omega$, the free homotopy class of $d$ in $S$ is uniquely determined up to Dehn twists about $\pi(c_1)$ (see §1.8).

The group $G$ extends to act on $\Omega \cup \mathbb{H}^3$; the quotient is a handlebody of genus two in which $\pi(c_1)$ and $\pi(c_2)$ are *compressing loops* that bound mutually disjoint *compressing disks* in the interior.

If, instead of circles, the pairs $(c_1, c_1')$ and $(c_2, c_2')$ are Jordan curves (which can always be assumed to be smooth), the resulting group is called more generally a (rank-two) *Schottky group*. According to [Chu], or [Z] in the handlebody interpretation, every set of free generators of a Schottky group (of the general kind!) corresponds to pairs of Jordan curves as described above.

Our method of construction in this paper will always yield classical Schottky groups in terms of designated generators. The extra knowledge that, for the designated generators, the loops can be taken as round circles is pleasing and convenient, but it is not really necessary for the proofs.

1.3. *Projective structures.* Let $R$ be a closed Riemann surface of genus at least two, and let $R = \mathbb{H}^2/G$ be its representation in the universal covering surface $\mathbb{H}^2$ (the two-dimensional hyperbolic plane) by a fuchsian covering group $G$. We will describe a projective structure first in the universal cover $\mathbb{H}^2$ and then intrinsically in $R$.

A *complex projective structure* with respect to $G$ is a meromorphic, locally univalent (i.e. locally injective) function $f : \mathbb{H}^2 \to f(\mathbb{H}^2) \subseteq \mathbb{S}^2$, for which there corresponds a homomorphism $\theta : G \to \Gamma \subset \mathrm{PSL}(2, \mathbb{C})$ such that $f(\gamma(t)) = \theta(\gamma)f(t)$ for any $t \in \mathbb{H}^2$ and any $\gamma \in G$. It follows that $f$ descends to a multivalued function $f_*$ on $R$, called the (*multivalued*) *developing map*; it "unrolls" $R$ onto the sphere. The *Schwarzian* derivative of $f$,

$$(1) \qquad S_t(f) := \left(\frac{f''}{f'}\right)' - \frac{1}{2}\left(\frac{f''}{f'}\right)^2 = \phi(t),$$

satisfies $\phi(\gamma(t))\gamma'^2(t) = \phi(t)$, and therefore descends to a holomorphic quadratic differential on $R$.

Conversely, given any holomorphic $\phi(t)$ in $\mathbb{H}^2$ with this invariance under $G$, there is a solution $f(t)$ of (1), uniquely determined up to post composition by Möbius transformations, which is a locally univalent meromorphic function that induces a homomorphism $\theta$ of $G$.

The Schwarzian equation is related to the second-order linear differential equation

$$(2) \qquad u''(t) + \tfrac{1}{2}\phi(t)u(t) = 0$$

as follows. The ratio $f(t) = u_1(t)/u_2(t)$ of any two linearly independent solutions $u_1$ and $u_2$ in $\mathbb{H}^2$ gives a solution $f$ of the Schwarzian; conversely, any solution $f$ of the Schwarzian can be so expressed, indeed

$$(3) \qquad u_2 = (f')^{-\frac{1}{2}}, \quad u_1 = fu_2,$$

if the Wronskian $\Delta(u_1, u_2)$, which is necessarily a constant, is normalized as $\Delta = 1$. Another pair $au_1 + bu_2$, $cu_1 + du_2$ of independent solutions corresponds to the solution $Bf$ of the Schwarzian, where $B(z) = (az + b)/(cz + d)$.

On the Riemann surface $R = \mathbb{H}^2/G$, a form of (2) that is invariant under change of local coordinates $z$ is,

$$(4) \qquad v''(z) + \tfrac{1}{2}\{\phi(z) + S_z(\pi^{-1})\}v(z) = 0,$$

where $\pi$ denotes the projection from $\mathbb{H}^2$. In interpreting this equation, the Schwarzian transforms as a connection under change of local coordinate $z \mapsto \zeta = \zeta(z)$ and $v$ transforms as a half-order differential (see [Ha-Sch]), specifically

$$v(\zeta(z))\zeta'(z)^{-\frac{1}{2}} = v(z).$$

The *monodromy group* and *monodromy representation* are computed as follows. Fix $\pi_1(R; O)$ with basepoint $O \in R$, and a solution $f_*(z)$ (or $v_1(z)/v_2(z)$) near $O$. Let $c \in \pi_1(R; O)$ be a simple loop based at $O$. Analytically continue $f_*$ (or $v_1/v_2$) around $c$, arriving back at a solution $\gamma f_*$ (or $\gamma(v_1/v_2)$), for $\gamma \in \mathrm{PSL}(2, \mathbb{C})$. Set $\theta(c) = \gamma$. In this manner the local solutions $f_*$ (or $v_1/v_2$) determine a monodromy epimorphism

$$\theta : \pi_1(R; O) \to \Gamma \subset \mathrm{PSL}(2, \mathbb{C}),$$

where $\Gamma$ is a monodromy group for the equation. A different local solution $Bf_*$ (or $B(v_1/v_2)$), coming possibly from a different choice of basepoint, determines a conjugate homomorphism $c \mapsto B\theta(c)B^{-1}$. Thus, the equation itself determines a conjugacy class of homomorphisms into $\mathrm{PSL}(2, \mathbb{C})$.

If $\mathcal{P}$ is a fundamental polygon for $G$ in $\mathbb{H}^2$, we can regard $f(\mathcal{P})$ as spread over the Riemann sphere, a *membrane* in Hejhal's terminology [He1]. The $\theta$-image of the edge pairing transformations of $\mathcal{P}$ will be edge-pairing transformations of the membrane $f(\mathcal{P})$, which therefore serves as an organizing principle for $\Gamma$.

From the *topological* point of view, a projective structure is defined by an orientation preserving local homeomorphism, called the (*multivalued*) *developing map*, of $R$ into $\mathbb{S}^2$ or, the (*single valued*) *developing map* of the universal cover $\tilde{R}$ into $\mathbb{S}^2$ which is equivariant with respect to the given homeomorphism $\theta$. From this perspective, the group $\Gamma$ is called the *holonomy* (or, more classically, *monodromy*) *group*. There is a unique complex structure on $R$ for which the local homeomorphism becomes conformal.

The fact that the Schwarzian equation can be replaced by the linear differential equation implies the following:

LEMMA 1.3.1. *If the homomorphism* $\theta : \pi_1(R; O) \to \Gamma \subset \mathrm{PSL}(2, \mathbb{C})$ *is induced by a projective structure on* $R$, *it can be lifted to a homomorphism* $\theta^* : \pi_1(R; O) \to \Gamma^* \subset \mathrm{SL}(2, \mathbb{C})$.

*Proof.* Consider an action of $G \cong \pi_1(R, O)$ on $\mathbb{H}^2$ given by the uniformization of the surface $R$, take an element $h \in G$. Then the solution pair (3) changes under analytic continuation from $t$ to $T = h(t)$ according to (see [Ha-Sch])

$$(5) \qquad \begin{pmatrix} 1/\sqrt{f'} \\ f/\sqrt{f'} \end{pmatrix}(T) = \sqrt{h'(t)} \begin{pmatrix} a & b \\ c & d \end{pmatrix} \begin{pmatrix} 1/\sqrt{f'} \\ f/\sqrt{f'} \end{pmatrix}(t),$$

where

$$\begin{pmatrix} a & b \\ c & d \end{pmatrix} \in \mathrm{SL}(2, \mathbb{C}), \quad \theta(h) = \frac{az+b}{cz+d}.$$

There are $2^{2g}$ possible choices for $\sqrt{h'(t)}$ over a set of canonical generators $\{h\}$ of $G$. After we make a choice we get the homomorphism

$$\theta^* : G \to \mathrm{SL}(2, \mathbb{C}), \quad \theta^*(h) = \begin{pmatrix} a & b \\ c & d \end{pmatrix} \in \mathrm{SL}(2, \mathbb{C}).$$

Note however that $\theta^*$ is not canonically determined by the differential equation (2).                                                                 $\square$

We emphasize that our notion of *lifting* does not require that the image $\Gamma$ of $\theta$ be isomorphic to the image of the lift $\theta^*$. For example, a lift to $\mathrm{SL}(2, \mathbb{C})$ of a half-rotation in $\mathrm{PSL}(2, \mathbb{C})$ has order four, not two.

We will refer to $\theta^*$ as a *linear monodromy representation* of the projective structure.

*Remark* 1.3.2. The projective structure associated with the equation $S_z(f) = \phi$ can be joined to the identity by means of solutions of $S_z(f) = t\phi$, for $t \in \mathbb{C}$.

1.4. *Branched projective structures.* A branched projective structure on a hyperbolic Riemann surface $R$ is a holomorphic mapping $f : \mathbb{H}^2 \to \mathbb{S}^2$ which is locally univalent except in a discrete subset of $\mathbb{H}^2$ and which is equivariant with respect to a homomorphism $\theta : G \to \mathrm{PSL}(2, \mathbb{C})$. We will say that such a structure is *singly branched* if $f'(z)$ has at most simple zeroes and the projection of the set $\{z : f'(z) = 0\}$ to $R$ is exactly one point $q$. These are the structures which appear in Theorem 1.1.1 and we will restrict our comments here to this special case. The more general case will be discussed separately in §11. Near such point $q$ (which we will identify with zero in local coordinates), the quadratic differential $\phi = S_z(f)$ has a Laurent expansion of the form,

$$(6) \qquad \phi(z) = \frac{-3}{2z^2} + \frac{b}{z} + \sum_{i=0}^{\infty} a_i z^i, \quad b^2 + 2a_0 = 0.$$

Conversely, if $\phi(z)$ has such an expression near $z = 0$, a solution of the Schwarzian will be of the form $f(z) = az^2(1 + o(1))$ near $z = 0$. With $\phi$ given by (6), the equation (2) has the two linearly independent solutions with expansions near $z = 0$ of the form

$$\begin{aligned} v_1(z) &= z^{3/2}(1 + o_1(1)), \\ v_2(z) &= z^{-1/2}(1 + o_2(1)). \end{aligned}$$

A circuit about $z = 0$ generates the monodromy

$$\begin{pmatrix} u_1 \\ u_2 \end{pmatrix} \mapsto J \begin{pmatrix} u_1 \\ u_2 \end{pmatrix}, \qquad \text{where} \quad J = \begin{pmatrix} -1 & 0 \\ 0 & -1 \end{pmatrix}.$$

The projectivized monodromy in $\text{PSL}(2, \mathbb{C})$ is just the identity.

Therefore the branched structure determines the homomorphism

$$\theta : \pi_1(R; O) \to \text{PSL}(2, \mathbb{C})$$

as in the unbranched case. However, $\theta$ cannot be lifted to a homomorphism into $\text{SL}(2, \mathbb{C})$. Indeed, given a standard presentation

$$\langle a_1, b_1, \ldots, a_g, b_g \mid \prod [b_i, a_i] = 1 \rangle$$

for $\pi_1(R; O)$, and matrix representations $A_i$ and $B_i$ for $\theta(a_i)$ and $\theta(b_i)$, we have

$$\theta^* \left( \prod [b_i, a_i] \right) = \prod [B_i, A_i] = J,$$

where $\theta^*(a_i) = A_i$ and $\theta^*(b_i) = B_i$.

We will discuss this matter further in §§11.5, 11.6.

1.5. *Parameter count.* The vector bundle $Q_g$ of quadratic differentials over Teichmüller space $\mathfrak{T}_g$ has complex dimension $6g - 6$. Likewise, the representation variety $V_g$ of homomorphisms $\theta : \pi_1(R; O) \to \text{PSL}(2, \mathbb{C})$, modulo conjugacy, has complex dimension $6g - 6$. Let $V_g' \subset V_g$ denote the subset of *nonelementary* representations, i.e. equivalence classes of homomorphisms whose images are nonelementary subgroups of $\text{PSL}(2, \mathbb{C})$. Theorem 1.1.1 asserts that the map $P_g$ of projective structures $Q_g \to V_g$ is surjective onto the component of $V_g'$ consisting of representations liftable to $\text{SL}(2, \mathbb{C})$. In fact, the image space $V_g'$ is itself a complex analytic manifold [Gu3], [He1]. According to [Go2], or as a consequence of Theorem 1.1.1, it has two components (one corresponds to liftable representations and the other one to unliftable representations). See [Ben-C-R] and [Li] for more information about representation varieties of surface groups.

According to Hejhal's holonomy theorem [He1] the map $P_g$ is a local homeomorphism which is shown in [E] to be locally biholomorphic. In particular, the set of points with a given monodromy $\theta$ is discrete. According to (1) in §1.6 below, there is at most one representative in the fiber over a particular Riemann surface. However $P_g$ is *not* a covering map [He1].

In Theorem 11.5.2 we will prove an analogue of Hejhal's holonomy theorem for singly branched projective structures; we prove that the holonomy mapping from the space of singly branched projective structures to $V_g$ is locally a fiber bundle with fiber of complex dimension 1.

1.6. *The global structure.* Recorded below are basic facts about projective structures. For the unbranched case, proofs are in [Gu1] and [Kra1, 2]. Other useful references are [Gu3] and [He1]; the latter includes extensive historical background.

Here is a brief proof that (in the unbranched case) the holonomy group $\Gamma = \theta(G)$ cannot be a unitary group, that is, cannot be conjugate to a group of isometries of $\mathbb{S}^2$. Assume otherwise. Then $\Gamma$ preserves the spherical metric $\rho$. Its pullback $f^*\rho$ is a $G$-invariant metric on $\mathbb{H}^2$ which is locally isometric to the sphere. Consequently $f^*\rho$ has constant curvature $+1$, in violation of the Gauss-Bonnet theorem.

For the case of a singly branched structures, property (1) below is a special case of [He1, Theorem 15], (2) will be established as Theorem 11.6.1, and (3) will be established as Corollary 11.6.1.

Below we consider projective structures $\sigma$ on $R = \mathbb{H}^2/G$ which have the holomorphic developing mapping $f : \mathbb{H}^2 \to \mathbb{S}^2$ and monodromy representation $\theta : G \to \mathrm{PSL}(2, \mathbb{C})$. Assume that $\sigma$ is either unbranched (i.e. $f$ is locally univalent) or is singly branched. Let $\theta(G) = \Gamma \subset \mathrm{PSL}(2, \mathbb{C})$ denote the monodromy group. The following three properties hold:

(1) If two developing mappings $f_1$ and $f_2$ determine the same homomorphism $\theta$, then $f_1 = f_2$.

(2) $\Gamma$ is a nonelementary group.

(3) The following statements are equivalent provided that, when $\sigma$ is branched, $f(\mathbb{H}^2)$ is not a round disk in $\mathbb{S}^2$:

   (i) $f(\mathbb{H}^2) \neq \mathbb{S}^2$;

   (ii) $\mathbb{H}^2 \to f(\mathbb{H}^2)$ is a possibly branched cover;

   (iii) $\Gamma$ acts discontinuously on $f(\mathbb{H}^2)$.

Property (1) does not rule out the possibility that the same target group $\Gamma$ may arise from different projective structures on $R$. Property (2) shows that the requirement in Theorem 1.1.1 that $\Gamma$ be nonelementary is necessary. The situation (3) has a rich structure as it is associated with the theory of covering surfaces; in particular it includes the theory of quasifuchsian groups and Schottky groups. In contrast, in the general case there is a bare minimum of structure because $\Gamma$ need not be discrete.

1.7. *Strategy of the proof.* Given a homomorphism

$$\theta : \pi_1(R; O) \to \Gamma \subset \mathrm{PSL}(2, \mathbb{C})$$

such that $\Gamma$ is nonelementary, the strategy consists of two parts.

*Part* A (§§3–5). Find a pants decomposition $\{P_i\}$ of $R$ with the property that $\theta(\pi_1(P_i))$, for $1 \le i \le 2g-2$, is a two-generator (classical) Schottky group.

We recall that a *pants* is a Riemann surface conformally equivalent to a three-holed sphere. A surface of genus $g \ge 2$ requires $3g-3$ simple loops to cut it into pants, and there results $2g-2$ pants. It has infinitely many homotopically distinct pants decompositions.

*Part* B (§§6–10). Find representations of the universal covers $\tilde{P}_i$ in the regular sets (i.e. domains of discontinuity) of $\theta(\pi_1(P_i)) \subset \mathbb{S}^2$. Glue them together as dictated by the combinatorics of $\{\tilde{P}_i\}$ in $\tilde{R}$, as relayed by $\theta$. In general there is a $\mathbb{Z}/2$-obstruction to such gluing. If there is no obstruction, we end up with a simply connected *pants configuration* $\tilde{S}$ over $\mathbb{S}^2$ that models the universal cover of a new Riemann surface $S$ homeomorphic to $R$. The projection of $\tilde{S}$ to $\mathbb{S}^2$ is a $\theta$-equivariant local homeomorphism. If there is an obstruction, introduce a single branch point of order 2 by applying a *twist*. This removes the obstruction to the construction and a new Riemann $S$ surface homeomorphic to $R$ can be assembled as before. The result is either unbranched or singly branched projective structure on $S$ with the monodromy representation $\theta$. According to Theorem 11.2.2 if $\theta$ lifts to $SL(2,\mathbb{C})$ then the structure has to be unbranched, if $\theta$ does not lift then the structure has to be singly branched; in other words, the $\mathbb{Z}/2$-obstruction to gluing is the 2nd Stiefel-Whitney class of the representation $\theta$. This proves Theorem 1.1.1.

The method used to assemble the pants configuration is a form of "grafting," first applied to kleinian groups in [Mas1].

1.8. *Terminology and notation.* Throughout this paper we will work on closed surface $R$, of genus $g \ge 2$. When convenient, we will assume that $R$ is a Riemann surface $R = \mathbb{H}^2/G$ in terms of its universal cover (which may be taken as the hyperbolic plane $\mathbb{H}^2$) and fuchsian cover group $G$. Fix $O \in R$ as the basepoint for its fundamental group $\pi_1(R; O)$. Let

$$\theta : \pi_1(R; O) \to \Gamma \subset PSL(2, \mathbb{C})$$

be the designated homomorphism with a nonelementary image $\Gamma$.

Throughout we will use lower case Latin letters $a, b, c, \ldots$ to denote elements of $\pi_1(R; O)$, and the corresponding Greek letters $\alpha, \beta, \gamma, \ldots$ to denote their $\theta$-images in $\Gamma$. A *nontrivial* loop is one not homotopic to a point.

We will write the compositions of both curves and transformations (and their associated matrices) starting at the right. Thus, $b$ follows $a$ in both $ba$ and $\theta(b)\theta(a) = \beta\alpha$.

By a *standard set of generators* $\{a_i, b_i\}$ of $\pi_1(R; O)$, where $1 \le i \le 2g$, we mean a set of oriented simple loops that generate the fundamental group and have the following properties. For each $i$, the loop $b_i$ crosses $a_i$ at $O$, from the right side of $a_i$ to the left, and is otherwise disjoint. For $j \neq i$, the simple loops $(a_j, b_j)$ are freely homotopic to simple loops disjoint from $(a_i, b_i)$. The product of the commutators

$$\prod_i b_i^{-1} a_i^{-1} b_i a_i$$

bounds a simply connected region lying to its left.

We will refer to a product $ba$ as a simple loop if it is homotopic to one (with fixed basepoint). Thus, for any $k \in \mathbb{Z}$, the loop $b_1 a_1^k$ is simple, and so are $a_2 b_1^{-1} a_1^k$ and $a_2 b_1 a_1$, but not $a_2^{-1} b_1 a_1^k$, or, for $k \neq 1$, the loop $a_2 b_1 a_1^k$. The curve $b_1^{-1} a_1 b_1$ is simple, but not $a_2 b_1^{-1} a_1 b_1$.

Often we will modify a simple loop $c \subset R$ by applying a *Dehn twist*, which can be described as follows. Let $A$ be an annular neighborhood about a (nontrivial) simple loop $a$. Orient $\partial A$ so that $A$ lies to its left. Hold one component of $\partial A$ fixed and rotate the other $|n|$-times in the positive or negative direction according to whether $n \geq 1$ or $n \leq -1$. This action extends to an orientation preserving homeomorphism $\delta^n$ of $A$, and then to all $R$, by setting $\delta^n = \mathrm{id}$ outside $A$. $\delta^n$, or more precisely its homotopy class on $R$, is called the Dehn twist of order $n$ about $a$. If $c$ is not freely homotopic to a curve disjoint from $a$, then $\delta^n(c)$ is not freely homotopic to $c$.

## 2. Fixed points of Möbius transformations

In this section we will collect the lemmas needed to control the type of composed transformations.

### 2.1. *Basic lemmas.*

LEMMA 2.1.1.

(i) *Suppose $\alpha$ is loxodromic and $\beta$ sends neither fixed point of $\alpha$ to the other. Given $M > 0$ there exists $N \geq 0$ such that $|\mathrm{tr}\, \beta \alpha^n| > M$ and $\beta \alpha^n$ is loxodromic for all $|n| \geq N$.*

(ii) *Suppose $\alpha$ is loxodromic and $\beta$ sends exactly one fixed point of $\alpha$ to the other. Given $M > 0$ there exists $N \geq 0$ such that $|\mathrm{tr}\, \beta \alpha^n| > M$ and $\beta \alpha^n$ is loxodromic for all $n \geq N$ (if $\beta$ sends repulsive to attractive) or for all $n \leq -N$ (if $\beta$ sends attractive to repulsive).*

(iii) *Suppose $\alpha$ is parabolic and $\beta$ does not share a fixed point with $\alpha$. Given $M > 0$ there exists $N \geq 0$ such that $|\mathrm{tr}\, \beta \alpha^n| > M$ and $\beta \alpha^n$ is loxodromic for all $|n| \geq N$.*

*Proof.* For (i) and (ii) we may assume

$$\alpha = \begin{pmatrix} \lambda & 0 \\ 0 & \lambda^{-1} \end{pmatrix} \quad \text{with } |\lambda| > 1, \qquad \beta = \begin{pmatrix} a & b \\ c & d \end{pmatrix} \quad \text{with } ad - bc = 1.$$

Then tr $\beta\alpha^n = \lambda^n a + \lambda^{-n} d$. Not both $a$ and $d$ can vanish, because $\beta$ does not interchange the fixed points of $\alpha$. The assertions now follow directly.

For (iii), we may assume

$$\alpha = \begin{pmatrix} 1 & 1 \\ 0 & 1 \end{pmatrix}, \qquad \beta = \begin{pmatrix} a & b \\ c & d \end{pmatrix} \quad \text{with } ad - bc = 1.$$

Then tr $\beta\alpha^n = (a+d) + nc$, where $c \neq 0$. Again, the desired conclusion follows. $\square$

LEMMA 2.1.2.    *Assume $\alpha$ is loxodromic with attractive fixed point $p^*$ and repulsive fixed point $p_*$.*

(i) *For any sequence $k \to +\infty$, the fixed points of $\beta\alpha^k$ converge to $\beta(p^*)$ and $p_*$. The fixed points of $\alpha^k\beta$ converge to $p^*$ and $\beta^{-1}(p_*)$.*

(ii) *For any sequence $k \to -\infty$, the fixed points of $\beta\alpha^k$ converge to $\beta(p_*)$ and $p^*$. The fixed points of $\alpha^k\beta$ converge to $p_*$ and $\beta^{-1}(p^*)$.*

*Proof.* Part (ii) follows from (i) upon replacing $\alpha$ by $\alpha^{-1}$. The computational proof is instructive. Set

$$\alpha = \begin{pmatrix} \lambda & 0 \\ 0 & \lambda^{-1} \end{pmatrix} \quad \text{and} \quad \beta = \begin{pmatrix} a & b \\ c & d \end{pmatrix},$$

where $|\lambda| \geq 1$ and $ad - bc = 1$. If $ac \neq 0$, the two fixed points of $\alpha^k\beta$ are $\lambda^{2k}a(1 \pm \sqrt{\Delta})/2c - d/2c$, where

$$\Delta = 1 + \frac{2d}{a\lambda^{2k}} + \frac{d^2}{a^2\lambda^{4k}} - \frac{4}{a^2\lambda^{2k}}.$$

The "+" fixed point approaches $\infty$ with $k$. The "−" fixed point has the form

$$\frac{2}{ac} - \frac{d}{c} - \frac{d^2}{2ac\lambda^{2k}}(1 + \sqrt{\Delta})^{-1} - \frac{d}{2c}.$$

This one approaches $-b/a = \beta^{-1}(0)$.

If $c = 0$, one fixed point of $\alpha^k\beta$ is $\infty$. The other one is $b/(d\lambda^{-2k} - a)$. This too approaches $-b/a = \beta^{-1}(0)$ with $k$.

If $a = 0$ the two fixed points are $(-d \pm \sqrt{d^2 - 4\lambda^{2k}})/2c$. Both approach $\infty$ with $k$. Here $\beta^{-1}(0) = \infty$.

The fixed points of $\beta\alpha^k = \beta(\alpha^k\beta)\beta$ converge to $\beta(p^*)$ and $\beta(\beta^{-1}p_*) = p_*$. $\square$

LEMMA 2.1.3. *Suppose $\gamma$ is loxodromic with attractive fixed point $p^*$, repulsive fixed point $p_*$.*

(i) *Suppose $\alpha(p^*) \neq p^*$ and $\beta(p_*) \neq p_*$. Given $M > 0$ there exists $N \geq 0$ such that $|\mathrm{tr}\,\gamma^{-n}\alpha\gamma^n\beta| > M$ and $\gamma^{-n}\alpha\gamma^n\beta$ is loxodromic and does not share a fixed point with $\alpha$ or $\beta$ for all $n \geq N$.*

(ii) *Suppose $\alpha(p_*) \neq p_*$ and $\beta(p^*) \neq p^*$. Given $M > 0$ there exists $N \geq 0$ such that $|\mathrm{tr}\,\gamma^{-n}\alpha\gamma^n\beta| > M$ and $\gamma^{-n}\alpha\gamma^n\beta$ is loxodromic and does not share a fixed point with $\alpha$ or $\beta$ for all $n \leq -N$.*

*Proof.* We may assume
$$\gamma = \begin{pmatrix} \lambda & 0 \\ 0 & \lambda^{-1} \end{pmatrix}, \qquad \alpha = \begin{pmatrix} a & b \\ c & d \end{pmatrix}, \qquad \beta = \begin{pmatrix} u & v \\ w & x \end{pmatrix},$$
with $|\lambda| > 1$, $ad - bc = 1$, and $ux - vw = 1$. We find that $\mathrm{tr}\,\gamma^{-n}\alpha\gamma^n\beta = \lambda^{2n}cv + \lambda^{-2n}bw + (au + dx)$.

In case (i) we have $c \neq 0$ (since $\alpha(p^*) \neq p^*$), and $v \neq 0$ (since $\beta(p_*) \neq p_*$); thus $cv \neq 0$ and $\gamma^{-n}\alpha\gamma^n\beta$ is loxodromic for all large $n$. Moreover, if $q$ is a fixed point of $\beta$, then $q \neq p_*$ but $\lim_{n \to +\infty} \gamma^{-n}\alpha\gamma^n\beta(q) = p_*$.

Suppose instead that $q$ is a fixed point of $\alpha$, and of $\gamma^{-n}\alpha\gamma^n\beta$ for all large $n$. First $q \neq p_*$ since $\gamma^{-n}\alpha\gamma^n\beta(p_*) = p_*$ implies $\beta(p_*) = p_*$. Then $\beta(q) \neq p_*$ for $\gamma^{-n}\alpha(p_*) = q$ holds for all large $n$ only if $q = p_*$ or $q = p^*$. Thus once again, $\lim_{n \to +\infty} \gamma^{-n}\alpha\gamma^n\beta(q) = p_* \neq q$.

In case (ii), $b \neq 0$ (since $\alpha(p_*) \neq p_*$), and $w \neq 0$ (since $\beta(p^*) \neq p^*$); hence $\gamma^{-n}\alpha\gamma^n\beta$ is loxodromic for all small $n$. Moreover if $q$ is a fixed point of $\beta$, we have $q \neq p^*$, but $\lim_{n \to +\infty} \gamma^{-n}\alpha\gamma^n\beta(q) = p_*$.

Suppose instead that $q$ is a fixed point of $\alpha$, and of $\gamma^{-n}\alpha\gamma^n\beta$ for all small $n$. Again $q \neq p^*$ and then $\beta(q) \neq p^*$. As above, $q$ cannot be a fixed point of $\gamma^{-n}\alpha\gamma^n\beta$ for all small $n$. $\qquad\square$

LEMMA 2.1.4. *Suppose $\alpha$ is a loxodromic transformation with fixed points $u$ and $v$.*

(i) *Given $p^* \neq u, v$ and $T > 2$, there exists $\varepsilon > 0$ such that if $\beta$ is any loxodromic transformation with fixed points $p, q$ satisfying $d(p, q) < \varepsilon$, $d(q, p^*) < \varepsilon$, and with trace satisfying $|\mathrm{tr}\,\beta| \geq T$, then $\alpha$ and $\beta$ generate a classical Schottky group.*

(ii) *Given $p, q \neq u, v$, there exists $T > 2$ such that if $\beta$ is any loxodromic transformation with fixed points $p, q$ and satisfying $|\mathrm{tr}\,\beta| \geq T$, and if $\alpha$ also satisfies $|\mathrm{tr}\,\alpha| \geq T$, then $\alpha$ and $\beta$ generate a classical Schottky group.*

*Proof.* A loxodromic transformation $\beta$ acts in $\mathbb{H}^3 \cup \partial \mathbb{H}^3$. If $P \subset \mathbb{H}^3$ is a plane orthogonal to its axis, so is $\beta(P)$. The two circles $\partial P$ and $\partial \beta(P)$ in $\partial \mathbb{H}^3$ that separate the fixed points $p$ and $q$ of $\beta$ bound what we will refer to as an *annular region $A$* for $\beta$. Given any point $q^* \neq p, q, u, v$ in $\partial \mathbb{H}^3$, there are annular regions for $\beta$ that contain $q^*$ in their interior.

Fix $p^* \subset \partial \mathbb{H}^3$ distinct from $q^*, p, q, u, v$. Let $(p_n, q_n)$ be a sequence with $p_n \neq q_n$ and $\lim p_n = \lim q_n = p^*$. Let $T_n$ be the transformation with fixed point $q^*$ such that $T_n(p) = p_n$, $T_n(q) = q_n$. Ultimately $T_n$ is loxodromic, its attractive fixed point converges to $p^*$, and $|\mathrm{tr}\, T_n| \to \infty$. Consider an annular domain $A$ for $\beta$ containing $q^*$ in its interior. $T_n A$ is an annular region for $T_n \beta T_n^{-1}$, all containing $q^*$. The sequence of bounding circles of $T_n A$ converge to the point $p^*$; that is, $T_n A$ converges to $\partial \mathbb{H}^3 \setminus \{p^*\}$. The analysis would be equally applicable to a family of transformations $\{\beta\}$, each with fixed points $p, q$, so long as they all satisfied $|\mathrm{tr}\, \beta| \geq T$ for some $T > 2$ (uniformly loxodromic).

Now let $A'$ be an annular domain for $\alpha$ containing $p^*$ in its interior. Ultimately the bounding circles of $T_n A$ also lie in the interior of $A'$. For such indices $n$, $\alpha$ and $T_n \beta T_n^{-1}$ generate a classical Schottky group. Part (i) follows at once.

To establish part (ii), note that both $\alpha$ and $\beta$ have annular domains whose boundaries are circles arbitrarily close to their fixed points, if $T$ is large enough. $\qquad \square$

COROLLARY 2.1.5. *Suppose $\gamma$ is loxodromic with attractive fixed point $p^*$, repulsive fixed point $p_*$, and $\alpha, \beta$ are loxodromic as well.*

(i) *If $\alpha(p^*) \neq p^*$ and $\beta(p_*) \neq p_*$ there exists $N \geq 0$ such that $\gamma^{-n} \alpha \gamma^n$ and $\beta$ generate a classical Schottky group for all $n \geq N$.*

(ii) *If $\alpha(p_*) \neq p_*$ and $\beta(p^*) \neq p^*$ there exists $N \geq 0$ such that $\gamma^{-n} \alpha \gamma^n$ and $\beta$ generate a classical Schottky group for all $n \leq -N$.*

*Proof.* This is a corollary also of Lemma 2.1.3. In case (i), the fixed points of $\gamma^{-n} \alpha \gamma^n$ are arbitrarily close to $p_*$ for large $n$, since $p^*$ is not fixed by $\alpha$, where $p_*$ is not fixed by $\beta$. In case (ii), the fixed points of $\gamma^{-n} \alpha \gamma^n$ are arbitrarily close to $p^*$, for small $n$. $\qquad \square$

2.2. *Lemmas regarding half-rotations.*

LEMMA 2.2.1. *Suppose $\alpha$ and $\beta$ each have two fixed points and $\beta$ sends one of the fixed points of $\alpha$ to the other. Then $\alpha$ likewise sends one of the fixed points of $\beta$ to the other if and only if*

$$\mathrm{tr}^2 \alpha = \mathrm{tr}^2 \beta.$$

*Proof.* We may assume that

$$\alpha = \begin{pmatrix} \lambda & 0 \\ 0 & \lambda^{-1} \end{pmatrix} \quad \text{and} \quad \beta = \begin{pmatrix} 0 & b \\ c & d \end{pmatrix},$$

where $\lambda \neq \pm 1$ and $bc = -1$. The fixed points of $\beta$ are $(-d \pm \sqrt{d^2 - 4})/2c$. Suppose $\alpha$ sends one to the other. Each case implies and is implied by one of the relations

$$d(\lambda^2 - 1) = \pm \sqrt{d^2 - 4}(\lambda^2 + 1).$$

Squaring, we get $d^2\lambda^2 = (\lambda^2 + 1)^2$, or

$$\operatorname{tr}\beta = d = \pm(\lambda + \lambda^{-1}) = \pm\operatorname{tr}\alpha. \qquad \square$$

LEMMA 2.2.2. *An element $J$ of order two interchanges the fixed points of an elliptic or loxodromic transformation $\gamma$ if and only if*

$$J\gamma J = \gamma^{-1},$$

*and fixes them if and only if*

$$J\gamma J = \gamma.$$

*It fixes the fixed point of a parabolic transformation $\gamma$ if and only if*

$$J\gamma J = \gamma^{-1}.$$

*Proof.* For the first part we may assume that

$$\gamma = \begin{pmatrix} \lambda & 0 \\ 0 & \lambda^{-1} \end{pmatrix} \quad \text{and} \quad J = \begin{pmatrix} 0 & b \\ -b^{-1} & 0 \end{pmatrix},$$

while for the second,

$$\gamma = \begin{pmatrix} 1 & b \\ 0 & 1 \end{pmatrix} \quad \text{and} \quad J = \begin{pmatrix} i & 0 \\ 0 & -i \end{pmatrix}.$$

The conclusion is verified by computation. $\qquad \square$

LEMMA 2.2.3. *Suppose $\alpha$ and $\beta$ are loxodromic without both fixed points in common. $J$ is an element of order two.*

(i) *If $J$ interchanges the fixed points of both $\alpha$ and $\beta$, $J$ neither interchanges nor fixes the fixed points of $\beta\alpha$.*

(ii) *If $J$ interchanges the fixed points of $\beta$ but not of $\beta\alpha^k$ for some $k \neq 0$, then $J\beta$ does not interchange the fixed points of $\alpha$.*

(iii) *If $J$ interchanges the fixed points of both $\beta\alpha^k$ and $\beta\alpha^{k+1}$ for some $k$, then $J$ interchanges the fixed points of $\beta\alpha^k$ for all $k$, but neither interchanges nor fixes the fixed points of $\alpha$, and does not interchange the fixed points of $\alpha^m\beta$ for $m \neq 0$.*

*Proof.* For (i), $J\beta\alpha J = \beta^{-1}\alpha^{-1} \neq \alpha^{-1}\beta^{-1}$, $\beta\alpha$.

For (ii), $J_1 = J\beta$ has order two, $J_1 \neq$ id. If $J_1\alpha^k J_1 = \alpha^{-k}$, then $J\beta\alpha^k J = \alpha^{-k}\beta^{-1}$, a contradiction.

For (iii), the hypotheses $J\beta\alpha^k J = \alpha^{-k}\beta^{-1}$ and $J\beta\alpha^{k+1} J = \alpha^{-k}\beta^{-1}J\alpha J$ imply in turn that

$$\alpha^k\beta^{-1}J\alpha J = \alpha^{-k-1}\beta^{-1},$$

or

$$J\alpha J = \beta\alpha^{-1}\beta^{-1}, \quad (\neq \alpha^{-1}, \alpha).$$

Then

$$\alpha^{-k}\beta^{-1} = J\beta\alpha^k J = J\beta J\beta\alpha^{-k}\beta^{-1},$$

or $\quad J\beta J = \beta^{-1}$. Now, for any $k$,

$$J\beta\alpha^k J = \beta^{-1}\beta\alpha^{-k}\beta^{-1} = \alpha^{-k}\beta^{-1}.$$

Finally, for any $m \neq 0$,

$$J\alpha^m\beta J = \beta\alpha^{-m}\beta\beta^{-1} = \beta\alpha^{-m}\beta^{-2} \neq \beta^{-1}\alpha^{-m}.$$

(Note the proof holds as well if some $\beta\alpha^k$ is parabolic, under appropriate interpretation; see Lemma 2.2.2.) □

LEMMA 2.2.4. *Suppose $\alpha$ has two fixed points but $\alpha^2 \neq$ id, while $J$ is an element of order two that does not interchange the fixed points of $\alpha$. Then $(\alpha J)^2 \neq$ id and $(J\alpha)^2 \neq$ id.*

*Proof.* We may assume that

$$\alpha = \begin{pmatrix} \lambda & 0 \\ 0 & \lambda^{-1} \end{pmatrix} \quad \text{and} \quad J = \begin{pmatrix} a & b \\ c & -a \end{pmatrix},$$

with $\lambda^2 \neq \pm 1$, $a^2 + bc = -1$. Then

$$(\alpha J)^2 = \begin{pmatrix} \lambda^2 a^2 + bc & \lambda^2 ab - ab \\ ac - \lambda^{-2}ac & bc + \lambda^{-2}a^2 \end{pmatrix}.$$

If $(\alpha J)^2 = $ id, then

$$ab(\lambda^2 - 1) = 0,$$
$$ac(1 - \lambda^{-2}) = 0.$$

Either $a = 0$ or $b = c = 0$. The former case is impossible by hypothesis. If $b = c = 0$, then since $a^2 = -1$, we get $\lambda^2 = \lambda^{-2} = 1$. This is again a contradiction. □

LEMMA 2.2.5. *Suppose both $J$ and $J_1J$ interchange the fixed points of the loxodromic or elliptic transformation $\gamma$. Then $J_1$ fixes the fixed points of $\gamma$.*

*Proof.* Under the hypothesis, if $p$, $q$ denote the fixed points of $\gamma$, we have $J(p) = J_1J(p)$ and $J(q) = J_1J(q)$. Hence $J(p) = q$ and $J(q) = p$ are fixed by $J_1$. □

*Remark* 2.2.6. Suppose $\alpha$ and $\beta$ are loxodromic without a common fixed point and $\beta$ does not send one fixed point of $\alpha$ to the other. If $\gamma\beta$ fixes or interchanges the fixed points of $\alpha$, then $\gamma\beta^{-1}$ has neither of these properties. In the latter case, $\gamma\alpha\beta\alpha^{-1}$ does not send one fixed point of $\alpha$ to the other. What will prevent us from making use of such facts as these is that if $\gamma\beta$, for example, is the $\theta$-image of a simple loop, then in general $\gamma\beta^{-1}$ and $\gamma\alpha\beta\alpha^{-1}$ are not.

## A. The Pants decomposition

## 3. Finding a handle

3.1. *Handles.* By a *handle $H = \langle a, b \rangle$* we mean two simple loops $a, b \in \pi_1(R; O)$, crossing at $O$ but otherwise disjoint, and such that $\alpha = \theta(a)$ and $\beta = \theta(b)$ are loxodromic and generate a nonelementary subgroup $\langle \alpha, \beta \rangle$ of $\Gamma$.

PROPOSITION 3.1.1. *There exists a handle in $R$.*

*Proof.* The proof will occupy the remainder of this chapter.

3.2. *Case* 1. There exists a simple, nondividing loop $a \in \pi_1(R; O)$ for which $\theta(a) = \alpha$ is loxodromic. Choose $b \in \pi_1(R; O)$ such that $b$ is a simple loop crossing $a$ exactly at $O$, and set $\beta = \theta(b)$.

Suppose first that $\beta$ neither interchanges the fixed points of $\alpha$ nor shares a fixed point with $\alpha$. Then, by Lemma 2.1.1, $\beta\alpha^k$ is loxodromic for some $k$. Moreover, $\langle \alpha, \beta\alpha^k \rangle$ is nonelementary. We can consequently choose $H = \langle a, ba^k \rangle$.

Next suppose that $\beta$ shares exactly one fixed point $p$ with $\alpha$. Because $\Gamma$ is not elementary, there is a simple loop $y \in \pi_1(R; O)$ that does not cross $a$ or $b$ and such that $\theta(y) = \eta$ does not fix $p$. Take $y$ with the orientation such that $ay$ is homotopic to a simple loop. For any $k$, the loop $ay$ is homotopic to a simple loop that crosses $ba^k$ exactly at $O$ (Figure 1).

Now $\alpha\eta$ does not share the fixed point $p$ of $\beta\alpha^k$. For at most one value of $k$, $\alpha\eta$ shares another fixed point $q$ of $\beta\alpha^k$. For if

$$\alpha\eta(q) = q = \beta\alpha^k(q) = \beta\alpha^{k+m}(q)$$

with $m \neq 0$, we have $\alpha(q) = q$, and then $\eta(q) = q = \beta(q)$, a contradiction since $q \neq p$. Nor can $\alpha\eta$ send the fixed point $p$ of $\beta\alpha^k$ to another fixed point $q = \alpha\eta(p)$ of $\beta\alpha^k$ for more than one $k$. For

$$\alpha\eta(p) = q = \beta\alpha^k(q) = \beta\alpha^{k+m}(q)$$

with $m \neq 0$ implies that $\alpha(q) = q$, and then $\beta(q) = q$. This is impossible since $q \neq p$. Thus there exists $k$ such that $\alpha\eta$ neither interchanges the fixed points of $\beta\alpha^k$, nor fixes any. By Lemma 2.1.1, we may also assume that $\beta\alpha^k$ is loxodromic.

Consequently we can return to the case above with $ba^k$ and $ay$.

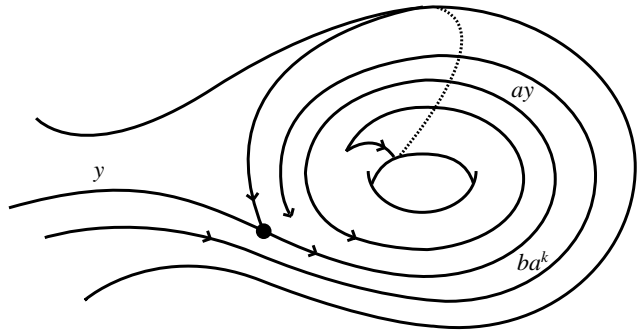

Figure 1.

Finally, suppose that $\beta$ either fixes both fixed points of $\alpha$ or interchanges them. Again find a simple loop $y$ that does not cross $a$ or $b$ and such that $\eta = \theta(y)$ neither fixes both fixed points of $\alpha$ nor interchanges them. Take the orientation of $y$ so that $yb$ is homotopic to a simple loop. Then $\eta\beta$ neither fixes both fixed points of $\alpha$ nor interchanges them. Consequently we can return to one of the cases above with $a$ and $yb$.

3.3.  *Case* 2. There is a simple, nondividing loop $a \in \pi_1(R; O)$ such that $\theta(a) = \alpha$ is parabolic. Let $b \in \pi_1(R; O)$ be a simple loop that crosses $a$ exactly at $O$.

If $\beta = \theta(b)$ does not fix the fixed point $p$ of $\alpha$, then $\beta\alpha^k$ is loxodromic for all large $|k|$, by Lemma 2.1.1. Thus we are back to Case 1.

Suppose instead that $\beta(p) = p$. There is a simple loop $y \in \pi_1(R; O)$, not crossing $a$ or $b$, and such that $\eta = \theta(y)$ does not fix $p$. We may take $y$ with the orientation for which $yb$ is homotopic to a simple loop, and hence also $yba^k$ is homotopic to a simple loop. Since $\eta\beta(p) \neq p$, we conclude that $\eta\beta\alpha^k$ is loxodromic for some $k$, and $yba^k$ brings us, once again, back to Case 1.

3.4. *Case 3.* Let $\{a_i, b_i\}$ be a canonical basis for $\pi_1(R; O)$, with $\theta(a_i) = \alpha_i$ and $\theta(b_i) = \beta_i$. Assume that all the elements $\alpha_i$, $\beta_i$, $\alpha_j \alpha_i$, $\beta_j \beta_i$, and $\beta_j \alpha_i$ are elliptic or the identity. As the basis of our analysis of this case, we will find a simple dividing loop $d$ for which $\theta(d)$ is loxodromic.

In this section we will establish some useful lemmas.

LEMMA 3.4.1. *If $\alpha$ and $\beta$ are elliptic, and their axes are not coplanar in $\mathbb{H}^3$, then $\beta \alpha$ is loxodromic.*

*Proof.* Let $P$ denote the plane in $\mathbb{H}^3$ spanned by the axis of $\alpha$ and the common perpendicular $l$ to that and the axis of $\beta$. Form the "open book" for $P$ with spine along the axis of $\alpha$ and angle half the rotation angle of $\alpha$. Then $\alpha = R_l R_{l_\alpha}$, where $R_{l_\alpha}$ and $R_l$ are half-rotations (180°) about the lines orthogonal to the axis of $\alpha$ indicated in Figure 2. Similarly, $\beta = R_{l_\beta} R_l$, where $l_\beta$ is the line orthogonal to the axis of $\beta$ at its intersection with $l$, and lies halfway between $l$ and $\beta(l)$.

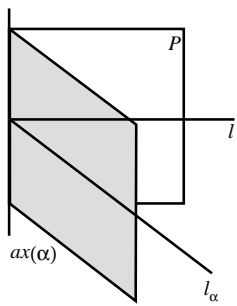

Figure 2. Open book for plane $P$

Consequently, $\beta \alpha = R_{l_\beta} R_{l_\alpha}$. Therefore $\beta \alpha$ is elliptic if and only if the lines $l_\alpha$ and $l_\beta$ intersect in $\mathbb{H}^3$: if instead they meet at $\partial \mathbb{H}^3$, then the composition $\beta \alpha$ is parabolic, and if they do not meet at all in $\mathbb{H}^3 \cup \partial \mathbb{H}^3$, then the composition is loxodromic. Since the axis of $\beta$ does not lie in $P$, $l_\beta$ does not lie in the plane spanned by $l_\alpha$ and $l$. Therefore $l_\alpha$ and $l_\beta$ cannot meet anywhere. $\square$

COROLLARY 3.4.2. *Under the hypotheses of Case 3, the axes of the non-identity elements of $\{\alpha_i, \beta_j\}$ either:*

(a) *all pass through some point $\zeta \in \mathbb{H}^3$, or*

(b) *all lie in a plane $P \subset \mathbb{H}^3$, or*

(c) *are all orthogonal to a plane $P \subset \mathbb{H}^3$.*

*Proof.* Apply Lemma 3.4.1 to the set $\{\alpha_i, \beta_j\}$.                    □

Note that, in case (c), the plane $P$ contains all the lines $l_{\alpha_i}$ and $l_{\beta_i}$. This is the fuchsian case: all elements of $\Gamma$ preserve $P$.

Case (a) does not arise for our situation since $\Gamma$ is nonelementary.

LEMMA 3.4.3. *Suppose $\alpha$ and $\beta$ are elliptic with distinct axes that lie in a plane $P \subset \mathbb{H}^3$. Assume $\beta\alpha$ is also elliptic. Its axis cannot lie in $P$.*

*Proof.* The axes of $\beta\alpha$ and $\alpha$ are different, so there is a fixed point $x$ of $\beta\alpha$ not lying in the axis of $\alpha$. Set $y = \alpha(x)$; then $\beta(y) = x$. Let the plane $P'$ be the perpendicular bisector of the line segment $[x, y]$. By construction, $x$ and $y$ are equidistant from $P'$. But $x$ and $y$ are also equidistant from the axis of $\alpha$, since $\alpha$ is a rotation about its axis. All points equidistant from $x$ and $y$ lie in $P'$, so the axis of $\alpha$ is contained in $P'$. Since $x$ and $y$ are also equidistant from the axis of $\beta$, this line, too, is contained in $P'$. We conclude that $P' = P$, so $x \notin P$.                    □

In fact, the proof shows that if the axis of $\beta\alpha$ meets $P$ or $\partial P$, it does so at, and only at, a point of intersection or common endpoint of the axes of $\alpha$ and $\beta$.

LEMMA 3.4.4. *Suppose $\alpha$, $\beta$, and $\gamma = \beta\alpha$ are elliptic with distinct axes, and that they preserve a plane $P \subset \mathbb{H}^3$. Then $\beta^{-1}\alpha^{-1}\beta\alpha$ is loxodromic.*

*Proof.* Let $a, b, c$ denote the fixed points in $P$ of $\alpha$, $\beta$ and $\gamma$. Replace $\alpha$ and $\beta$ by the inverses, if necessary, so that they rotate counterclockwise about $a$ and $b$. Let $R_1 = J$, $R_2$ and $R_3$ denote the reflection in the lines through $[a, b]$, $[b, c]$ and $[c, a]$, respectively. Then $\alpha = R_1 R_3$, $\beta = R_2 R_1$, and $\gamma = R_2 R_3$. The vertex angles of the triangle in Figure 3 represent the half-rotation angles. Then

$$\beta^{-1}\alpha^{-1}\beta\alpha = J\gamma J\gamma.$$

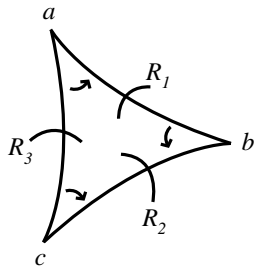

Figure 3. Reflection triangle for $\alpha$, $\beta$, $\beta * \alpha$

In order to better study $J\gamma J\gamma$, we take the line $l$ through $a$ and $b$ to be the real diameter in the disk model of $P$ (Figure 4). $J$ is reflection in $l$; let $R$ denote reflection in the vertical line through $c$ and $Jc$. Let $\theta$ denote the half-rotation angle of $\gamma$. Let $l_1$ denote the line through $c$ subtending angle $\theta$ with the vertical, and set $l_2 = Rl_1$. Let $R_1'$ denote reflection in $l_1$ and $R_2'$ reflection in $l_2$.

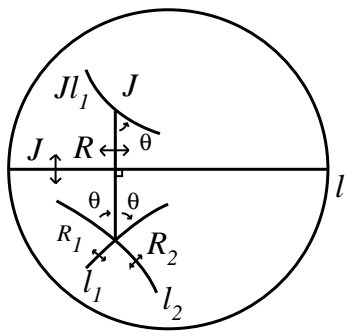

Figure 4. Reflection in $Jl_1$ and $l_2$

Now we have $\gamma = R_1'R = RR_2'$ and

$$J\gamma J = JR_1'JJRJ = R_1'^*R,$$

where $R_1'^* = JR_1'J$ is reflection in the line $Jl_1$. Consequently,

$$J\gamma J\gamma = R_1'^*RRR_2' = R_1'^*R_2'.$$

The lines $Jl_1$ and $l_2$ cannot intersect in $P \cup \partial P$. Therefore the composition of reflections in them, $R_1'^*R_2'$, is loxodromic (hyperbolic). □

Note, however, that $R_1'^*R_1' = J\gamma J\gamma^{-1} = (\beta\alpha^2\beta)^{-1}$ can sometimes be elliptic.

3.5. *Case* 3 (*continued*). Suppose that the elements $\{\alpha_i, \beta_i\}$, which are all elliptic or the identity, preserve a plane $P \subset \mathbb{H}^3$ (Case (c) of Corollary 3.4.2). We may assume that $\alpha_1 \neq \mathrm{id}$.

Consider first the case that $\beta_1$ is elliptic and its fixed point in $P$ differs from that of $\alpha_1$. Then the transformation $\delta = \beta_1^{-1}\alpha_1^{-1}\beta_1\alpha_1$, which corresponds to the simple loop $d = b_1^{-1}a_1^{-1}b_1a_1$, is hyperbolic (Lemma 3.4.4). Because $d$ divides $R$, there exists an element $c$ of $\{a_2, b_2, \ldots, a_g, b_g\}$ with $\gamma = \theta(c) \neq \mathrm{id}$. Apply the Dehn twist of order $n$ about $d$ to the simple loop $cb_1$, to get $cd^nb_1d^{-n}$. Its image $\gamma\delta^n\beta_1\delta^{-n}$ is loxodromic for all large $|n|$ by Lemma 2.1.3, since the fixed points on $\partial P$ of the hyperbolic $\delta$ are necessarily different from those of the elliptics $\gamma$ and $\beta_1$ in $P$. Since $cd^nb_1d^{-n}$ is a simple, nondividing loop, we can return with it to Case 1 (§3.2).

Consider next the case where $\beta_1$ has the same fixed point in $P$ as $\alpha_1$, or is the identity. We can find $c$ in $\{a_2, b_2, \ldots, a_g, b_g\}$ such that $\gamma = \theta(c)$ does not have the same fixed point in $P$ as $\alpha_1$. If $\theta(cb_1a_1)$ is not elliptic, return with $cb_1a_1$ to Case 1 or 2. Otherwise, set $d = (cb_1)^{-1}a_1^{-1}(cb_1)a_1$, and apply the Dehn twist about $d$ to $ca_1$ for a sufficiently high power. As above, return to Case 1 with the result.

Next, suppose that the axes of the elliptic elements $\{\alpha_i, \beta_i\}$, which are all elliptic or the identity, lie in a plane $P \subset \mathbb{H}^3$ (Case (b) of Corollary 3.4.2). We may assume that $\alpha_1 \neq \mathrm{id}$.

Assume first that the axes of $\alpha_1$ and $\beta_1$ differ. If they cross at $p \in P$, or meet at $p \in \partial P$, there is an element $c$ of $\{a_2, b_2, \ldots, a_g, b_g\}$ such that the axis of $\gamma = \theta(c)$ does not contain $p$. By Lemma 3.4.3, the axis of $\beta_1\alpha_1$ does not lie in $P$, but it crosses $P$ at $p$ or meets $\partial P$ at $p$. Consequently this axis is not coplanar with the axis of $\gamma$, which lies in $P$. Now Lemma 3.4.1 says that $\gamma\beta_1\alpha_1$ is loxodromic. Return to Case 1 with $cb_1a_1$.

On the other hand, suppose that the axis $l_1$ of $\alpha_1$ and the axis $l_2$ of $\beta_1$ are disjoint in $P \cup \partial P$. If the axis $l$ of $\beta_1\alpha_1$ is not coplanar with $l_1$ or $l_2$, the situation is again as above. If $l$ is coplanar with each of $l_1$ and $l_2$, it cannot meet $P \cup \partial P$. The plane $P'$ orthogonal to $l$ and to $P$ is necessarily orthogonal to $l_1$ and $l_2$. If the axes of all nonidentity elements of $\{\alpha_1, \beta_1, \ldots, \}$ are orthogonal to $P'$, we can return to the first subcase of this section. Otherwise the axis of some $\gamma \in \{\alpha_2, \beta_2, \ldots\}$ is not orthogonal to $P'$. Then $\gamma\beta_1\alpha_1$ is loxodromic, since the axis of $\gamma$ is contained in $P$.

Finally we need to consider the situation where $\beta_1 = \mathrm{id}$ or the axes of $\alpha_1$ and $\beta_1$ coincide. Find $\delta$ in $\{\alpha_2, \beta_2, \ldots\}$ distinct from the identity and having an axis distinct from that of $\alpha_1$. Replace $\beta_1$ by $\delta$ in the analysis above. The triple of loops in $\pi_1(R; O)$ giving rise to the loxodromic element found there also corresponds to a simple nondividing loop, and it is only this property that is needed.

In light of Corollary 3.4.2, the analysis of Case 3 is complete. A handle exists, and Proposition 3.1.1 is proved.

## 4. Cutting the handles

4.1. We have found a special handle $H$ as specified in §3. The next step is to cut all the other (topological) handles, ending up with a (connected) surface of genus one with $2(g - 1)$ boundary components. In cutting the handles, we will require that the $\theta$-image of each cutting loop is loxodromic.

Although $H$, or rather the established properties of the $\theta$-image of its fundamental group, serves to anchor the cutting process, in fact $H$ itself will have to undergo successive changes. It will become more and more complicated in terms of an initial basis of $\pi_1(R; O)$. Roughly speaking, we will be applying Dehn twists of possibly high order to felicitous combinations of simple loops. The process will be governed by the applicability of the lemmas of §2 to yield loxodromic transformations, yet still arising under $\theta$ from simple loops in $\pi_1(R; O)$.

4.2. Let $H = \langle a, b \rangle$ denote the special handle found in Chapter 3, and set $\alpha = \theta(a)$, $\beta = \theta(b)$. We claim that after replacing $H = \langle a, b \rangle$ by another handle of the form $\langle ab^q, b \rangle$ or $\langle b, ab^q \rangle$, if necessary, we can assume that $\beta$ does not send one fixed point of $\alpha$ to the other.

For suppose $\beta$ sends one fixed point of $\alpha$ to the other. Then, using Lemma 2.1.1(ii), find $q$ so that $\alpha\beta^q$ is loxodromic and $\mathrm{tr}^2\alpha\beta^q \neq \mathrm{tr}^2\beta$. Necessarily, $\alpha\beta^q$ does not share either of its fixed points with $\beta$. By Lemma 2.2.1, at least one of the following statements is true: $\beta$ does not send one fixed point of $\alpha\beta^q$ to the other, or $\alpha\beta^q$ does not send one fixed point of $\beta$ to the other.

4.3. Now suppose that $\langle x, y \rangle$ is another pair of loops in $\pi_1(R; O)$ of the form $x = u^{-1}x'u$, $y = u^{-1}y'u$, where $x'$ and $y'$ are simple loops disjoint from $a$ and $b$, with one intersection point where they cross, and $u$ is a simple arc from $a \cap b = O$ to $x' \cap y'$, otherwise disjoint from $a, b, x', y'$ (see Figure 5).

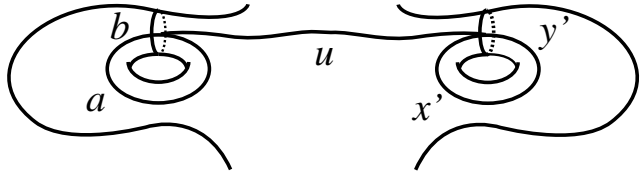

Figure 5. Connection to handle $H$

Consider $d = yba^k$ and its $\theta$-image $\delta = \eta\beta\alpha^k$, for some $k$. Set $\xi = \theta(x)$ and $\eta = \theta(y)$. The effect of a Dehn twist of order $n$ about $d$ is

$$\begin{aligned}
\langle \alpha, \beta\alpha^k \rangle &\mapsto \langle \delta^n\alpha, \beta\alpha^k \rangle, \\
\langle \xi, \eta \rangle &\mapsto \langle \delta^n\xi, \eta \rangle.
\end{aligned}$$

We claim that there exist $k$ and $n$ such that:

(i) $\beta\alpha^k$ is loxodromic;

(ii) $\delta = \eta\beta\alpha^k$ is loxodromic;

(iii) $\delta^n\alpha$ is loxodromic without a common fixed point with $\beta\alpha^k$;

(iv) $\delta^n\xi$ is loxodromic;

or that, after necessary relabeling and rearrangement to be spelled out below, analogous properties hold. This claim will be established in the four steps of §4.4.

Once this is accomplished, we will replace the handle $H = \langle a, b \rangle$ by the handle $\langle d^n a, ba^k \rangle$, and cut $R$ along $d^n x$, represented by a freely homotopic simple loop. This operation will also serve as the basis of an induction procedure.

Note that it may well be that $\eta = \mathrm{id}$, or $\xi = \mathrm{id}$, or both. In the former case, property (ii) is satisfied with (i), and in the second, property (iv) is satisfied with (ii).

4.4.  *Step* (i). By §4.2 and Lemma 2.1.1(i), there exists $K \geq 0$ such that $\beta \alpha^k$ is loxodromic for all $|k| \geq K$.

*Step* (ii). If $\eta \beta$ does not interchange the fixed points of $\alpha$, then by Lemma 2.1.1 we may take $K$ in step (i) so large that $\delta = \eta \beta \alpha^k$ is loxodromic for all $k \geq K$ or for all $k \leq -K$.

If, however, $\eta \beta$ does interchange the fixed points of $\alpha$ but $\xi \beta$ does not, interchange $\eta$ and $\xi$ and return to the paragraph above.

Finally, if both $\eta \beta = J$ and $\xi \beta = J_1$ interchange the fixed points of $\alpha$, then $\eta \xi^{-1} = JJ_1$ fixes them (and is either loxodromic or the identity). In this case replace $\langle x, y \rangle$ by $\langle x, yx^{-1} \rangle$, and $\eta$ by $\eta \xi^{-1}$, and revert to the original notation. For this case, then, $\eta \beta \alpha^k$ is loxodromic for all $|k| \geq K$, for some $K$.

*Step* (iii).  First note that for no $k \in \mathbb{Z}$ and no $m \neq 0$ can both $\delta^n \alpha$ and $\delta^{n+m} \alpha$ have fixed points in common with $\beta \alpha^k$. For the relations

$$\delta^n \alpha(p) = p = \beta \alpha^k(p),$$
$$\delta^{n+m} \alpha(p) = p = \delta^n \delta^m \alpha(p)$$

imply that $\alpha(p)$ is a fixed point of $\delta$, then that $p$ is a fixed point of $\alpha$, and finally that $p$ is a fixed point of $\beta$. The last consequence is impossible.

For any sequence $k \to +\infty$, according to Lemma 2.1.2 the fixed points of $\delta = \eta \beta \alpha^k$ converge to $\eta \beta(q)$ and $p$, where $q$ and $p$ denote the attractive and repulsive fixed points of $\alpha$, respectively. Thus, if $\alpha$ sends one fixed point of $\delta$ to the other for this sequence, then $\eta \beta(q) = p$. Similarly, for a sequence $k \to -\infty$, the fixed points of $\delta$ converge to $\eta \beta(p)$ and $q$. If, for this sequence, $\alpha$ sends one fixed point of $\delta$ to the other, then $\eta \beta(p) = q$. By our construction, $\eta \beta$ does not interchange the fixed points of $\alpha$, so $\alpha$ cannot interchange the fixed points of $\delta = \eta \beta \alpha^k$ for both a sequence $k \to +\infty$ and a sequence $k \to -\infty$.

Now if $\eta \beta(p) = q$, so that $\delta$ is loxodromic for all $k \geq K$ (step (ii)), then for sufficiently large $K$, $\alpha$ cannot send one fixed point of $\delta = \eta \beta \alpha^k$ to the other. Likewise, if $\eta \beta(q) = p$ so that $\delta$ is loxodromic for $k \leq -K$, again $\alpha$ cannot send one fixed point of $\delta$ to the other, for sufficiently large $K$. If $\eta\beta$ sends neither fixed point of $\alpha$ to the other, $\alpha$ sends neither fixed point of $\delta$ to the other, for all large $|k|$.

We conclude that there exists $K \geq 0$ such that $\delta = \eta\beta\alpha^k$ is loxodromic for any $k \geq K$ or any $k \leq -K$, or both. Furthermore, $\alpha$ does not send one fixed point of $\delta$ to the other. Given $k$ in the admissible range, there exists $N = N(k) \geq 0$ such that $\delta^n\alpha$, for all $|n| \geq N$, is loxodromic and does not have a fixed point in common with $\beta\alpha^k$.

*Step* (iv).   Consider $\xi$ and $\delta = \eta\beta\alpha^k$ for fixed $k \geq K$ or $k \leq -K$, according to (iii).

If $\xi$ does not interchange the fixed points of $\delta$, we can take $N$ so large that either $\delta^n\xi$ or $\delta^n\xi^{-1}$ is loxodromic for $n \geq N = N(k)$.

Suppose instead that $\xi$ interchanges the fixed points of $\delta$ but not of $\eta\beta\alpha^{k+1} = \delta'$. Then replace $\delta$ by $\delta'$.

However, not both $\xi$ and $\eta\xi$ (nor equivalently, $\xi$ and $\eta^{-1}\xi$) can interchange the fixed points of both $\eta\beta\alpha^k$ and $\eta\beta\alpha^{k+1}$. For, if so, we apply Lemma 2.2.5 to $J = \xi$, $J_1 = \eta$ (or $\eta^{-1}$) and to both $\eta\beta\alpha^k$ and $\eta\beta\alpha^{k+1}$. That implies that the fixed points of both $\eta\beta\alpha^k$ and $\eta\beta\alpha^{k+1}$ coincide with fixed points $p, q$ of $\eta$. For this to occur, $\alpha$ fixes both $p$ and $q$, and then $\eta\beta$ must do so as well. But since $\eta$ itself fixes them, $\beta$ must also fix them. This is impossible.

We may assume one of $yx$ or $y^{-1}x$ is a simple loop. Depending on which, replace $\langle x, y\rangle$ by $\langle yx, y\rangle$ or $\langle y^{-1}x, y\rangle$. Correspondingly, replace $\xi$ by $\eta\xi$ or $\eta^{-1}\xi$. This returns us to one of the previous cases for $\delta = \eta\beta\alpha^k$ or $\eta\beta\alpha^{k+1}$.

4.5. *Cutting the surface.*  The loop $d^n x$ is freely homotopic to a simple, nondividing loop $d'$, disjoint from $d^n a$ and $ba^k$. Cutting $R$ along $d'$ results in a new surface $R^1$ with a handle $H = \langle d^n a, ba^k\rangle$ and two boundary components freely homotopic to $d^n$ and $yx^{-1}d^{-n}y^{-1}$ (or $y^{-1}x^{-1}d^{-n}y$). The corresponding transformations are $\delta^n\xi$ and $\eta\xi^{-1}\delta^{-n}\eta^{-1}$ (or $\eta^{-1}\xi^{-1}\delta^{-n}\eta$), which have the same trace. The common trace, however, can be made as large as desired (Lemma 2.1.1).

If the genus of $R^1$ exceeds one, repeat the process using the new $H$, and so on. At the end, we will have a surface $R^{g-1}$ with a handle $H = \langle a, b\rangle$ (using again the original notation) and $2(g-1)$ boundary components.

Orient all the boundary components so that $R^{g-1}$ lies to their right. Let $x, y, \ldots$ denote simple loops from the basepoint $O$ parallel to them but otherwise disjoint from each other and from $a$ and $b$. Our construction allows us to assume that the $\theta$-images $\theta(x)$, $\theta(y)$, $\ldots$ are all loxodromic. Pairwise they have the same trace, but the traces of different pairs can be assumed to be different.

## 5. The pants decomposition

5.1. We carry on from the situation left in §4.5. To start, adjust the special handle $H = \langle a, b \rangle$ as in §4.2 so that $\beta = \theta(b)$ does not send one fixed point of $\alpha = \theta(a)$ to the other. Orient $b$ so that it crosses $a$ from the right side of $a$ to the left; then the boundary of $R^{g-1}$ lies to the left of $c = b^{-1}a^{-1}ba$ and we have oriented the boundary so that $c$ lies to its right.

Choose simple loops $x, y \in \pi_1(R^{g-1}; O)$, each parallel to a boundary component and disjoint from each other and $a, b$, except at $O = a \cap b$. The orientations are such that $yba^k$ and $xba^k$ (but not $yb^{-1}a^k$ or $xb^{-1}a^k$ for $k \neq 1$) are homotopic to simple loops for all $k$ (see Figure 6). From §4.5 we know that $\xi = \theta(x)$ and $\eta = \theta(y)$ are loxodromic.

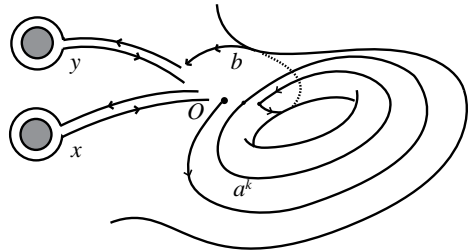

Figure 6. Connection of boundary to handle $H$

5.2. We begin by sorting out the following possibilities.

(1) If exactly one of $\eta\beta$ and $\xi\beta$ interchanges the fixed points of $\alpha$, assume that the one that does is $\xi\beta$. In this case, we claim that, for all sufficiently large $|k|$, the composition $\delta = \eta\beta\alpha^k$ does not fix either fixed point $p, q$ of $\xi$.

For if $\eta\beta\alpha^k$ fixes $p$ for two values of $k$, then $p$ itself must be fixed by $\alpha$, and then by $\eta\beta$ as well as by $\xi$. On the other hand, since $\xi\beta$ interchanges the fixed points $p$ and $p'$ of $\alpha$, we get $\xi\beta(p') = p = \xi(p)$, so $\beta(p') = p$. This contradiction to the known properties of the handle $H$ establishes the claim.

(2) If neither $\eta\beta$ nor $\xi\beta$ interchanges the fixed points of $\alpha$, then by interchanging $\xi$ and $\eta$ and relabeling if necessary, we may assume that for all sufficiently large $|k|$, the composition $\delta = \eta\beta\alpha^k$ does not fix both fixed points $p, q$ of $\xi$.

For suppose $\eta\beta\alpha^k$ fixes $p, q$ for two values of $k$, and, correspondingly, $\xi\beta\alpha^k$ fixes the two fixed points of $\eta$ for two other values of $k$. The first supposition implies that $p$ and $q$ are fixed by $\alpha$, then by $\eta\beta$, and of course by $\xi$. The second implies that $p$ and $q$ are fixed in addition by $\xi\beta$ and $\eta$. But $\eta\beta(p) = p = \eta(p)$ implies that $\beta$ itself fixes $p$, a contradiction.

It is important to note that if, in addition, $\eta\beta$ sends one fixed point of $\alpha$ to the other, then we may assume that $\delta = \eta\beta\alpha^k$, for $|k|$ large, does not fix even one fixed point of $\xi$. This is another application of the reasoning of (1).

(3) We defer consideration until §5.5 of the remaining case that both $\eta\beta$ and $\xi\beta$ interchange the fixed points of $\alpha$.

5.3. In this section and the next we will work with cases (1) and (2) of §5.2. That is, $\eta\beta$ does not interchange the fixed points of $\alpha$, and, for all large $|k|$, the composition $\delta = \eta\beta\alpha^k$ does not fix both fixed points of $\xi$, and if $\eta\beta$ sends one fixed point of $\alpha$ to the other, $\eta\beta\alpha^k$ does not fix either fixed point of $\xi$. Consider $d = yba^k$, for some $k$, and its $\theta$-image $\delta$. The effect of a Dehn twist of order $n$ about $d$ is

$$\langle \alpha, \beta\alpha^k \rangle \quad \mapsto \quad \langle \delta^n\alpha, \beta\alpha^k \rangle,$$
$$\langle \xi, \eta \rangle \quad \mapsto \quad \langle \delta^{-n}\xi\delta^n, \eta \rangle.$$

We will find $k$ and $n$ such that:

(i) $\beta\alpha^k$ is loxodromic;

(ii) $\delta = \eta\beta\alpha^k$ is loxodromic;

(iii) $\delta^n\alpha$ is loxodromic and has no common fixed point with $\beta\alpha^k$;

(iv) $\delta^{-n}\xi\delta^n\eta$ is loxodromic and has no common fixed point with $\eta$;

(v) $\delta^{-n}\xi\delta^n$ and $\eta$ generate a classical Schottky group.

(vi) $|\mathrm{tr}\,\delta^{-n}\xi\delta^n\eta|$ is unbounded in $|n|$.

Once this is accomplished, we will replace the handle $\langle a, b \rangle$ with $\langle d^n a, ba^k \rangle$, then remove from $R^{g-1}$ the pants determined by

$$(d^{-n}xd^n, y, d^{-n}xd^ny),$$

and repeat the process.

5.4. *Step* (i). The properties of the special handle $H$ (§5.1) and Lemma 2.1.1(i) imply that there is $K \geq 0$ such that $\beta\alpha^k$ is loxodromic for all $|k| \geq K$.

*Step* (ii). Since $\eta\beta$ does not interchange the fixed points of $\alpha$, we can choose $K$ above so large that $\eta\beta\alpha^k$ is loxodromic for $k \geq K$, $k \leq -K$, or both. In addition, for the admissible range of $k$, the composition $\eta\beta\alpha^k$ does not fix both fixed points of $\xi$.

*Step* (iii). This is identical with step (iii) of §4.3. There exists $K \geq 0$ such that $\delta = \eta\beta\alpha^k$ is loxodromic for any $k \geq K$, or any $k \leq -K$, or both.

Given $k$ in the admissible range, there exists $N = N(k) \geq 0$ such that, for all $|n| \geq N$, the element $\delta^n \alpha$ is loxodromic and has no fixed point in common with $\beta \alpha^k$.

*Step* (iv). Note that we may assume that $K$ is sufficiently large so that $\delta = \eta \beta \alpha^k$ has no fixed point in common with $\eta$ for $|k| \geq K$. For, if $\eta \beta \alpha^k$ fixes a fixed point $p$ of $\eta$ for two values of $k$, then $\alpha$ itself must fix $p$, and then $\beta$ must as well.

*Case* (1). $\eta \beta$ sends one fixed point of $\alpha$ to the other. In this case $\delta$ has no fixed point in common with $\xi$ (§5.3). Thus, by Lemma 2.1.3, the composition $\delta^{-n} \xi \delta^n \eta$ is loxodromic for all large $|n|$, while we must have $|n| \geq N$ to ensure that $\delta^n \alpha$ is loxodromic.

*Case* (2). $\eta \beta$ does not send one fixed point of $\alpha$ to the other. Then $\delta^n \alpha$ is loxodromic for all $|n| \geq N$, while $\delta^{-n} \xi \delta^n \eta$ is loxodromic either for $n \geq N$ or for $n \leq -N$, for sufficiently large $N$.

Finally, $\delta^{-n} \xi \delta^n \eta$ and $\eta$ have a fixed point in common only if $\delta^{-n} \xi \delta^n$ and $\eta$ do. The fixed points of $\delta^{-n} \xi \delta^n$ are $\delta^{-n}(p)$ and $\delta^{-n}(q)$, where $p$ and $q$ are the fixed points of $\xi$. If neither point is fixed by $\delta$, then, for sufficiently large $|n|$, neither $\delta^{-n}(p)$ nor $\delta^{-n}(q)$ will be fixed by $\eta$. On the other hand, if $p$, say, is fixed by $\delta$, the same conclusion holds because $\delta$ and $\eta$ do not share a fixed point.

*Step* (v). Since $\delta$ and $\eta$ have no fixed point in common, it follows from Lemma 2.1.3 and Corollary 2.1.5 that $\delta^{-n} \xi \delta^n$ and $\eta$ generate a classical Schottky group for sufficiently large $N$. Also the trace of $\delta^{-n} \xi \delta^n \eta$ can be made arbitrarily large, for sufficiently large $N$.

5.5. Now we turn to the case, left aside in §5.2, where both $\eta \beta = J$ and $\xi \beta = J_1$ interchange the fixed points of $\alpha$. Then $\eta = JJ_1 \xi$, where $JJ_1$ is loxodromic or the identity, and fixes the fixed points of $\alpha$.

At the start we arranged matters so that $yba^k$ is homotopic to a simple loop for all $k$. This is equally true of $(ba^k)^{-1} y(ba^k)$, and of $x(ba^k)^{-1} y(ba^k)$, which is homotopic to a simple loop bounding a triply connected region (pants) with boundary components corresponding to $x$ and $y$.

We claim that, in the present case, there exists $K \geq 0$ such that the corresponding transformation

$$\delta = \xi(\beta \alpha^k)^{-1} \eta(\beta \alpha^k) = \xi \alpha^{-k} \beta^{-1} J \alpha^k$$

is loxodromic for all $|k| \geq K$.

For $\xi$ has no fixed point in common with $\alpha$: Indeed, $\xi(p) = p = \alpha(p)$ would imply that $J_1 \beta^{-1}(p) = p$, in other words that $\beta(q) = p$, where $q$ is the other fixed point of $\alpha$. Similarly, $\beta^{-1}J$ has no fixed point in common with $\alpha$. Hence the assertion follows from Lemma 2.1.3.

We can take $K$ so large that, in addition, $\xi$ does not have a fixed point $p$ in common with

$$(\beta\alpha^k)^{-1}\eta(\beta\alpha^k) = \alpha^{-k}\beta^{-1}J\alpha^k$$

for $|k| \geq K$. For, since neither $\xi$ nor $\beta^{-1}J$ has a fixed point in common with $\alpha$, we have on the complement of $q_*, q^*$,

$$\lim_{k\to+\infty} \alpha^{-k}\beta^{-1}J\alpha^k = q_* \quad \text{and} \quad \lim_{k\to-\infty} \alpha^{-k}\beta^{-1}J\alpha^k = q^*,$$

where $q_* \neq p$ and $q^* \neq p$ denote the repulsive and attractive fixed points of $\alpha$.

For sufficiently large $K$ as dictated by Corollary 2.1.5, $\xi$ and $(\beta\alpha^k)^{-1}\eta(\beta\alpha^k)$ generate a classical Schottky group for all $|k| \geq K$.

As in §5.4(v), the trace magnitude of the transformation $(\beta\alpha^k)^{-1}\eta(\beta\alpha^k)\xi$ corresponding to the new boundary component of the pants can be made arbitrarily large, in particular in comparison with that of $\xi$ and $\eta$, which correspond to the boundary components on which the new pants was built.

Replace the handle $\langle a, ba^k \rangle$ by its conjugate $\langle (ba^k)^{-1}a(ba^k), ba^k \rangle$. The new pants is determined by $\langle x, (ba^k)^{-1}y(ba^k) \rangle$.

5.6. *The pants decomposition.* In §§5.3–5.4 we showed that, given any two boundary components of $R^{g-1}$, we could construct a pants with them as boundary components and such that the transformation corresponding to the third boundary component has trace of arbitrarily large magnitude. The surface remaining after this pants is removed is again of genus one, but with one fewer boundary components. Again choosing any two boundary components, we can construct another pants, and so on until all that remains is a surface of genus one with one boundary component: a handle.

For later requirements, we will specify the initial steps of the decomposition as follows: Group the $2(g-1)$ boundary components of $R^{g-1}$ into pairs, where the two components of each pair arise from cutting a handle of $R$. Construct first $(g-1)$ pants, one corresponding to each pair, which then comprise two of its boundary components. After this is done, finish the construction with any possible succession of pairings.

Each pair of boundaries of $R^{g-1}$ corresponds to transformations of the same trace, but we may assume from §4.5 that different pairs correspond to transformations of different traces. When each new boundary component forming a new pants is inserted, we can ensure by §5.4(v) and §5.5 that the trace magnitude of its corresponding transformation exceeds that corresponding to all previously inserted boundaries.

The combinatorics of the decomposition and a corresponding reorganization of the generating set for $\pi_1(R; O)$ will be discussed in §5.8 below.

5.7. *The final cut.* We are left with $2g-3$ pants and a handle $H$. Yet more adjustment to $H$ is necessary before gaining the assurance that one final cut will produce a pants decomposition $\{R_i\}$ called for in §1.7. Consider the handle $H = \langle a, b \rangle$ remaining at the end of the process. A simple loop $c \sim b^{-1}a^{-1}ba$ bounds $H$ on its right side (starting in §5.1 we specifically assumed that $b$ crosses $a$ from the right side of $a$ to the left). On the left side of $c$ is a pants with boundary components oriented so that the pants is to their left. Also, $c^{-1} \sim yx$, where $x$ and $y$ are simple loops parallel to the boundary components, and are disjoint from $c$, $b$ and $a$ except for a shared basepoint.

Set $\alpha = \theta(a)$, $\beta = \theta(b)$, $\sigma = \theta(c)$, $\xi = \theta(x)$, $\eta = \theta(y)$. We know that $\langle \xi, \eta \rangle$ is a Schottky group and that $\alpha$ and $\beta$ are loxodromic without a common fixed point. As in §4.2, we can assume that $\beta$ does not send one fixed point of $\alpha$ to the other.

By construction (see §5.4(iv)–(v) and §5.5), the trace magnitude of $\eta\xi$ exceeds that of $\eta$ and $\xi$. In particular, neither $\xi$ nor $\eta$ can be conjugate within $\mathrm{PSL}(2, \mathbb{C})$ to $\eta\xi = \alpha^{-1}\beta^{-1}\alpha\beta = \sigma^{-1}$.

We may assume that $\xi\beta = J$ does not interchange the fixed points of $\alpha$, and is not the identity. Otherwise, replace $\xi$ and $\eta$ by their conjugates $\sigma^m\xi\sigma^{-m}$ and $\sigma^m\eta\sigma^{-m}$, where $m$ is chosen so that $\sigma^m\xi\sigma^{-m}\beta = J_m$ neither interchanges the fixed points of $\alpha$ nor is the identity. To see that such an $m$ exists, consider $\sigma^m\xi\sigma^{-m}\xi^{-1} = J_m J$, which either has the same fixed points as $\alpha$, or has order two. The latter case is impossible because $\langle \xi, \eta \rangle$ is a Schottky group. Since $\xi$ and $\sigma = \eta\xi$ have no fixed points in common, the former is impossible as well, except perhaps for a finite number of values of $m$. Consequently, replace $\langle \xi, \eta \rangle$ by the conjugate group $\langle \sigma^m\xi\sigma^{-m}, \sigma^m\eta\sigma^{-m} \rangle$, and correspondingly $\langle x, y \rangle$ by the conjugate pants $\langle c^m x c^{-m}, c^m y c^{-m} \rangle$. Return again to the original notation.

Now we are ready to cut the handle $H$. But first, apply a Dehn twist of order $k$ about $a$. This changes $H$ to $\langle a, ba^k \rangle$.

Next, apply a Dehn twist of order $n$ about a simple loop $d \sim xba^k$. This results in the changes

$$\langle a, ba^k \rangle \quad \mapsto \quad \langle d^n a, ba^k \rangle.$$

Finally, cut the resulting handle along a simple loop freely homotopic to $ba^k$. This results in a pants whose fundamental group is

$$\langle ba^k, (d^n a)^{-1}(ba^k)^{-1}(d^n a) \rangle.$$

We claim that $k$ and $n$ can be chosen so that the groups representing the adjacent pants are now both Schottky groups

$$\langle \gamma, (\delta^n \alpha)^{-1}\gamma^{-1}(\delta^n \alpha) \rangle \quad \text{and} \quad \langle \xi, \delta^{-n}\eta\delta^n \rangle,$$

where $\gamma = \beta\alpha^k$ and $\delta = \xi\gamma$.

(i) There exists $K \geq 0$ such that $\delta = \xi\beta\alpha^k$ is loxodromic either for $k \geq K$ or for $k \leq -K$; for definiteness assume the former is true. Indeed we have already arranged matters so that $\xi\beta$ does not interchange the fixed points of $\alpha$.

(ii) For sufficiently large $K$ and $k \geq K$ the composition $\delta = \xi\gamma$ has no fixed points in common with $\xi$ or $\gamma$, and $\alpha\gamma\alpha^{-1}$ has no fixed points in common with $\alpha\beta\alpha^{-1}\beta^{-1}\xi^{-1} \neq \mathrm{id}$.

First note that neither $\beta\alpha^k$ nor $\alpha\beta\alpha^{k-1}$ can have the same fixed points for two values of $k$. For example, $\beta\alpha^k(p) = p = \beta\alpha^m(p)$ for $m \neq k$ implies that $\alpha(p) = p$, and then $\beta(p) = p$, which is impossible. Consequently, for sufficiently large $K$ and $k \geq K$, the element $\gamma = \alpha\beta^k$ does not share a fixed point with $\xi$ or $\eta$, nor $\alpha\gamma\alpha^{-1}$ with $\alpha\beta\alpha^{-1}\beta^{-1}\xi^{-1}$, provided this latter is not the identity. It follows that neither $\xi\gamma$ and $\gamma$, nor $\xi\gamma$ and $\xi$, can share fixed points either. Finally, $\alpha\beta\alpha^{-1}\beta^{-1}\xi^{-1} \neq \mathrm{id}$ because $\xi$ is not conjugate to $\eta\xi = \alpha^{-1}\beta^{-1}\alpha\beta$ (because they have unequal traces, as we have seen earlier in §5.7).

(iii) We show that $\langle \xi, \delta^{-n}\eta\delta^n \rangle$ is a Schottky group, either for all $n \geq N$ or all $n \leq -N$, for some $N \geq 0$; for definiteness we will assume the former.

For if $\delta$ has both its fixed points in common with $\eta$, then $\delta$ and $\eta$ commute and the group remains $\langle \xi, \eta \rangle$. If $\delta$ has one fixed point in common with $\eta$, say its repulsive fixed point $p$, the fixed points of $\delta^{-n}\eta\delta^n$ converge to $p$ as $n \to +\infty$. Since $p$ is not also a fixed point of $\xi$, the group is Schottky for large $n$. If $\delta$ has no fixed points in common with $\eta$, it is Schottky for all large $|n|$.

(iv) We show that $\langle \gamma, (\delta^n\alpha)^{-1}\gamma^{-1}(\delta^n\alpha) \rangle$ is a Schottky group for all $|n| \geq N$, for sufficiently large $N$ in (iii) and fixed $k \geq K$ from (i) and (ii).

For the fixed points of $(\delta^n\alpha)^{-1}\gamma^{-1}(\delta^n\alpha)$ are the images under $\alpha^{-1}\delta^{-n}$ of the fixed points of $\gamma$. As $n \to +\infty$ or $n \to -\infty$, these images converge to $\alpha^{-1}(p)$, where $p$ is the repulsive or attractive fixed point of $\delta$, since $\delta$ and $\gamma$ have no fixed points in common. If $\alpha^{-1}(p)$ is not a fixed point of $\gamma$, Corollary 2.1.5 implies that the group is Schottky for large $|n|$.

Suppose to the contrary that $\gamma\alpha^{-1}(p) = \alpha^{-1}(p)$, while $\xi\gamma(p) = p$. Then $\alpha\gamma\alpha^{-1}(p) = p$, while $p = \alpha\gamma\alpha^{-1}\gamma^{-1}\xi^{-1}(p) = \alpha\beta\alpha^{-1}\beta^{-1}\xi^{-1}(p)$. This does not occur, by (ii).

*Remark* 5.7.1. Had we not been so concerned about the final cut forming two of the boundary components of a single pants corresponding to a Schottky group, we would have proceeded more simply, as follows. Cut $H = \langle a, b \rangle$ along $a$ resulting in a pants $\langle a, b^{-1}a^{-1}b \rangle$. Pair boundary components of this with those of neighboring pants $\langle x, y \rangle$, to get two new pants $\langle a, y \rangle$ and $\langle b^{-1}a^{-1}b, x \rangle$. Apply to these the Dehn twist of order $m$ about $c \sim b^{-1}a^{-1}ba$. For all large $|m|$, the corresponding groups are easily seen to be Schottky.

5.8. *The combinatorics of pants decomposition.* We will systematically organize a generating set for the fundamental group of $R$ in terms of the pants decomposition $\{P_i\}$.

Start by fixing points $O_i, O'_i, O''_i$ on each component of $\partial P_i$, and disjoint simple auxiliary arcs from $O_i$ to $O'_i$ and $O''_i$. In terms of these auxiliary arcs, there is a unique path in $P_i$ between any two boundary components. Also, a component of $\partial P_i$ with an assigned orientation uniquely determines a loop from $O_i$, which we will take as the basepoint of $\pi_1(P_i; O_i)$. If $a_i$ and $b_i$ are two boundary components of $P_i$, an orientation of $a_i$ uniquely determines an orientation of $b_i$ such that $b_i a_i$ is homotopic to a simple loop around the third (here making use of the auxiliary arcs).

If the components $a$ of $\partial P_i$ and $a'$ of $\partial P_j$ correspond to the same simple loop on $R$, choose the points $O \in a$ and $O' \in a'$ to correspond to the same point on $R$.

Let $\mathbb{T}$ denote the trivalent graph of genus $g$ corresponding to the pants decomposition $\{P_i\}$: each vertex of $\mathbb{T}$ corresponds to one of the pants $P_i$, and each edge corresponds to a pair $(a, a')$ of boundary components, one on each pants corresponding to an endpoint. Two boundary components are paired $(a, a')$ if and only if they correspond to the same simple loop on $R$.

$\mathbb{T}$ has $2g - 2$ vertices and $3g - 3$ edges. Exactly $g$ of the vertices have one-edge loops attached to them; this is a consequence of the particular combinatorics of the decomposition. We call these vertices *extreme*.

Remove from $\mathbb{T}$ those $g$ one-edge loops; the result $\mathbb{T}_0$ is a maximal (connected) tree. The extreme vertices of $\mathbb{T}$ are those that are extreme in $\mathbb{T}_0$ in the sense that only one edge of $\mathbb{T}_0$ is attached to the vertex.

Designate one of these extreme vertices as the *root* $v_0$ of $\mathbb{T}_0$: for example, the vertex corresponding to the last handle we cut. There is a unique simple path in $\mathbb{T}_0$ from any vertex to the root.

Denote the pants corresponding to the vertex $v$ by $P(v)$. Consider the vertices $v' \neq v$ whose shortest path to $v_0$ contains $v$. Mark the boundary components of $P(v)$ where these shortest paths first cross; we will use these shortest paths below. If $v$ is not extreme, two of the three boundary components of $P(v)$ will be marked. If $v$ is extreme but $v \neq v_0$, none of the boundary components will be marked. Exactly one of the boundary components of $P(v_0)$ will be marked.

Making use of the auxiliary arcs in the $\{P_i\}$, the simple edge-arc in $\mathbb{T}_0$ from the vertex $v_i = P_i$ to $v_j = P_j$ uniquely determines a simple arc in $R$ from $O_i$ to $O_j$. Likewise, a simple edge-loop in $\mathbb{T}$ uniquely determines a simple loop in $R$.

Let $P_0$ be the pants corresponding to the root $v_0$, and $O = O_0$ the designated basepoint for its fundamental group. Take also $O$ as the basepoint of the fundamental group of $R$. As we have seen, $\mathbb{T}_0$ uniquely determines a simple arc $c_i$ in $R$ from $O$ to each $O_i$. Thus, a simple loop $a_i \in \pi_1(P_i; O_i)$ can be uniquely associated with $c_i^{-1}a_i c_i \in \pi_1(R; O)$. Suppose $e_i$ is one of the $g$ edge-loops of $\mathbb{T}$, with both end points on the same vertex $v_i$. Likewise with the help of the auxiliary arcs in $P_i = P(v_i)$, the edge $e_i$, with an assigned orientation, uniquely determines a loop $c_i'^{-1}e_i c_i' \in \pi_1(R; O)$.

The totality of elements $c_i^{-1}a_i c_i$ from oriented boundary components of pants $\{P_i\}$ plus $g$ elements $c_i'^{-1}e_i c_i'$ from edges $e \notin \mathbb{T}_0$ generate $\pi_1(R; O)$.

## B. Pants configurations from Schottky groups

## 6. Joining overlapping plane regions

6.1. In this section we will describe a method of using covering surfaces to separate two overlapping plane regions which are acted on by a common Möbius transformation. It is no restriction to describe the process with the loxodromic transformation $\alpha : z \mapsto \lambda^2 z$, with $|\lambda| > 1$ and fixed points $0$ and $\infty$. Let $T$ or $T(\alpha)$ denote the quotient torus

$$T = (\mathbb{C} \setminus \{0\})/\langle \alpha \rangle,$$

and $\pi$ the projection from $\mathbb{C} \setminus \{0\}$. Denote the simple compressing loop $\pi(\{z : |z| = 1\})$ in $T$ by $c$. A noncontractible simple loop on $T$ lifts to a closed curve in $\mathbb{C} \setminus \{0\}$ if and only if it is freely homotopic (or homologous) to $\pm c$.

If a simple loop $a$ is not of this type, $a^* = \pi^{-1}(a)$ is a simple $\alpha$-invariant arc. If $a$ is given the orientation dictated by $\alpha$, the arc $a^*$ is directed toward the attractive fixed point.

Conversely, if $a_1^*$ is a simple, $\alpha$-invariant arc in $\mathbb{C}$ directed toward the attractive fixed point, $a_1 = \pi(a_1^*)$ is a simple loop freely homotopic (or homologous) to the result of applying to $a$ the Dehn twist about $c$ of some order $n$: namely, $a_1 \sim a + nc$.

6.2. Let $S_N$ denote the $N$-sheeted cover of the sphere $S_1 = \mathbb{S}^2$, branched over the fixed points $0$ and $\infty$ of $\alpha$. Topologically, $S_N$ is again a sphere. The map $z \mapsto z^{1/N} = w$ sends $S^N$ back to $S_1$; it is conformal except at $0$ and $\infty$. The cyclic group of cover transformations is conjugated to the group of rotations $\langle w \mapsto e^{2\pi i/N}w \rangle$.

The transformation $\alpha$ lifts to an automorphism $\alpha^*$ of $S_N$, determined up to composition with cover transformations. It is conjugated to the loxodromic transformation $w \mapsto \lambda^{2/N}w$, which in turn is determined only up to composition with cover transformations.

Consider the torus $T_N = T_N(\alpha)$, defined by

$$T_N = (S_N \setminus \{0, \infty\})/\langle \alpha^* \rangle.$$

It is the $N$-sheeted torus over $T$, uniquely determined by the properties that $a$ lifts to exactly $N$ mutually disjoint simple loops and $c^N$ lifts to one simple loop.

For the following lemma $a$ and $c$ are simple loops on $T$ as before: $c$ is the projection of the unit circle and $a$ is the projection of a simple, $\alpha$-invariant arc $a^* \subset \mathbb{C}$, positively oriented by $\alpha$. If the simple arc $a_1$ crosses $a$ transversely at every point of intersection, the geometric intersection number is defined as the number of points of intersection. We assume this number is finite.

LEMMA 6.2.1.   *Suppose $a_1$ is freely homotopic and transverse to $a$, with geometric intersection number $n$.  Set $N = 2n + 1$.  Then there is a lift $a'$ of $a$ and a lift $a_1'$ of $a_1$ that are disjoint, freely homotopic simple loops in $T_N$. Correspondingly, there is a lift $a^*$ of $a$ and $a_1^*$ of $a_1$ to $S_N$ that are disjoint, $\alpha^*$-invariant simple arcs.*

*Remark* 6.2.2.  A more precise measure of intersection would be to set

$$n = \max |m(\tau)|,$$

where $\tau \subset a_1$ is a segment whose endpoints don't lie on $a$, $m(\tau)$ is the algebraic intersection number of $\tau$ and $a$, and the maximum is taken over all such connected segments $\tau$ of $a_1$.

*Proof.*  Fix a lift $a'$ of $a$ to $T_N$ or a lift $a^*$ to $S_N$.  We can label the $N = 2n + 1$ sheets on $T_N$ or on $S_N$ over $T \setminus \{a\}$ in cyclic order, starting to the left of $a'$ or of $a^*$.  A point in the $(n + 1)$-st sheet can be connected to one on $a'$ or on $a^*$ only by crossing $n$ other lifts of $a$.  Fix $p \in a_1 \setminus a$, and the point $p'$ or $p^*$ lying over $p$ in the $(n + 1)$-st sheet.  The endpoint of the arc $\tilde{a}_1$ lying one-to-one over $a_1 \setminus \{p\}$ and starting at $p'$ or $p^*$ also lies in the $(n + 1)$-st sheet, because $a_1$ is freely homotopic to $a$; in $T_N$, the arc $\tilde{a}_1$ closes up to form a simple loop.  The conclusion is a direct consequence.   □

Note that without the condition that $a_1$, positively oriented by $\alpha$, be freely homotopic to $a$, the conclusion of the lemma is false.  Instead, the following is true.

COROLLARY 6.2.3.   *Suppose, more generally, that the simple loop $a_1 \subset T$, transverse to $a$, is the projection of an $\alpha$-invariant arc in $\mathbb{C} \setminus \{0\}$.  There exists $N = N(a, a_1) \geq 1$ and $m \in \mathbb{Z}$ such that $\delta^m a_1$ and $a$ have disjoint lifts on $T_N$ and $S_N$, where $\delta$ denotes the Dehn twist about $c$.*

6.3.

LEMMA 6.3.1.   *Suppose $a$ and $a_1$ are $\alpha$-invariant simple arcs in $\mathbb{C} \setminus \{0\}$, the lifts of freely homotopic transverse loops in $T(\alpha)$ with geometric intersection number $n$. Suppose $a$ is contained in the boundary of a simply connected* *region* $P \subset \mathbb{C} \setminus \{0\}$ *lying to its left, while* $a_1$ *is contained in the boundary of a simply connected region* $P_1$ *lying to its right. Set* $N = 2n + 4$. *Then on* $S_N \setminus \{0, \infty\}$ *there are disjoint lifts* $a^*$ *of* $a$ *and* $a_1^*$ *of* $a_1$ *with the property that the corresponding lifts* $P^*$ *of* $P$ *and* $P_1^*$ *of* $P_1$ *that contain* $a^*$ *and* $a_1^*$ *in their respective boundaries are disjoint as well.*

*Proof.* Fix a lift $a^*$ of $a$ in $S_N \setminus \{0, \infty\}$ and let $E$ be the generator of the order-$N$ cyclic group of cover transformations with the property that $a^*$ and $Ea^*$ bound to the left of $a^*$ a lift $\sigma^*$ of $\mathbb{C} \setminus \{a\}$, which we will refer to as the first sheet of the covering. In cyclic order to the left of $a^*$ the lifts of $a$ are $Ea^*, \ldots, E^{N-1}a^*$, and the corresponding sheets are $\sigma^*, E\sigma^*, \ldots, E^{N-1}\sigma^*$.

Denote by $P^*$ the lift of $P$ adjacent to $a^*$ on its left side. Necessarily $P^*$ lies entirely in the first sheet $\sigma^*$.

If $a_1$ is disjoint from $a$, then $N = 4$ (although $N = 3$ will do). Let $a_1^*$ be the lift of $a_1$ lying in the third sheet $E^2\sigma^*$, and $P_1^*$ the lift of $P_1$ adjacent to $a_1^*$ on its right side. $P_1^*$ lies in the sector bounded by $a_1^*$ and $E^{-1}a_1^*$, which lies in the second sheet $E\sigma^*$. Hence $P_1^*$ is disjoint from $P^*$ (see Figure 7).

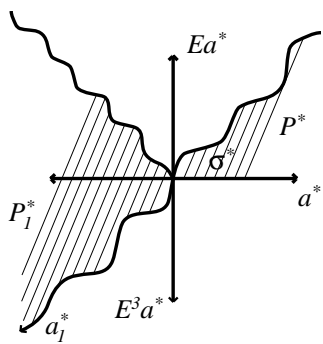

Figure 7. Separation of regions when $N = 4$

More generally, choose $p \in a_1 \cap a$ and let $p^*$ denote the point over $p$ on $E^{n+2}a^*$. Let $a_1^*$ denote the lift of $a_1^*$ through $p^*$; then $a_1^*$ does not intersect $E^2a^*$ or $E^{-2}a^*$. Consequently, $E^{-1}a_1^*$ does not intersect $Ea^*$. Let $P_1^*$ be the lift of $P_1$ adjacent to $a_1^*$ on its right side; $P_1^*$ lies in the sector between $a_1^*$ and $E^{-1}a^*$. Therefore $P_1^*$ is disjoint from $P^*$. □

Note that we have not optimized the choice of $N$, which can be done in particular cases.

COROLLARY 6.3.2. *In the hypotheses of Lemma* 6.3.1, *assume that not only* $a$ *and* $a_1$ *but also* $P$ *and* $P_1$ *are* $\alpha$-*invariant in* $\mathbb{C} \setminus \{0\}$. *There is a lift* $\alpha^*$ *of* $\alpha$ *to* $S_N$ *that leaves* $P^*$ *and* $P_1^*$ *invariant.*

*Proof.* Let $\alpha^*$ be the lift of $\alpha$ that maps the first sheet $\sigma^*$ onto itself, and hence $P^*$ onto itself. Necessarily $\alpha^*$ maps every sheet $E^k\sigma^*$ onto itself, and hence $P_1^*$ onto itself as well.                                            □

6.4. *Joining overlapping regions.* In this section we will build a prototype for the procedure that forms the basis of §8. It is typical of tricks used in classical function theory and is a generalization of a technique applied to Möbius groups called grafting [Mas1], [He1], [Go1].

Consider the hypotheses of Lemma 6.3.1: $a$ and $a_1$ are $\alpha$-invariant simple arcs in $\mathbb{C}\setminus\{0\}$ directed toward $\infty$, and one does not spiral around the other (an informal way of saying that they arise from freely homotopic loops in $T$). The region $P$ lies to the left of $a$, and $P_1$ to the right of $a_1$. Like $a$ and $a_1$ themselves, $P$ and $P_1$ can badly overlap each other.

However, on $S_N$, $P^*$ and $P_1^*$ are disjoint. Let $Q^*$ be the region on $S_N$ that lies to the right of $a^*$ and to the left of $a_1^*$: then $P_1^* \cup a_1^* \cup Q^* \cup a^* \cup P^*$ is a simply connected region $R^*$ in $S_N \setminus \{0, \infty\}$. According to Corollary 6.3.2, if $P$ and $P_1$ are $\alpha$-invariant, $\alpha^*$ is a conformal automorphism of $R^*$.

Let $g : \mathbb{H}^2 \to R^*$ be a Riemann map, where the hyperbolic plane $\mathbb{H}^2$ is realized as the unit disk. Then $g^{-1}\alpha^* g$ is a hyperbolic Möbius transformation $\alpha_0^*$ in $\mathbb{H}^2$. Let $\pi : S_N \to \mathbb{S}^2$ denote the projection. Then $f = \pi \circ g$ is a locally univalent meromorphic function on $\mathbb{H}^2$ with the property that

$$f\alpha_0^*(z) = \alpha f(z)$$

for all $z \in \mathbb{H}^2$. That is, $f$ determines a complex projective structure on $\mathbb{H}^2$ that induces the isomorphism between cyclic groups $\langle\alpha_0^*\rangle \to \langle\alpha\rangle$.

We have *joined together* the annular regions $P/\langle\alpha\rangle = P^*/\langle\alpha^*\rangle$ and $P_1/\langle\alpha\rangle = P_1^*/\langle\alpha^*\rangle$ by means of the annulus $Q^*/\langle\alpha^*\rangle$, which attaches to the boundary components $a/\langle\alpha\rangle$ and $a_1/\langle\alpha\rangle$.

## 7. Pants within rank-two Schottky groups

7.1. Suppose $\langle\alpha_1, \alpha_2\rangle$ is a two-generator classical Schottky group acting on its regular set $\Omega \subset \mathbb{S}^2$. The quotient surface $R = \Omega/\langle\alpha_1, \alpha_2\rangle$ has genus two (and bounds the handlebody $R_+ = \mathbb{H}^3/\langle\alpha_1, \alpha_2\rangle$ if the group is extended to hyperbolic three-space).

There are round circles $b_1^*$ and $b_2^*$, mutually disjoint in $\Omega$, with the following property: The two pairs of circles $(b_1^*, \alpha_1 b_1^*)$ and $(b_2^*, \alpha_2 b_2^*)$ are mutually disjoint with mutually disjoint interiors in $\mathbb{S}^2$, and $\alpha_i$ maps the exterior of $b_i^*$ onto the interior of $\alpha_i b_i^*$. The circles $b_1^*$ and $b_2^*$ are lifts of mutually disjoint, nondividing simple loops $b_1$ and $b_2$ in $R$. These bound disks in $R_+$ and for that reason are called *compressing loops.*

Let $a_1$ and $a_2$ be simple, nondividing loops in $R$ such that $a_1 \cap (a_2 \cup b_2) = \emptyset$ and $a_2 \cap (a_1 \cup b_1) = \emptyset$, while $a_i$ crosses $b_i$ transversely at a single point. Then $a_1$ and $a_2$ have lifts $a_1^*$ and $a_2^*$ to $\Omega$ uniquely determined by the condition that they are $\alpha_1$-invariant and $\alpha_2$-invariant simple arcs, respectively. Let $\delta_i$ denote the Dehn twist about $b_i$. Then, for example, $\delta_1^n a_1$ can be used in place of $a_1$: it, too, can be taken to be a simple loop disjoint from $a_2$ and $b_2$, meeting and there crossing $b_1$ at a single point. It too has a uniquely determined $\alpha_1$-invariant lift $(\delta^n a_1)^*$ in $\Omega$. (More generally, the simple loop $a_1'$ has an $\alpha_1$-invariant lift if and only if $a_1'$ is freely homotopic to $a_1$ within the handlebody $R_+$.)

7.2. *Finding pants.* Assign $a_1$ and $a_2$ their *positive orientation*, that is, the one that directs $a_1^*$ and $a_2^*$, their $\alpha_1$- and $\alpha_2$-invariant lifts, toward the attractive fixed points of $\alpha_1$ and $\alpha_2$. We can join $a_1$ and $a_2$ to a common basepoint $O \in R$ so that the resulting simple loops $a_1'$ and $a_2'$ have the property that $a_2' a_1'$ is homotopic to a simple loop $a_3'$; this loop $a_3'$ is then freely homotopic to a simple loop $a_3$ that, together with $a_1$ and $a_2$, divides $R$ into two pants $P$ and $P'$; also, $a_3$ has an $\alpha_2\alpha_1$-invariant lift $a_3^*$ and an $\alpha_1\alpha_2$-invariant lift $\alpha_1 a_3^*$ (Figure 8).

Note that the free homotopy class of $a_3$ on $R$ is not uniquely determined by that of $a_1$ and $a_2$: we can change $a_3$ by applying Dehn twists about a simple dividing loop homotopic to $b_1'^{-1} a_1'^{-1} b_1' a_1'$ without affecting $a_1$ or $a_2$. We can also change $a_3$ by applying Dehn twists about $b_1$ or $b_2$, but that will change $a_1$ or $a_2$ as well. In any case there is an $\alpha_2\alpha_1$-invariant lift of $a_3$.

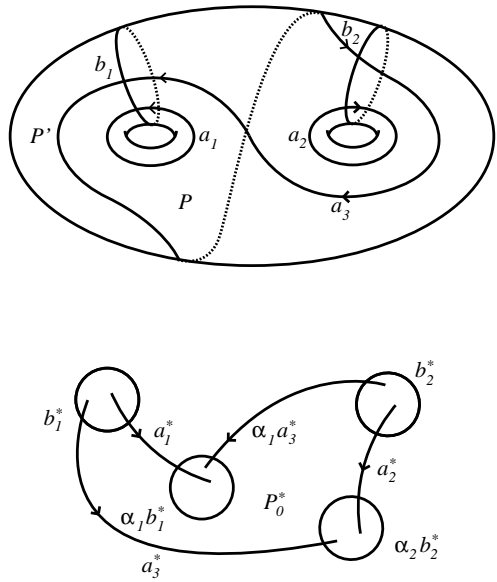

Figure 8. Pants determined by Schottky group

Let $P$ denote the pants lying to the right of $a_1$ and $a_2$, and to the left of $a_3$. Of course we may assume that $b_1 \cap P$ and $b_2 \cap P$ are simple arcs. There is a lift $P_0^*$ of $P$ to $\Omega$ that is an "octagon" bounded by connected segments of $a_1^*, a_2^*, a_3^*, \alpha_1 a_3^*$ and $b_1^*, b_2^*$ (see Figure 8, bottom). The orbit of $P_0^*$ (adding its boundary arcs on $b_1^*$ and $b_2^*$) under $\langle \alpha_1, \alpha_2 \rangle$ is a simply connected region $P^*$ that is the universal cover of $P$.

### 7.3. *Isomorphisms are geometric.* We summarize the analysis of §7.2 as follows.

LEMMA 7.3.1.   *Let $Q$ be a pants with oriented boundary components* $(d_1, d_2, d_3)$, *and choose generators* $d_1', d_2', d_3'$, $d_3' \sim d_2' d_1'$, *for* $\pi_1(Q; O)$ *such that* $d_i'$ *is parallel to* $d_i$, $O \in Q$. *Suppose $\theta$ is the isomorphism of $\pi_1(Q; O)$ onto the Schottky group* $\langle \alpha_1, \alpha_2 \rangle$ *determined by the correspondence* $\theta(d_1') = \alpha_1$, $\theta(d_2') = \alpha_2$. *Then there is a pants $P$ in $R = \Omega / \langle \alpha_1, \alpha_2 \rangle$ bounded by simple loops* $(a_1, a_2, a_3)$, *positively oriented by $\alpha_1, \alpha_2, \alpha_2 \alpha_1$, and a homeomorphism* $h : \bar{Q} \to \bar{P}$ *taking $d_i$ (with its orientation) to $a_i$, $i = 1, 2, 3$, which induces $\theta$: there is a point $O^* \in P^* \subset \Omega$ over $h(O) \in P$ such that the lift of $h(d_i')$ from $O^*$ terminates at $\alpha_i(O^*)$, for $i = 1, 2, 3$.*

*Proof.* In §7.2 we observed the following convention for finding pants $P$ and $P'$ in a Schottky group with designated generators $\alpha_1$ and $\alpha_2$. The three boundary components have $\alpha_1$-, $\alpha_2$-, and $\alpha_2 \alpha_1$-invariant lifts, positively oriented by $\alpha_1$, $\alpha_2$ and $\alpha_1 \alpha_2$, respectively. If $a_1$ and $a_2$ are represented by generators $a_1'$ and $a_2'$ in $\pi_1(P; O)$ or $\pi_1(P'; O')$, then $a_2' a_1'$ is homotopic to a simple loop parallel to $a_3$. The two pants $P$ and $P'$ are distinguished in that one lies to the right of $a_1$ and $a_2$, and to the left of $a_3$, while the opposite holds for the other.

The orientations of the $d_i$ can temporarily be reversed as necessary so that $Q$ lies to the right of $d_1, d_2$ and left of $d_3$. Make the corresponding temporary replacements of $\alpha_i$ by $\alpha_i^{-1}$. Now find a pants $P$ meeting the requirements, and then return to the original designations.   □

### 7.4. *Two groups with a common generator: Compatibility.* Consider two Schottky groups $\langle \alpha_1, \alpha_2 \rangle$ and $\langle \alpha_2, \alpha_3 \rangle$ with a common generator $\alpha_2$. Denote the regular sets in $\mathbb{S}^2$ by $\Omega$ and $\Omega'$, and set $R = \Omega / \langle \alpha_1, \alpha_2 \rangle$ and $R' = \Omega' / \langle \alpha_2, \alpha_3 \rangle$. Choose simple loops $(a_1, b_1, a_2, b_2)$ in $R$ and $(a_2', b_2', a_3', b_3')$ in $R'$ as in §7.1; here $a_j'$ and $a_i$ are taken with their positive orientations from $\alpha_j$ and $\alpha_i$. Find, as in §7.2, a pants $P \subset R$ lying, say, to the right of $a_2$, and then a pants $P' \subset R'$ lying to the left of $a_2'$.

*Definition* 7.4.1. As above, suppose $a_2$ and $a'_2$ are simple loops on $R$ and $R'$, with $\alpha_2$-invariant lifts $a_2^*$ and $a_2'^*$ to $\Omega$ and $\Omega'$, respectively. The loops $a_2$ and $a'_2$ are *compatible* (with respect to $\alpha_2$) if the projections of $a_2^*$ and $a_2'^*$ (that is, the embeddings of $a_2$ and $a'_2$) in the torus $T(\alpha_2)$, are freely homotopic there.

Recall that $T(\alpha_2) = (\mathbb{S}^2 \setminus \{p, q\})/\langle \alpha_2 \rangle$, where $p$ and $q$ are the fixed points of $\alpha_2$. Let $\delta_2$ denote the Dehn twist about $b'_2$ on $R'$. In general $a'_2$ will not be compatible with $a_2$. However,

LEMMA 7.4.2. *The loop $a'_2$ on $R'$ can be made compatible with $a_2$ on $R$: $a_2$ is compatible with $\delta_2^m a'_2$ for a unique value of $m$.*

*Proof.* Let $\delta_2$ denote the Dehn twist about $b'_2$ on $R'$. Note that $b'_2$ embeds as a simple loop on $T(\alpha_2)$, so that $\delta_2$ can be taken to act on $T(\alpha_2)$ as well as on $R'$. For exactly one value of $m$, the loop $\delta_2^m a'_2$ will be compatible with $a_2$. □

*Remark* 7.4.3. We emphasize that the process of making $a_2$ and $a'_2$ compatible affects only one of the surfaces: say $R'$. And, on $R'$, it affects only $a'_2$, not $a'_3$. However, the third boundary component $c'$ of the pants $P'$ is affected. Indeed, there is a lift $c'^*$ to $\Omega'$ invariant under $\alpha_3 \alpha_2$. The simple loop $b'_2$ which crosses $c'$ once also embeds in $T(\alpha_3 \alpha_2)$, and the twist $\delta_2$ equally can be taken to act on the torus $T(\alpha_3 \alpha_2)$. Thus, under the action of $\delta_2^m$ on $R'$, the loop $c'$ changes to $\delta_2^m c'$; the pants $\delta_2^m P'$ is bounded by $\delta_2^m a'_2$, $\delta_2^m c'$, and $a'_3$.

7.5. *Compatibility conditions on one pants.* Consider a Schottky group $\langle \alpha, \beta \rangle$ and a pants $P$ in $\Omega/\langle \alpha, \beta \rangle$, as in Figure 9.

Denote the boundary components of $P$ by $a, b, c$, with the orientation indicated. With respect to these curves, we can find compressing loops $x$ and $y$ (which lift to closed loops in $\Omega$), with the orientations and intersections indicated in the figure.

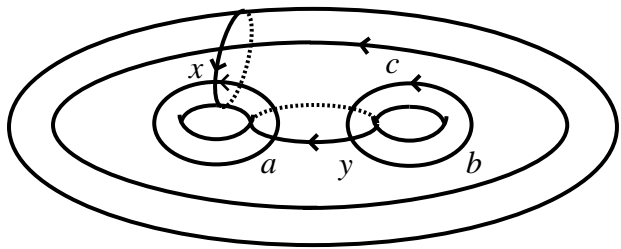

Figure 9. Pants and compressing loops

A Dehn twist of order $p$ about $x$ composed with a twist of order $q$ about $y$ has the following effect on $a, b, c$:

$$a \mapsto \delta^{p+q}a, \quad b \mapsto \delta^{-q}b, \quad c \mapsto \delta^{p}c.$$

Here we use the notation $\delta^{k}t$ to denote the effect on the simple oriented loop $t$ of a twist of order $k$ about an oriented simple *compressing* loop crossing $t$ once, from its right side to its left; geometrically the result is realized and accounted for on the torus $T(\tau)$ that is associated with $t$.

Suppose $b$ and $c$ are to be paired with boundary components $b'$ and $c'$ on other pants, where $\delta^{m}b$ is compatible with $b'$ and $\delta^{n}c$ is compatible with $c'$. This can be fulfilled simultaneously in $P$ by setting $p = n$ and $q = -m$. The effect on $a$ is to replace it by $\delta^{n-m}a$. That is, compatibility for two boundary components of $P$ can always be achieved, but then the state of the third boundary component is determined.

Suppose instead that $c$ is to be paired with $c'$ on another pants, with compatibility requirement $c' = \delta^{n}c$, while $b$ is to be paired with $a$ with compatibility requirement $a = \delta^{m}b$. In terms of §7.4, this means that there is a transformation $\gamma$ with $\alpha = \gamma\beta\gamma^{-1}$, where $a$ and $b$ have been determined by $\alpha$ and $\beta$, respectively. That is, there is an $\alpha$-invariant lift $a^{*}$ and a $\beta$-invariant lift $b^{*}$, and the two can be compared in terms of the $\alpha$-invariant arcs $a^{*}$ and $\gamma b^{*}$.

Therefore $p = n$, while $q$ is determined by the condition

$$-q = n + q + m, \quad \text{or} \quad q = -\tfrac{1}{2}(m+n).$$

A solution $q \in \mathbb{Z}$ exists if and only if $m + n$ is even, that is, if $m$ and $n$ have the same parity.

In other words, the algebraic sum $[(p+q) - q + p] = 2p$ of the Dehn twists that can be applied effectively to the boundary components of a pants is even. Consequently, if the requirements for compatibility in a pants demand that the algebraic sum be odd, those requirements cannot be met.

*Remark* 7.5.1. There is also a compressing loop $u$ in $\Omega/\langle\alpha, \beta\rangle$ that divides, separating $a$ and $b$ while crossing $c$ twice. A Dehn twist about $u$ leaves $a$ and $b$ unchanged, but changes the homotopy type of $c$ and $P$ on the surface $\Omega/\langle\alpha, \beta\rangle$. Yet it leaves unchanged the free homotopy type of the projection $c_{*}$ of $c$ to its associated torus $T(\beta\alpha)$.

For on $T(\beta\alpha)$, there are two representatives of $u$, $u_{*1}$ and $u_{*2}$. They are disjoint, parallel and oriented opposite one another: one crosses $c_{*}$ from right to left, the other from left to right. A Dehn twist about $u$ is reflected by twists on $T(\beta\alpha)$ about $u_{*1}$ and $u_{*2}$. But, because of their opposite orientations, these twists cancel, leaving the free homotopy class of $c_{*}$ unchanged.

In short, twists about $u$ have no effect on compatibility questions.

7.6. *A compatibility condition on identical pants.* For later application in §8, consider the following augmentation to the second situation of §7.5, where $\alpha = \gamma\beta\gamma^{-1}$. In the conjugate group $\gamma\langle\alpha,\beta\rangle\gamma^{-1}$, consider the pants $P_1$ that corresponds precisely to $P$, distinguishing corresponding elements by the subscript. Suppose, as before, that $c$ and $c_1$ are to be paired with $c'$ and $c'_1$ on other pants $P'$, $P'_1$, but now with the same compatibility requirements: $c' = \delta^n c$ and $c'_1 = \delta^n c_1$. Instead of pairing $b$ with $a$ as before, pair $b_1$ with $a$ and $b$ with $a_1$. Because the two groups are virtually identical, the compatibility requirements are $a = \delta^m b_1$ and $a_1 = \delta^m b$.

The result of Dehn twists of order $p$ and $q$ about $x$ and $y$, and of order $p_1$ and $q_1$ about $x_1$ and $y_1$, is as calculated in §7.5. We must have $p = p_1 = n$. That leaves, for $q$ and $q_1$, the equation

$$-q_1 = n + q + m, \quad \text{or} \quad q + q_1 = -(m + n).$$

In this case there are always solutions: for example, $q = -m$ and $q_1 = -n$.

# 8. Building the pants configuration

8.1. *What remains to be done?* In §5.7 the combinatorics of the pants decomposition $\{P_i\}$ of $R$ found in Part A was described as a trivalent graph $\mathbb{T}$ arising from a tree $\mathbb{T}_0 \subset \mathbb{T}$ by the addition of $g$ edges, one attached to each extreme vertex. The universal cover of $\mathbb{T}$ is reflected in the combinatorics of their lifts $\{Q_i^*\}$ in the universal cover $\mathbb{H}^2$ of $R$, that is, in how the lifts fit together.

Corresponding to each lift $Q_i^*$ is the Schottky group $\theta(\mathrm{Stab}\, Q_i^*)$, which in turn stabilizes the lift $P_i^*$ of a pants $P_i$ in its quotient surface. Using the technique of §6, our goal is to follow the information in $\mathbb{T}$, or the combinatorics of $\{Q_i^*\}$ in $\mathbb{H}^2$, to build a simply connected Riemann surface $\mathcal{J}$. This will be the universal cover of a surface obtained by joining together the pants $\{P_i\}$ by attaching auxiliary cylinders.

However, to join a boundary component $a$ of $P_i$ to $a'$ of $P_j$ (or perhaps to $a'$ of $P_i$), it is necessary that $a$ and $a'$ be compatible in the sense of §7.4. It is not necessarily true that the totality of compatibility conditions can be satisfied.

In §§8.2–8.4, typical cases of joining pants will be described, before we draw the general conclusions in §§8.4–8.6. In §9, we will show how to add branch points when needed.

8.2. *Joining pants.* We continue with the situation of §7.4. There, we found simple loops $a_3 \sim a_2 a_1$ on $R$ and $a'_4 \sim a'_3 a'_2$ on $R'$ such that $(a_1, a_2, a_3)$ bound pants $P \subset R$ lying to the right of $a_1$ and $a_2$, while $(a'_2, a'_3, a'_4)$ bound pants $P' \subset R'$ lying to the left of $a_2'$ and $a_3'$. Here $(a_1, a_2, a_2', a_3')$ are positively oriented by generators $\alpha_1$, $\alpha_2$ and $\alpha_3$. According to Lemma 7.4.2, the loop $a_2'$ can be taken compatible with $a_2$. We will now show how to join the pants $P$ and $P'$ by attaching a cylinder to the left side of $a_2$ and the right side of $a_2'$.

Let $P^*$ denote the region in $\Omega$ over $P$ and $P'^*$ the region over $P'$ in $\Omega'$. Both $P^*$ and $P'^*$ are simply connected, as they represent the respective universal covers.

We are in a position to apply Lemma 6.3.1 to $P^*$ and $P'^*$. There exists an $N$-sheeted $S_N$ of $\mathbb{S}^2$, branched over the fixed points of $\alpha_2$, on which there are disjoint lifts $a_2^{**}$ and $a_2'^{**}$ of $a_2^*$ and $a_2'^*$ that border disjoint lifts $P^{**}$ and $P'^{**}$ of $P^*$ and $P'^*$: the projections $P^{**} \to P^*$ and $P'^{**} \to P'^*$ are homeomorphisms. Equally well, $P^{**}$ and $P'^{**}$ represent the universal covers of $P$ and $P'$.

Next, take the sector $Q^{**}$ on $S_N$ lying between the left side of $a_2^{**}$ and the right side of $a_2'^{**}$, and form

$$Q_1^{**} = P^{**} \cup a_2^{**} \cup Q^{**} \cup a_2'^{**} \cup P'^{**}.$$

Then $Q_1^{**}$ is invariant under a lift $\alpha_2^*$ of $\alpha_2$ to $S_N$. It comes with a conformal structure and a projection $\pi^*$ into $\mathbb{S}^2$ which is a locally injective meromorphic function.

Construct the orbit of $Q_1^{**}$ under the group $\Gamma^{**}$ generated by the cover transformations of $P^{**}$ over $P$ and $P'^{**}$ over $P'$; $\Gamma^{**}$ is the free product of these groups with amalgamation over $\langle \alpha_2^* \rangle$. This can be done as follows. Suppose, for example, that $\alpha^* \notin \langle \alpha_2^* \rangle$ is a cover transformation of $P^{**}$ over $P$, so that $\alpha^*$ is the lift of a cover transformation $\alpha$ of $P^*$ over $P$. In particular, $\alpha^*$ sends the edge $a_2^{**}$ of $P^{**}$ to the edge $\alpha^* a_2^{**}$, which is invariant under the conjugate $\alpha^* \alpha_2^* \alpha^{*-1}$ of $\alpha_2^*$.

But the configuration $Q_1^{**}$ extends beyond $P^{**}$ at $a_2^{**}$. We correspondingly attach $\alpha^*(Q_1^{**})$ to extend beyond $\alpha^* a_2^{**}$. Moreover there is a projection $\pi^*$ of $Q_1^{**}$ into $\mathbb{S}^2$ which is a local homeomorphism, the extension of the restriction of $\pi^*$ to $P^{**}$. Extend $\pi^*$ from $P^{**}$ to $\alpha^*(Q_1^{**})$ by

$$\pi^*(z) = \alpha \pi^*(z_0), \quad z = \alpha^*(z_0), \quad z_0 \in Q_1^{**}.$$

The cover transformation $\gamma^*$ of $P^{**}$ or $P'^{**}$ over $P$ or $P'$ is conjugated to the cover transformation $\alpha^* \gamma^* \alpha^{*-1}$ of $\alpha^*(P^{**}) = P^{**}$ or $\alpha^*(P'^{**})$ over $P$ or $P'$. The transformation $\alpha^* \gamma^* \alpha^{*-1}$ itself is the lift of the cover transformation $\alpha \gamma \alpha^{-1}$ of $P^*$ over $P$ or of $\alpha(P'^*)$ over $P'$.

$Q_1^{**}$ and then $Q_1^{**} \cup \alpha^*(Q_1^{**})$ are simply connected Riemann surfaces that inherit their complex structure from $\mathbb{S}^2$ via $\pi^*$.

Continuing on, we construct a *pants configuration*

$$\mathfrak{J}(P, a_2; a_2', P')$$

which is a simply connected Riemann surface with a group of conformal automorphisms $\Gamma^{**}$. It has a meromorphic projection $\pi^*$ into (usually onto) $\mathbb{S}^2$, which is a local homeomorphism. The projection $\pi^*$ induces a homomorphism of $\Gamma^{**}$ onto the group generated by $\langle\alpha_1, \alpha_2\rangle$ and $\langle\alpha_2, \alpha_3\rangle$.

Consequently, with the group $\Gamma^{**}$, the abstract configuration

$$\mathcal{J}(P, a_2; a_2', P')$$

is a model for the universal covering of the Riemann surface

$$P \cup a_2 \cup (Q^{**}/\langle\alpha_2^*\rangle) \cup a_2' \cup P'.$$

It is a four-holed sphere; the pants $P$ and $P'$ have been connected by the cylinder $Q^{**}/\langle\alpha_2^*\rangle$, which joins $a_2$ and $a_2'$.

The Riemann mapping

$$g : \mathbb{H}^2 \to \mathcal{J}(P, a_2; a_2', P')$$

conjugates $\Gamma^{**}$ to a fuchsian group $G$ in $\mathbb{H}^2$. The function $f = \pi^* g : \mathbb{H}^2 \to \mathbb{S}^2$ is meromorphic and locally univalent in $\mathbb{H}^2$. It gives a projective structure on the four-holed sphere $\mathbb{H}^2/G$ with the associated homomorphism sending $G$ to the group generated by $\langle\alpha_1, \alpha_2\rangle$ and $\langle\alpha_2, \alpha_3\rangle$.

8.3. *Adding to the join of two pants.* At the level of the pants $P$ in $R$ and $P'$ in $R'$, the construction of §8.2 only involved neighborhoods of the boundary components $a_2$ of $P$ and $a_2'$ of $P'$, and the sector of $S_N$ between their two lifts.

Thus, suppose there is another Schottky group $\langle\alpha_3, \alpha_4\rangle$ sharing the generator $\alpha_3$ with $\langle\alpha_2, \alpha_3\rangle$. We can join the boundary component $a_3'$ of $P'$ to a compatibly chosen boundary component $a_3''$ of a pants $P''$ in $R'' = \Omega''/\langle\alpha_3, \alpha_4\rangle$, lying to the right of $a_3''$, by constructing the appropriate $S_N$. A lift of $P'^*$ appears in both configurations $\mathcal{J}(P', a_3'; a_3''; P'')$ and $\mathcal{J}(P, a_2; a_2', P')$, and these two lifts of $P'^*$ can be identified.

Join together these two configurations by identifying the two lifts $P'^{**}$ and $P_1'^{**}$ of $P'^*$. After that, further construct its orbit under $\Gamma_2^{**}$. Now $\Gamma_2^{**}$ is the free product of $\Gamma^{**}$ and the corresponding group $\Gamma_1^{**}$ of $\mathcal{J}(P', a_3'; a_3'', P'')$ with amalgamation over the common subgroup $\text{Stab}(P'^{**}) = \text{Stab}(P_1'^{**})$, which is just the lift of the covering group of $P'^*$ over $P'$. We end up with an abstract configuration

$$\mathcal{J}(P, a_2; a_2', P', a_3'; a_3'', P'') = \mathcal{J}_2,$$

which is a simply connected Riemann surface with a group $\Gamma_2^{**}$ of conformal automorphisms. There is a meromorphic projection $\pi^*$ into $\mathbb{S}^2$ that is a local homeomorphism and induces a homomorphism of $\Gamma_2^{**}$ onto the group generated by $\langle\alpha_1, \alpha_2\rangle$, $\langle\alpha_2, \alpha_3\rangle$ and $\langle\alpha_3, \alpha_4\rangle$. Also, $\mathcal{J}_2$ is a model of the universal cover for a five-holed sphere formed by connecting $P$ to $P'$ as in §7.5, and the result to $P''$ with an appropriate cylinder connecting $a_3'$ and $a_3''$.

8.4. *Making handles.* In §8.2, suppose that instead of a second Schottky group, we are presented with a transformation $\beta$ such that $\beta\alpha_1\beta^{-1} = \alpha_2$. We can as well join the group $\langle\alpha_1, \alpha_2\rangle$ to its conjugate $\beta\langle\alpha_1, \alpha_2\rangle\beta^{-1} = \langle\alpha_2, \beta\alpha_2\beta^{-1}\rangle$, to have the effect that the boundary component $a_2$ of the pants $P$ in $R$ is joined to $a_1$. We must start by ensuring that $a_2^*$ is compatible with $\beta a_1^*$ with respect to $\alpha_2$; this may require replacing $a_2$ by the result of applying some power of a Dehn twist about $b_2$.

As before, we can find an $S_N$ that holds disjoint lifts of $P^*$ and $\beta P^*$. Then a configuration $\mathcal{J}$ is constructed with a group of automorphisms $\Gamma^{**}$ isomorphic to the HNN extension of $\mathrm{Stab}(P^*)$ by a suitable lift $\beta^*$ of $\beta$. This $\mathcal{J}$ is a simply connected Riemann surface with a locally univalent meromorphic projection into $\mathbb{S}^2$. It is a model for the universal covering surface for the one-holed torus obtained by attaching the cylinder obtained from $S_N$ to the boundary components $a_1$ and $a_2$ of $P$.

8.5. Recall from §5.7 the trivalent graph $\mathbb{T}$ and the maximal tree $\mathbb{T}_0 \subset \mathbb{T}$. There, we chose one of the extreme vertices of $\mathbb{T}_0$ as the root. Let $\mathbb{T}_r$ denote the graph resulting from $\mathbb{T}$ after removing the one-edge loop hanging from the root. Thus $\mathbb{T}_r$ represents a surface $S \subset R$ of genus $g-1$ with two boundary components. Let $\Sigma_r$ denote the subgroup of $\pi_1(R; O)$ that is the fundamental group of $S$.

LEMMA 8.5.1. *There exists a pants configuration $\mathcal{J}(\mathbb{T}_r)$ modeled on $\mathbb{T}_r$. It is a simply connected Riemann surface, the universal cover of a Riemann surface $S$ of genus $g-1$ with two boundary components. Let $g : \mathbb{H}^2 \to \mathcal{J}(\mathbb{T}_r)$ be a Riemann mapping, and $\pi : \mathcal{J}(\mathbb{T}_r) \to \mathbb{S}^2$ the meromorphic projection. Then $f = \pi g$ is a projective structure for $S$ for the homomorphism $\theta : \Sigma_r \to \theta(\Sigma_r) \subset \Gamma$.*

*Proof.* First we check that the compatibility conditions can be satisfied. Denote by $P(v)$ the pants corresponding to the vertex $v$, and by $\Gamma(v)$ the Schottky group with regular set $\Omega(v)$.

In §5.8 we marked the boundary components of $P(v)$ according to the following rule. There is a unique path in $\mathbb{T}_r$ from any vertex $v'$ to the root $v_0$. The unmarked boundary component $a$ of $P(v)$ is the one on the path from $v$ itself. If $v$ is not an extreme point of $\mathbb{T}_0$, it has two immediate predecessors $v_1$ and $v_2$, and $P(v)$ has two marked boundary components $b$ and $c$, lying on their paths to $v_0$. Following the notation of §7.5, let $x$ and $y$ denote compressing loops (which lift to simple loops in $\Omega(v)$) such that $x$ crosses $c$ and $a$, and $y$ crosses $b$ and $a$.

Now move down the tree $\mathbb{T}_0$. Start at the extreme vertices $v \neq v_0$. Two of the boundary components $b$ and $c$ of $P(v)$ are to be paired. Make them compatible by a twist about either $x$ or $y$.

Continue down the tree. Do not go to a vertex before dealing with all its predecessors. Arriving at a vertex $v$ and $P(v)$ with marked borders $b$ and $c$, replace them by the result of twists about $x$ and $y$, so as to be compatible with the (unmarked) borders $b'$ and $c'$ associated with the immediate predecessors $v_1$ and $v_2$. When the root $v_0$ is reached, the one marked border of $P(v_0)$ is made compatible with its immediate predecessor.

Finally, use the technique illustrated in §§8.2–8.4 to join the pants $\{P_i\}$ together with auxiliary cylinders to build a Riemann surface of genus $g - 1$ with two boundary components remaining from the pants $P(v_0)$. This is done by building a pants configuration $\mathcal{J}(\mathbb{T}_r)$, which is its universal cover.  □

8.6  *The final handle or the two-sheeted covering.* Having constructed $\mathcal{J}(\mathbb{T}_r)$, all attention is focused on $P(v_0)$, with its three boundary components $a, b, c$ and compressing loops $x, y$ as in §7.5. Since $P(v_0)$ has been attached to its predecessor, say by establishing the compatibility of $c$ with its partner $c'$, no more twisting about $x$ is possible. Can we make $a$ compatible with $b$, allowing attachment of the final handle? As we have seen in §7.5, this is possible if and only if one can do the job with an even number of twists. If so, we can finish the construction of $\mathcal{J}(\mathbb{T})$, the pants configuration reflecting the full trivalent graph $\mathbb{T}$, which will then be a simply connected Riemann surface with a group of conformal automorphisms making it the universal cover of a surface of genus $g$.

If not, keeping in mind the alternate construction of §7.6, we will construct instead a pants configuration $\mathcal{J}$ that models a two-sheeted unbranched covering of the reference surface $R$.

Suppose $a$ and $b$ have arisen from cutting $R$ along a curve $b_0'$, freely homotopic to the nondividing simple loop $b_0 \in \pi_1(R; O)$. Set $R_0 = R \setminus \{b_0'\}$, and find the simple loop $a_0 \in \pi_1(R; O)$ such that $b_0$ and $a_0 b_0^{-1} a_0^{-1}$ give rise to $\pi_1(P(v_0); O)$. The group

$$N = \langle a_0^2,\ \pi_1(R_0; O), a_0 \pi_1(R_0; O) a_0^{-1} \rangle$$

is a normal subgroup of index two in $\pi_1(R; O)$. It defines a two-sheeted unbranched covering $\tilde{R}$ of $R$ that is a compact surface of genus $2g - 1$.

The surface $\tilde{R}$ is explicitly constructed as follows. Label the boundary components of $R_0$ as $b_0^+$ and $b_0^-$, corresponding to the two sides of $b_0'$ in $R$. Take another copy $R_0'$ of $R_0$. Then $\tilde{R}$ is the surface obtained by identifying $b_0^+$ and $b_0^-$ on $R_0'$ with $b_0^-$ and $b_0^+$, respectively on $R_0$. The cover transformation is determined by $a_0$.

Let $\mathbb{T}_2$ denote the trivalent graph built likewise by taking two copies of $\mathbb{T}_r$ and attaching two new edges $e_1$ and $e_2$, as follows. The endpoints of the new edges are the two vertices corresponding to $v_0$ (and pants $P(v_0)$), and they serve to pair the boundary components $a$ and $b$ on one copy of $P(v_0)$ with $b$ and $a$, respectively, on the other.

Correspondingly, take two copies of $\mathcal{J}(\mathbb{T}_r)$. Because of the compatibility established in §7.6, they can be joined together following the combinatorics of $\mathbb{T}_2$ and the restriction of $\theta$ to $N$. The resulting pants configuration $\mathcal{J}(\mathbb{T}_2)$ is again a simply connected Riemann surface with a group of conformal automorphisms isomorphic to $N$, making it the universal cover of a surface of genus $2g - 1$.

Because of the asymmetry in satisfying the compatibility for the two copies of $P(v_0)$ (see §7.6), $\mathcal{J}(\mathbb{T}_2)$ does not have conformal automorphisms that represent the sheet interchange of $\tilde{R}$. If, however, $\mathcal{J}(\mathbb{T})$ can be constructed, and then $\mathcal{J}(\mathbb{T}_2)$ constructed in addition, $\mathcal{J}(\mathbb{T}_2)$ will have that symmetry: it will represent the universal cover of the two-sheeted cover of the Riemann surface corresponding to $\mathcal{J}(\mathbb{T})$.

A Riemann mapping $g : \mathbb{H}^2 \to \mathcal{J}(\mathbb{T})$ or $g_2 : \mathbb{H}^2 \to \mathcal{J}(\mathbb{T}_2)$ conjugates the cover transformations to a fuchsian group $G$ isomorphic to $\pi_1(R; O)$ or to a fuchsian group $G_2$ isomorphic to the index two subgroup $N$. Let $\pi$ denote the projection of $\mathcal{J}(\mathbb{T})$ or $\mathcal{J}(\mathbb{T}_2)$ to $\mathbb{S}^2$. The meromorphic function $f = \pi g$ or $f_2 = \pi g_2$ determines a projective structure that induces $\theta : G \to \Gamma$ or $\theta : G_2 \to \theta(G_2) \subset \Gamma$.

We cannot exclude the possibility that $\theta(G_2) = \Gamma$. Although the transformation in $\Gamma$ that makes the conjugation corresponding to the pairing of the boundary components $a$, $b$ of $P(v_0)$ is not the identity (because $P(v_0)$ arises from a two-generator Schottky group), it may already belong to $\theta(G_2)$. In any case, if $\theta : \pi_1(R) \to \Gamma$ cannot be lifted to $\mathrm{SL}(2, \mathbb{C})$, $\theta : N \to \theta(N)$ can be so lifted.

## 9. Attaching branched disks to pants

9.1. One can attach a disk to any surface with boundary by introducing a single branch point. Explicitly for our situation, consider a pants $P$ embedded in $\mathbb{C}$ and a boundary component $a$ oriented so that $P$ lies to its right. Suppose $d$ is an oriented simple loop bounding a disk $\Delta$ lying to its right. Suppose that $d$ crosses $a$ at a point $p$, and that $z_1$ and $z_2$ are points separated by both $a$ and $d$, with $z_1 \in P \cap \Delta$. Assume that there exists a simple arc $\sigma$ between $z_1$ and $z_2$ that crosses both loops at $p$ and is otherwise disjoint from them. Set $\sigma_0 = \sigma \cap P \cap \Delta$.

Attach the $\Delta$ to $P$ as follows. Denote the opposite sides of $\sigma_0$ by $\sigma_0^+$ and $\sigma_0^-$. Identify the side $\sigma_0^+$ of $\Delta \setminus \sigma_0$ with the side $\sigma_0^-$ of $P \setminus \sigma_0$, and the side $\sigma_0^-$ of $\Delta \setminus \sigma_0$ with the side $\sigma_0^+$ of $P \setminus \sigma_0$. This determines a new Riemann surface $P'$ that is conformally equivalent to a new pants. Its boundary $\partial P'$ consists of $a \cup d$ (here $d$ lies "over" $P$) and the remaining components of $\partial P$. The natural holomorphic projection $\pi : P' \to P \cup \Delta$ is a local homeomorphism except at the point over $z_1$, where it behaves like $z \mapsto z^2$. See Figure 10.

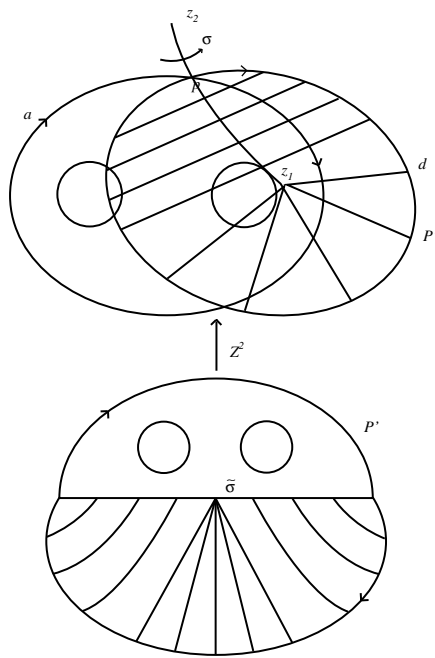

Figure 10. Attachment of branched disk

Note that the construction does not essentially depend on a choice for $\sigma$. Instead we can work in the two-sheeted cover of $\mathbb{S}^2$, branched over $z_1$ and $z_2$

The same construction can be applied to attach an $(n-1)$-sheeted disk to $P$, for any $n \geq 2$.

9.2. *Application to pants in a Schottky group.* Suppose that $\langle \alpha, \beta \rangle$ is a Schottky group acting on $\Omega \subset \mathbb{C}$, and $P \subset \Omega/\langle \alpha, \beta \rangle$ is a pants with boundary components $a, b, c$ oriented so that $P$ lies to the right of $a$ and $b$, which have $\alpha$- and $\beta$-invariant lifts $a^*$ and $b^*$ in $\Omega$. Let $d$ be a compressing curve on the handlebody surface $\Omega/\langle \alpha, \beta \rangle$ that crosses $a$ exactly once, at a point $p$.

Introduce a simple arc $\sigma$ in $\Omega/\langle \alpha, \beta \rangle$ that joins a point $z_1 \in P$ to $z_2$ in its complement, and crosses the loops $a$ and $d$ at $p$, otherwise being disjoint from them. Set $\sigma_0 = \sigma \cap P$.

Let $d^*$ be a simple loop in $\Omega$ lying over $d$, which crosses $a^*$ (necessarily once). Orient $d$ and thus $d^*$ so that the disk $\Delta$ lying to its right contains the lift of $\sigma_0$ that is adjacent to $d^*$.

Attach $\Delta$ to $P$ by means of the slit $\sigma_0$. Neither the resulting pants $P_1$ is embedded in $\mathbb{C}$ nor $\Delta$ is embedded in $\Omega/\langle \alpha, \beta \rangle$. Nevertheless, any annular neighborhood of $d$ in $\Omega/\langle \alpha, \beta \rangle$ is conformally equivalent to its lift about $d^*$. Thus the conformal structure of $P_1$ is well defined.

Equivalently, the universal cover $\tilde{P}$ of $P$ is embedded in $\Omega$, and the universal cover $\tilde{P}_1$ of $P_1$ arises from that by attaching $\Delta$ by means of the lift of $\sigma_0$ that is adjacent to $d^*$, and then taking the orbit under $\langle \alpha \rangle$ of the attachment. We need to examine this construction more closely.

The attachment of $\Delta$ to $a^*$ at $p^* \in a^*$ over $p$ leads to the attachment of the loop $d^*$ to $a^*$ at $p^*$: as we move along $a^*$ toward the attracting fixed point of $\alpha$, when we reach $p^*$ we take a detour along $d^*$ in its positive direction, returning to $p^*$ and then continuing along $a^*$. Since $d^*$ intersects $a^*$ only at $p^*$, the resulting arc is essentially a simple arc, and so is its $\langle \alpha \rangle$-orbit, which covers the point set $a^* \cup \alpha^k(d^*)$.

The essentially simple, $\alpha$-invariant arc $a^* \cup \alpha^k(d^*)$ can equally be described as follows. It is the lift of the result of applying to $a$ on $\Omega/\langle \alpha, \beta \rangle$, or its representation in the torus $T(\alpha)$, a Dehn twist about $d$.

9.3. *Another alternative to the geometric obstruction of Section* 8.6. In §8.6 we faced the question of adding the final handle to the pants configuration $\mathcal{J}(\mathbb{T}_r)$. If that was not possible, we showed that we could instead construct a pants configuration corresponding to a two-sheeted, unbranched cover of the surface of genus $g$.

Alternatively, using the construction of §9.2, we can carry out the final construction after introducing a branch point of order two (or any even order). That is, we can construct a pants configuration $\mathcal{J}_b(\mathbb{T})$ representing the universal covering of a Riemann surface of genus $g$. If $g : \mathbb{H}^2 \to \mathcal{J}_b(\mathbb{T})$ is a Riemann map, and $\pi : \mathcal{J}_b(\mathbb{T}) \to \mathbb{S}^2$ is the natural projection, then $f = \pi \circ g$ is a meromorphic function. It is locally injective except at the conjugacy class of branch points of order two, and still induces the homomorphism $\theta : \pi_1(R; O) \to \Gamma$.

## 10. The obstructions

10.1. *The modulo* 2 *construction invariant.* An *admissible* pants decomposition $\{P_i\}$ for the homomorphism $\theta : \pi_1(R; O) \to \Gamma$ is one for which the restriction of $\theta$ sends each $\pi_1(P_i)$ to a Schottky group. Its combinatorics are associated with a trivalent graph $\mathbb{T}$. To each vertex $v$ of $\mathbb{T}$ is associated a Schottky group $S(v) = \langle \alpha_v, \beta_v \rangle$ acting on $\Omega(v) \subset \mathbb{S}^2$. To each $S(v)$ is associated a pants $P(v) \subset \Omega(v)/S(v)$ with boundary components $a, b, c$ that have $\alpha_v$-, $\beta_v$- and $\beta_v\alpha_v$-invariant lifts in $\Omega(v)$. In terms of corresponding elements of $\pi_1(P(v))$, we have $c' \sim b'a'$ in $\Omega(v)/S(v)$. The orientation of $P(v)$ with respect to $a$ and $b$, and hence for $c$, has been dictated by that of the corresponding $P_i$ with respect to its boundary components and carried over to $\mathbb{T}$ by $\theta$.

Each edge $e$ of $\mathbb{T}$ corresponds to a common generator $\alpha$ of the two Schottky groups $S(v_1)$ and $S(v_2)$ if the endpoints of $e$ lie on $v_2 \neq v_1$. If $v_2 = v_1$, then $e$ is associated with a pair of boundary components of $P(v_1)$, which in turn correspond to generators $\alpha_v$ and $\beta_v$ related by $\beta_v = \gamma_v \alpha_v^{-1} \gamma_v^{-1}$ for some element $\gamma_v \in \Gamma$. In any case the pair of boundary components corresponding to $\alpha$ project to a pair of simple loops on the torus $T(\alpha)$. The two boundary components are called *compatible* if their projections, appropriately oriented, are freely homotopic on $T(\alpha)$.

We will call $\mathbb{T}$ *compatible* if all pairs of boundary components of the associated pants $\{P(v)\}$ are compatible.

Recall that on each torus $T(\alpha)$ there is a free homotopy class of simple loops called compressing loops (§6.1), each of which lifts to a simple loop in $\mathbb{S}^2$.

LEMMA 10.1.1. *Suppose that on each $T(\alpha)$ one of the boundary projections is freely homotopic to the result of a Dehn twist of order $n(\alpha)$ (about a compressing loop) applied to the other. Set $n(\mathbb{T}) = \left( \sum_\alpha n(\alpha) \right) \mathrm{mod}\, 2$. There is a compatible pants decomposition $\{P(v)\}$ corresponding to $\mathbb{T}$ if and only if $n(\mathbb{T}) = 0$.*

*Proof.* To each pair of pants one can apply Dehn twists about compressing loops on $\Omega(v)/S(v)$. The algebraic sum $n(P(v))$ of their effect on the three boundary components of $P(v)$ is an even number. Thus

$$\sum n(P(v)) = 0 \pmod 2.$$

Hence the values of $n(\mathbb{T})$ cannot be changed by repositioning the pants $P(v)$ in the surfaces $\Omega(v)/S(v)$.

For the graph $\mathbb{T}$ of §8.4 that represents the "localization" of the obstruction to the construction, the question of compatibility rested on the compatibility of the two paired boundary components in the root pants $P(v_0)$ (§8.6). This was precisely the question of whether or not $n(\mathbb{T}) = 0$. That is, if $n(\mathbb{T}) = 0$ we can distribute the twists so that $\mathbb{T}$ is compatible.

For other graphs $\mathbb{T}$, we refer to Corollary 10.5.1. $\qquad \square$

10.2. *Lifting Schottky groups. Lifting* refers to the property that a given homomorphism $\theta \colon \pi_1(R; O) \to \mathrm{PSL}(2, \mathbb{C})$ lifts to a homomorphism $\theta^* \colon \pi_1(R; O) \to \mathrm{SL}(2, \mathbb{C})$. The image groups are not necessarily isomorphic.

It is helpful to recall the case where $H = \langle \alpha, \beta \rangle$ is a two-generator, purely loxodromic fuchsian group. As such it represents either a handle or a pants. Let $A$ and $B$ be matrix representations of $\alpha$ and $\beta$. Then $H$ is isomorphic to $\langle A, B \rangle$. The commutator matrix $[A, B]$ is independent of the choice of lift of $\alpha$ and $\beta$. The two cases, handle or pants, can be distinguished according to whether $[\alpha, \beta]$ represents a simple loop or not, or whether no axis in its conjugacy class separates the axes of $\alpha$ and $\beta$ or does, or whether $\mathrm{tr}[A, B] < -2$ or $\mathrm{tr}[A, B] > 2$. Moreover, in the case of a handle, the free homotopy class in the torus $T([\alpha, \beta])$ determined by a loop parallel to the handle boundary is uniquely determined, independent of Dehn twists about compressing loops when regarding $\langle \alpha, \beta \rangle$ as a Schottky group.

More generally, any Schottky group $\langle \alpha, \beta \rangle$ can be lifted to an isomorphic group in $\mathrm{SL}(2, \mathbb{C})$ by designating matrix representatives for $\alpha$ and $\beta$.

10.3. *The modulo 2 lifting obstruction.* Let $\mathbb{T}$ be a trivalent graph as in §10.1. Lift to $\mathrm{SL}(2, \mathbb{C})$ the Schottky groups corresponding to its vertices.

Let $e$ be an edge of $\mathbb{T}$ with endpoints $v_1$ and $v_2$. If $v_1 \neq v_2$, the edge $e$ corresponds to a common generator $\alpha$ of $S(v_1)$ and $S(v_2)$. The lifting will be called *compatible* on $e$ if the lifted $\alpha$ in $S(v_1)$ and lifted $\alpha$ in $S(v_2)$ have the same trace. If $v_1 = v_2$, the compatibility condition is that the designated lifts of $\alpha$ and $\gamma \alpha \gamma^{-1}$ from $S(v_1)$ have the same trace. The lifting of $\mathbb{T}$ will be called compatible if it is compatible on each edge.

Suppose $\mathbb{T}$ is the graph of §8.4 with its maximal tree $\mathbb{T}_0$. Start at the extreme vertices of $\mathbb{T}_0$ and work down towards the root: Exactly in analogy to the construction of §8.4, choose at each step a lift of a generator of a Schottky group to be compatible with the lifts previously chosen. We end up with a compatible lift of $\mathbb{T}_r$. The lift of $\mathbb{T}_r$ is determined by the two choices made at the $g - 1$ extreme points of $\mathbb{T}_0$ other than the root, and one choice at the root.

LEMMA 10.3.1. *Suppose $\mathbb{T}$ is the trivalent graph corresponding to an admissible pants decomposition. Then $\mathbb{T}$ has a compatible lift to $\mathrm{SL}(2, \mathbb{C})$ if and only if the homomorphism $\theta$ can be lifted to $\mathrm{SL}(2, \mathbb{C})$.*

*Proof.* The graph $\mathbb{T}$ corresponds to a presentation of $\pi_1(R)$.     □

10.4. *Localization of the lifting obstruction.* Denote by $\langle \alpha^{-1} \beta \alpha, \beta^{-1} \rangle$ the Schottky group corresponding to the root. We recall from §5.7 that the "handle group" $H = \langle \alpha, \beta \rangle$ is nonelementary with $\alpha$ and $\beta$ loxodromic, even though it may not be discrete.

Applying the technique of §8, we can build a pants configuration $\mathcal{J}_h$ on which $H$ acts so that $\mathcal{J}_h/H$ is a handle. Likewise the graph $\mathbb{T}'_h$ resulting from removing from $\mathbb{T}$ the root and attached edges determines a pants configuration $\mathcal{J}'_h$ acted on by a group $H'$ so that $\mathcal{J}'_h/H'$ is a surface of genus $g - 1$ with one boundary component.

Choose matrix representatives $A$ and $B$ for $\alpha$ and $\beta$; then $[B, A]$ is a representative for $[\beta, \alpha]$, which corresponds to the boundary component of the handle.

The graph $\mathbb{T}'_h$ can be lifted to $\mathrm{SL}(2, \mathbb{C})$ as in §10.4, which yields a matrix $C$ representing $[\beta, \alpha] \in H'$. Therefore $C = \pm[B, A]$.

LEMMA 10.4.1. *The homomorphism $\theta$ lifts to $\mathrm{SL}(2,\mathbb{C})$ if and only if $C = [B, A]$. In particular, $\theta$ lifts if $\mathfrak{J}_h$ and $\mathfrak{J}'_h$ can be joined to form a pants configuration for $\mathbb{T}$.*

*Proof.* The first assertion follows from Lemma 10.3.1. The second follows as a consequence of the existence of a projective structure (see, for example, Lemma 1.3.1). ∎

10.5. *Equivalence of obstructions.*

PROPOSITION 10.5.1. *The procedure of §8 succeeds in constructing a projective structure associated with the given homomorphism $\theta : \pi_1(R; O) \to \mathrm{PSL}(2,\mathbb{C})$ if and only if $\theta$ can be lifted to a homomorphism into $\mathrm{SL}(2,\mathbb{C})$.*

*Proof.* From §1.3 we already know lifting is a necessary condition. Now suppose $\theta$ can be lifted, yet the construction cannot be completed. That is, in the notation of §10.4, $\mathfrak{J}_h$ cannot be attached to $\mathfrak{J}'_h$. But then, as in §9, we can introduce a single branch point of order two and construct instead a branched projective structure associated with $\theta$. According to §1.4, $\theta$ cannot be then lifted to $\mathrm{SL}(2,\mathbb{C})$, in contradiction with the assumption. ∎

COROLLARY 10.5.2. *If the construction of a projective structure works for one admissible pants decomposition for $\theta$, it works for any admissible decomposition.*

## C. Ramifications

## 11. Holomorphic bundles over Riemann surfaces, the 2nd Stiefel-Whitney class, and branched complex projective structures

The purpose of this chapter is to place Theorem 1.1.1 in a more general setting, and to use that to clarify the role played by branched structures in Part B. We will also discuss relations between instability of holomorphic vector bundles over Riemann surfaces and branched complex projective structures. In §11.5 we establish the local character of the map between singly branched structures over Teichmüller space and the representation variety. In §11.6, we again use holomorphic vector bundles to prove that for singly branched structures too the monodromy representation is necessarily nonelementary.

11.1. *The 2nd Stiefel-Whitney class of sphere bundles over Riemann surfaces.* Suppose that $\eta : P \to R$ is a holomorphic $\mathbb{CP}^1$-bundle over a closed Riemann surface $R$. It is known (see for instance [Beau, Prop. III.7]) that $P$ can be obtained as the projectivization of a holomorphic (rank 2) vector bundle $\xi : V \to R$. Let $\det(V)$ denote the determinant bundle of $V$, this is a holomorphic line bundle over the surface $R$. The bundle $V$ is not uniquely determined by the projective bundle $P \to R$, and to obtain an isomorphic projective bundle, we can alter $V$ by multiplying it by a holomorphic line bundle $\Lambda$ over $R$. Then $\deg(\det(V \otimes \Lambda)) = \deg(\det(V)) + 2\deg(\Lambda)$. Thus we can always choose $V$ so that $\det(V)$ has degree 0 or 1.

Let $p : V \to P(V) = P$ be the projectivization. We shall think of $p$ as a holomorphic line bundle over the base $P$. It is well-known that there are exactly two topologically distinct orientable $\mathbb{S}^2$-bundles over the surface $R$ (see [Mel]) and they are distinguished by the 2nd Stiefel-Whitney class $w_2(P)$ of the bundle $P \to R$.

Note that if $\deg(\det(V)) = 0$ then the determinant bundle $\det(V)$ is topologically trivial. In this case the bundle $V$ is associated to an $SL(2, \mathbb{C})$-bundle over $R$ which is henceforth topologically trivial. We conclude that $w_2(P)$ equals $\deg(\det(V)) \pmod 2$.

Let $\sigma : R \to P(V)$ be a holomorphic section of $P(V)$. It defines a holomorphic line bundle $L \to R$ by pull-back $\sigma^*(p)$ of the line bundle $p$. The line bundle $L$ is canonically embedded as a holomorphic subbundle of the bundle $\xi : V \to R$ with the image $p^{-1}(\sigma(R))$.

LEMMA 11.1.1. (1) $\sigma^2 = \deg(\det(V)) - 2\deg(L)$, where the left-hand side is the self-intersection number of the cycle $\sigma(R)$ in $P(V)$. (2) The number $\sigma^2 \pmod 2$ equals the 2nd Stiefel-Whitney class $w_2(P)$ of the bundle $\eta : P \to R$.

*Proof.* The first assertion is a particular case of a general result proven in [La, §1]. Since $w_2(P)$ equals $\deg(\det(V)) \pmod 2$, the second assertion follows.

Nevertheless we will provide a elementary proof of the first assertion for the sake of completeness. We first consider the case $\deg(\det(V)) = 0$ and then we shall reduce the general case to this one. If $\deg(\det(V)) = 0$ then both bundles $V$ and $P$ are topologically trivial. Hence there is an orientation preserving diffeomorphism $P(V) \to R \times F$, where $F = \mathbb{S}^2$. By the Künneth formula, the homology class $[\sigma]$ can be written as

$$[\sigma] = n[F] + [R],$$

and we get: $\sigma^2 = 2n$. There are two possible cases: $n \geq 0$ (if $\sigma^2 \geq 0$) and $n < 0$ (if $\sigma^2 < 0$). We consider the former; the later case is analogous (one just has to work with anti-holomorphic functions instead of the holomorphic ones). We can think of $\sigma : R \to R \times F$ as a graph of a smooth function $f : R \to F = \mathbb{S}^2$ which has nonnegative degree $n$. The function $f$ is not holomorphic, however (after deforming the section $\sigma$ within its homotopy class) we can assume that $f^{-1}(\infty) = Z := \{z_1, ..., z_n\} \subset R$ and $f$ is holomorphic near each point $z_j$ so that $f'(z_j) \neq 0$, $1 \leq j \leq n$. Now we realize $F = \mathbb{C} \cup \{\infty\}$ as the complex projective line $\mathbb{CP}^1$ so that the point $\infty$ has the homogeneous coordinates $[1:0]$. Then we lift the function $f$ to the meromorphic function

$$\tilde{f} : R \to \mathbb{C}^2, \tilde{f}(z) = (f(z), 1)$$

which does not have zeroes and is holomorphic in a punctured neighborhood of each point $z_j \in Z$ and has a simple pole at each $z_j \in Z$. Thus $\tilde{f}$ corresponds to a smooth meromorphic section of the line bundle $L \subset V$ which has $n$ simple poles and no zeroes. Hence $\deg(L) = -n = -\sigma^2/2$.

Now we consider the case when $\deg(\det(V)) = 2k$ is an even number. Take a complex line bundle $\Lambda$ over $R$ so that $\deg(\Lambda) = -k$, then $\deg(\det(\Lambda \otimes V)) = 0$. The section $\sigma : R \to P$ defines complex line subbundle of $\Lambda \otimes V$ which is isomorphic to $\Lambda \otimes L$. As we proved above, $\sigma^2 = \deg(\det(\Lambda \otimes V)) - 2 \deg(\Lambda \otimes L)$ which in turn equals to $\deg(\det(V)) - 2 \deg(L)$. This completes the proof in the case when $\deg(\det(V))$ is even.

In the case when $\deg(\det(V))$ is odd take a 2-fold unramified covering $\tilde{R} \to R$. Then the bundle $V \to R$ pulls back to a bundle $\tilde{V} \to R$ and $\deg(\det(\tilde{V})) = 2 \deg(\det(V))$ is even. Similarly, the section $\sigma$ determines a section $\tilde{\sigma} : \tilde{R} \to P(\tilde{V})$ and $\tilde{\sigma}^2 = 2\sigma^2$. The pull-back of the line bundle $L$ to $\tilde{L} \subset \tilde{V}$ has degree equal to $2 \deg(L)$. We get:

$$\tilde{\sigma}^2 = \deg(\det(\tilde{V})) - 2 \deg(\tilde{L})$$

which implies

$$\sigma^2 = \deg(\det(V)) - 2 \deg(L).$$

This concludes the proof in the general case. □

11.2. *Branched structures.* Consider a Riemann surface $R = \Omega/\pi_1(R)$ where $\Omega$ is the universal cover $\tilde{R}$ of $R$ and is either the unit disk, or the complex plane, or the Riemann sphere and the group $\pi_1(R)$ of Möbius transformations acts freely and discontinuously on $\Omega$.

Suppose that $\theta : \pi_1(R) \to \Gamma \subset \mathrm{PSL}(2, \mathbb{C})$ is a homomorphism, and $f : \Omega \to f(\Omega) \subseteq \mathbb{S}^2$ is a meromorphic function (without essential singularities) which is $\theta$-equivariant and defines a *branched (complex) projective structure* $\sigma$ on $R$ as in §1.4. Alternatively one can define a branched projective structure on $R$ as a collection of locally defined holomorphic (but not necessarily univalent) mappings $\phi_\alpha$ from $R$ to $\mathbb{S}^2$ so that different mappings are related by Möbius transformations $\gamma_{\alpha,\beta}$:

$$\phi_\alpha = \gamma_{\alpha,\beta} \circ \phi_\beta$$

(see for instance [Man1]).

The homomorphism $\theta$ is the (*projective*) *monodromy representation* of the branched projective structure, and in the terminology of §1.3 the projection $f_* : R \to \mathbb{S}^2$ is the *(multivalued) developing map*. We define the *branching divisor* $D_f$ as follows. Consider the discrete set $\tilde{D}_f \subset \tilde{R}$ consisting of critical points of $f$. Thus (after holomorphic change of variables), near such a critical point $z_j$ the function $f(z)$ can be written as

$$f(z) = z^k, \quad 2 \le k < \infty.$$

The number $k$ is the *order* of branch point $z_j$. Since the function $f$ is $\theta$-equivariant we conclude that for any $\gamma \in \pi_1(R)$ the point $\gamma(z_j)$ is again a branch point with the same order $k$. Hence the projection of $\tilde{D}_f$ to the surface $R$ is a finite collection of points, to each such point $w_j$ we have the associated the number $\mathrm{ord}(w_j) = k_j > 1$ which is its order. Define the (additive) *branching divisor* $D = D_f$ of the structure $\sigma$ as

$$\sum_{w_j} (k_j - 1) w_j.$$

The number

$$d = \sum_{w_j} (k_j - 1) = \deg(D_f) \ge 0$$

is *the degree* of this divisor. The number $k_j - 1$ is the *local degree* $\deg_{w_j}(D)$ of the divisor $D$ at the point $w_j$. The *multiplicity* $|D|$ of the divisor $D$ is just the number of points in it. If $\deg(D) = 0$, the divisor $D$ is empty and there is no branching.

For reasons that we shall see later on, it is convenient to define the divisor $D$ by subtracting 1 from the order of each branch point. In addition we will consider the branching divisor as a *topological* object, not an analytic one. Thus we will say that two branching divisors $D, D'$ on $R$ are *equivalent* if there exists a bijective order-preserving map $D \to D'$ between them. This is the only meaningful equivalence relation in our situation since we will have to change the complex structure on $R$ in order to find a branched projective structure with the prescribed monodromy.

Next we review the relation between branched projective structures and Schwarzian differential equations as in §1.4. Let $D$ be a positive divisor on the Riemann surface $R$. Suppose that $\phi(z)dz^2$ is a meromorphic quadratic differential on $R$ which is holomorphic on $R - D$ and near each point $w_j \in D$ has a Laurent expansion of the form

$$(7) \qquad \phi(z) = \frac{(1 - k_j^2)}{2z^2} + \frac{b}{z} + \sum_{i=0}^{\infty} a_i z^i.$$

Here we use local coordinates such that $w_j = 0$ and $k_j - 1 = \deg_{w_j} D$ is the local degree of $D$. If

$$(8) \qquad f(z) = z^{k_j} h(z)$$

where $h(z)$ is a holomorphic function such that $h'(0) \neq 0$, then the Schwarzian derivative $S_f(z)$ near zero has Laurent expansion of the form (7). Conversely, to have a solution in the form (8) the quadratic differential $\phi(z)dz^2$ must satisfy an extra condition of *integrability*; see [He1] or [Man2].

Let $Q_D(R)$ denote the space of meromorphic quadratic differentials on $R$ with at most simple poles at points of $D$. Suppose that $\psi_0$ is a fixed quadratic differential of the form (7), then all other such quadratic differentials can be written as $\phi = \phi_0 + \psi$, where $\psi \in Q_D(R)$. Let $n$ denote the multiplicity of $D$. There exists a collection of $n$ polynomials $K_j$ on the $(3g-3+n)$-dimensional complex vector space $Q_D$ so that $\phi$ is integrable if and only if the differential $\psi$ belongs to the zero set of all the polynomials $K_j$. If $\deg(D) \leq 2g-2$ then the algebraic variety

$$I(R,D) := \{K_j(\psi) = 0, j = 1, ..., n\}$$

has generic dimension $3g-3$. In the case of a single-order two branch point at the orbit of $z = 0 \in \mathbb{H}^2$, $I(R,D)$ is given by the polynomial equation

$$(9) \qquad u^2 + 2bu + 2v = 0$$

where $u$ is the coefficient of the $z^{-1}$ term and $v$ is the constant term in the Laurent expansion of $\psi$ at $z = 0$. The number $b$ is given by §1.4(6). We refer to [Man1, 2, 3] for more details.

Now we go back to the linear differential equation

$$(10) \qquad u'' + \tfrac{1}{2}\phi u = 0$$

expressed in a local coordinate system on the surface $R$. With $\phi \in Q_D(R) + \phi_0$ and satisfying the integrability condition as above, the equation (9) has two linearly independent solutions. If $z_j$ is a singular point of $\phi$ and we choose local coordinates so that $z_j = 0$, near this point these solutions have the form

$$\begin{cases} u_1(z) = z^{(1+k_j)/2}(1 + o(1)) \\ u_2(z) = z^{(1-k_j)/2}(1 + o(1)). \end{cases}$$

A circuit about $z = 0$ generates the linear monodromy

$$\begin{pmatrix} u_1 \\ u_2 \end{pmatrix} \mapsto J^{k_j-1} \begin{pmatrix} u_1 \\ u_2 \end{pmatrix}, \qquad \text{where} \quad J = \begin{pmatrix} -1 & 0 \\ 0 & -1 \end{pmatrix}.$$

The projectivization of this monodromy in $\mathrm{PSL}(2,\mathbb{C})$ is just the identity.

LEMMA 11.2.1. *On the surface $R-D$ with a base-point $O$, the differential equation* (9) *has a linear monodromy representation*

$$\theta^* : \pi_1(R-D, O) \to \mathrm{SL}(2, \mathbb{C}).$$

*Proof.* This is a consequence of the fact that the Wronskian of two solutions is a constant (see Corollary 1.3.1). □

Let $U \subset R$ be a closed disc which contains all the singular points $z_j \in D$ and fix a base point $O \in \partial U = \ell$. The matrix $\theta^*(\ell)$ that results from analytic continuation along $\ell$ equals $J^d$ where $d = \deg(D)$ is the degree of this divisor and $J = -1$. The representation $\theta^*$ projects to a homomorphism $\theta : \pi_1(R) \to \mathrm{PSL}(2, \mathbb{C})$. We conclude that $\theta$ can be lifted to a linear representation

$$\tilde{\theta} : \pi_1(R) \to \mathrm{SL}(2, \mathbb{C})$$

if and only if the number $d$ is even, in particular if $d = 0$ as in Chapter 1. It is instructive to see a topological proof of this fact as well.

Let $P$ denote the $\mathbb{S}^2$-bundle over $R$ associated with the monodromy representation $\theta$ of a complex projective structure $\tau$ on $R$. It carries a natural flat connection. Let $w_2(\theta) := w_2(P)$. The developing map $f$ of the structure $\tau$ defines a holomorphic section $\sigma$ of the holomorphic bundle $P \to R$. We will treat $\sigma$ as a 2-cycle in $P$.

PROPOSITION 11.2.2.  *Under the above conditions we have*:

$$\langle \sigma(R), \sigma(R) \rangle = 2 - 2g + \deg(D),$$

*where $\langle \cdot, \cdot \rangle$ is the intersection pairing on the 4-manifold $P$.*

*Proof.* Note that the polynomial $z^n$ admits arbitrarily small deformations $p_\varepsilon$ in the space of polynomials of degree $n$ so that $p'_\varepsilon(z)$ has only simple roots near zero. Thus, after perturbing the projective structure a little bit and keeping the homomorphism $\theta$ the same, we assume that the order of each critical point of the meromorphic function $f : \Omega \to \mathbb{S}^2$ is 2. It is clear that this perturbation does not change $\langle \sigma(R), \sigma(R) \rangle$ and $d = \deg(D)$. The developing section $\sigma$ is transversal to the flat connection over all points of $R$ except at the singular points $\xi_1, ..., \xi_d$ of the structure. Let $D$ be the divisor of this singular locus. There exists a smooth vector field $X$ on $R$, which has $n = 2g + 2$ nondegenerate zeros, where $g$ is the genus of $R$: it has 1 sink, 1 source, and $2g$ saddle-type points. (For instance, take a Morse function $\mu : R \to \mathbb{R}$ which has one minimum, one maximum and $2g$ saddle points, then using a Riemannian metric on $R$ let $X := \mathrm{grad}(\mu)$.) Denote zeroes of $X$ by $\zeta_1, ..., \zeta_n$ where the last two points have index 1. We can choose $X$ so that

$$\{\zeta_1, \ldots, \zeta_n\} \cap \{\xi_1, \ldots, \xi_d\} = \emptyset.$$

Thus the vector field $X$ is a section of the tangent bundle $T_R$ which is transversal to the zero section. Now using the developing section $\sigma : R \to P$ we lift the vector field $X$ to a tangent vector field $Y = \sigma_*(X)$ along the surface $\Sigma = \sigma(R) \subset P$. The vertical directions in $P$ define the normal bundle $N(\Sigma)$ as in subsection 11.1. The flat connection on $P$ defines the projection $\nabla : T_x(P) \to V_x(P)$ where $V_x(P)$ is the distribution of vertical planes in $P$. The vector field $Q = \nabla(Y)$ is a section of the normal bundle $N(\Sigma)$. The section $\sigma$ is transversal to the flat connection on $P$ everywhere except at the set $\{\xi_1, \ldots, \xi_d\}$. Thus the set of zeros of the field $Q$ is

$$\sigma\{\xi_1, \ldots, \xi_d, \zeta_1, \ldots, \zeta_n\}.$$

A direct computation shows that the section $Q$ of the normal bundle $N(\Sigma)$ is transversal to the zero section $0_\Sigma$. Moreover, the intersection $Q(\Sigma) \cap 0_\Sigma$ is positive at the points $\{\xi_1, \ldots, \xi_d, \zeta_{n-1}, \zeta_n\}$ and is negative at the points $\{\zeta_1, \ldots, \zeta_{n-2}\}$. Hence the algebraic intersection number $\langle Q(\Sigma), 0_\Sigma \rangle$ (which is equal to $\langle \Sigma, \Sigma \rangle$) equals

$$d + 2 - (n - 2) = d + 2 - 2g$$

which proves the proposition. $\qquad\qquad\qquad\qquad\qquad\qquad\qquad\qquad\qquad\qquad\square$

COROLLARY 11.2.3.   *The degree* $\deg(D) = d$ *is even if and only if the representation* $\theta$ *lifts to* $\mathrm{SL}(2, \mathbb{C})$. *Equivalently,* $\theta$ *is liftable if and only if the second Stiefel-Whitney class satisfies the equation* $w_2(P) = \deg(D) = 0(\mathrm{mod}\ 2)$.

*Proof.* The representation $\theta$ lifts to $\mathrm{SL}(2, \mathbb{C})$ if and only if the bundle $P$ is trivial (equivalently, $w_2(P) = 0$); see [Go2]. As in the previous proposition we have the developing section $\sigma$ of the bundle $P \to R$. We proved that $\langle \sigma(R), \sigma(R) \rangle = 2 - 2g + \deg(D)$; hence

$$\langle \sigma(R), \sigma(R) \rangle = \deg(D) \ (\mathrm{mod}\ 2).$$

On the other hand, according to Lemma 11.1.1 we have:

$$\langle \sigma(R), \sigma(R) \rangle = w_2(P) \ (\mathrm{mod}\ 2)$$

and the corollary follows. $\qquad\qquad\qquad\qquad\qquad\qquad\qquad\qquad\qquad\qquad\qquad\square$

Now we are ready with the promised refinement of Theorem 1.1.1.

THEOREM 11.2.4.   *Suppose the surface* $R$ *and homomorphism* $\theta$ *satisfy the hypothesis of Theorem 1.1.1. Suppose that* $D$ *is a nonnegative divisor on* $R$ *such that* $w_2(\theta) = d \ (\mathrm{mod}\ 2)$, *where* $d = \deg(D)$. *Then there exists a complex projective structure on* $R$ *that has the monodromy* $\theta$ *and branching divisor equivalent to* $D$.

*Proof.* The proof is a straightforward generalization of the proof of Theorem 1.1.1. Let $P$ denote the $\mathbb{S}^2$-bundle over the surface $R$ associated with the homomorphism $\theta$. We first construct a decomposition of the surface $R$ into a union of pairs of pants so that the restriction of $\theta$ to the fundamental group of each pair of pants is a Schottky representation. We use these representations to build a complex projective structure on a pants configuration. But there is a "topological" $\mathbb{Z}/2$-obstruction to forming the final handle. This obstruction is a Dehn twist along a compressing loop. Suppose first that $w_2(P) = 0$. If the obstruction is nontrivial, then we can still construct a projective structure for the pants configuration which has exactly one branch point of order 1 and the monodromy $\theta$. However the existence of such a structure contradicts Corollary 11.2.2. Thus the "topological" obstruction to the existence of an unbranched structure was trivial to begin with. In parallel, we conclude that if $w_2(P) \neq 0$, then the pants configuration admits a branched structure with a single branch point of order 1. Now consider the general case assuming that $w_2(P) = 0$. By adding to the pants configuration (for example to a single pants in the configuration) branch points equivalent to the divisor $D$, we do not change the "topological" $\mathbb{Z}/2$-obstruction to completing the construction. Since $\deg(D) = 0 \pmod 2$, adding the branch points has the effect of twisting one of the boundary curves an even number of times. Hence for the resulting branched pants configuration there is no obstruction to completing it to a closed surface. The construction in the case $w_2(P) \neq 0$ is similar. $\qquad\square$

11.3. *The algebro-geometric interpretation.* Let $R$ be a closed Riemann surface $R$ of genus $g \geq 2$. In this section we shall consider holomorphic vector bundles $W$ over $R$ such that $\operatorname{rank}(W) = 2$ and $\det(W) = 1$ (i.e. the determinant bundle is trivial). Let $V^*(R)$ denote the collection of holomorphic vector bundles $W$ over $R$ such that $W$ admits a holomorphic flat connection. According to Weil's theorem (see [At], [Gu2], [W]), elements of $V^*(R)$ can be characterized intrinsically as follows:

Suppose that $W = \oplus_j W_j$ is the holomorphic direct sum decomposition of $W$ into (holomorphically) indecomposable vector bundles. Then the bundle $W$ admits a holomorphic flat connection if and only if $\deg(\det(W_j)) = 0$ for all $j$.

Let

$$F^*(R) := \{(\xi, \nabla) : \xi \in V^*(R), \nabla \text{ is a holomorphic flat connection on } \xi\}$$

be the space of *local systems* on $R$. We have the *Riemann-Hilbert correspondence*:

$$RH_R^* : F^*(R) \to Y(\pi_1(R), \operatorname{SL}(2, \mathbb{C})) := \operatorname{Hom}(\pi_1(R), \operatorname{SL}(2, \mathbb{C}))/\operatorname{SL}(2, \mathbb{C})$$

given by the conjugacy class of the monodromy of the flat connection $\nabla$. It is clear that the mapping $RH_R^*$ is bijective (since every flat bundle over $R$ has a canonical complex structure). The space $Y(\pi_1(R), \operatorname{SL}(2, \mathbb{C}))$ has a natural (non-Hausdorff) topology, we topologize $F^*(R)$ so that $RH_R^*$ is a homeomorphism.

We also have the natural projection

$$\pi_R^* : F^*(R) \longrightarrow V^*(R), \quad \pi_R^*(\xi, \nabla) := \xi.$$

Recall that each holomorphic vector bundle $W$ has the *degree of instability* $u(W)$ defined as follows:

> $u(W) = d$ is the maximal number such that $W$ contains a holomorphic line subbundle $L \subset W$ such that $\deg(L) = d$.

In general, $u(W) = d - \deg(\det(W))$.

For all bundles $W \in V^*(R)$,

$$-g \leq u(W) \leq g - 1$$

(see for instance [Gu2]), and *stable* (resp. *semistable*) bundles $W$ are defined by the condition $u(W) < 0$ (resp. $u(W) \leq 0$). Stable and semistable bundles and their moduli spaces have been extensively studied by algebraic geometers since the seminal paper of Narasimhan and Seshadri [N-S]. In contrast, our main objects are *maximally unstable* bundles $W$ which are defined by the condition $u(W) = g - 1$. Gunning [Gu1] proves that projectivizations of all maximally unstable bundles over $R$ are holomorphically isomorphic to each other. We let $M_R$ denote the corresponding projective bundle over $R$. It gives rise to a finite subset $M_R^*$ of $V_R^*$ that consists of $2^{2g}$ vector bundles that can be described as follows. Let $K$ denote the canonical bundle on $R$. Choose a holomorphic line bundle $L$ on $R$ such that $L^2 = K$. Then $\deg(L) = g - 1$. There are $2^{2g}$ characters $\chi : \pi_1(R) \to \{\pm 1\} \subset \mathbb{C}$. Each character gives rise to a holomorphic line bundle over $R$ which we shall denote by the same letter $\chi$. Then the collection of *square roots* $\sqrt{K}$ of the bundle $K$ consists of $2^{2g}$ bundles $\chi \otimes L$. For each $\Lambda = \chi \otimes L \in \sqrt{K}$ there is a unique holomorphically indecomposable bundle $W = W_\chi$ for which there is a short exact sequence

$$1 \to \Lambda \to W \to \Lambda^{-1} \to 1$$

of holomorphic morphisms of holomorphic bundles. Notice that $W_\chi = \chi \otimes W_1$ where $1 : \pi_1(R) \to \{\pm 1\}$ is the trivial homomorphism. Then

$$M_R^* = \{W_\chi, \chi : \pi_1(R) \to \{\pm 1\}\}.$$

Also in [Gu1], Gunning establishes the basic relation between maximally unstable bundles and complex projective structures on the surface $R$. He proves that

$$RH_R^*((\pi_R^*)^{-1}(M_R^*)) \subset Y(\pi_1(R), \mathrm{SL}(2, \mathbb{C}))$$

consists of (conjugacy classes of) linear monodromy representations of complex projective structures on the Riemann surface $R$. The relation between (branched) complex projective structures and instability of holomorphic vector bundles is further explored in [Man1, 2, 3].

The results of the previous two sections imply the following:

COROLLARY 11.3.1.  *Suppose that $\theta : \pi_1(R) \to \mathrm{PSL}(2, \mathbb{C})$ is the monodromy representation of a branched projective structure with branching divisor $D$. Let $P \to R$ denote the associated $\mathbb{S}^2$-bundle over $R$ which is the projectivization of a holomorphic vector bundle $V \to R$. Then $u(V) \geq g - 1 + [\deg(\det(V)) - \deg(D)]/2$.*

*Proof.* The developing map of the projective structure defines a section $\sigma$ of the bundle $P$, let $L \subset V$ be the corresponding line subbundle. Then Lemma 11.1.1 and Proposition 11.2.2 imply that

$$\deg(L) = g - 1 + [\deg(\det(V)) - \deg(D)]/2. \qquad \square$$

From now on it will be convenient to projectivize all vector bundles, connections and representations. Let

$$Y(\pi_1(R), \mathrm{PSL}(2, \mathbb{C})) := \{p(\rho), \rho \in \mathrm{Hom}(\pi_1(R), \mathrm{SL}(2, \mathbb{C}))\}/\mathrm{PSL}(2, \mathbb{C}) \subset V_g,$$

where $p(\rho)$ is the projectivization of $\rho$. Denote the spaces of projectivized holomorphic bundles and local systems over $R$ by $V(R)$ and $F(R)$ respectively. Let $RH_R : F(R) \to Y(\pi_1(R), \mathrm{PSL}(2, \mathbb{C}))$ denote the induced Riemann-Hilbert correspondence. Similarly define the projection $\pi_R$ by projectivizing the mapping $\pi_R^*$.

Our next step is to allow the complex structure on the surface $R$ to vary. We let $S$ be the oriented smooth surface underlying $R$. Let $\mathfrak{T}(S)$ denote the Teichmüller space of $S$. Consider the spaces

$$V_{\mathrm{top}}(S) := \bigcup_{R \in \mathfrak{T}(S)} V(R),$$

$$F_{\mathrm{top}}(S) := \bigcup_{R \in \mathfrak{T}(S)} F(R),$$

and mappings,

$$\pi : F_{\mathrm{top}}(S) \to V_{\mathrm{top}}(S) \ , \ RH : F_{\mathrm{top}}(S) \to Y(\pi_1(S), \mathrm{PSL}(2, \mathbb{C})),$$

whose restrictions to the fibers $F(R)$ are $\pi_R : F(R) \to V(R)$ and $RH_R$.

*Remark* 11.3.2. The space $F_{\mathrm{top}}(S)$ is naturally identified with the product

$$F_{\mathrm{top}}(S) = \mathfrak{T}(S) \times Y(\pi_1(S), \mathrm{PSL}(2, \mathbb{C})).$$

The projection $F_{\mathrm{top}}(S) \to \mathfrak{T}(S)$ which maps $F(R, \phi)$ to $(R, \phi) \in \mathfrak{T}(S)$ is the projection of $F_{\mathrm{top}}(S)$ to the first factor of the product decomposition.

Indeed, suppose $(R, \phi) \in \mathfrak{T}(S)$ is a *marked Riemann surface* with the marking $\phi : \pi_1(S) \to \pi_1(R)$ (which is an isomorphism defined up to an inner automorphism). Then $\phi$ indices an natural isomorphism

$$Y(\pi_1(R), \mathrm{PSL}(2, \mathbb{C})) \to Y(\pi_1(S), \mathrm{PSL}(2, \mathbb{C}))$$

given by precomposition of representations with $\phi$. Note that we have to work with the Teichmüller space of $S$ rather than with the moduli space $\mathfrak{M}(S)$, otherwise the natural projection to $\mathfrak{M}(S)$ would be a nontrivial fibration (in the orbifold sense).

The projection

$$\Pi : V_{\text{top}}(S) \to \mathfrak{T}(S), \Pi : V(R) \to R$$

has a section

$$\mu : R \mapsto M_R \in V(R) \subset V_{\text{top}}(S),$$

where $M_R$ is the projectivization of maximally unstable vector bundles over $R$. Let

$$Y_{ne}(\pi_1(S), \text{PSL}(2, \mathbb{C})) \subset V'_g$$

denote the collection of conjugacy classes of all projectivized nonelementary representations into $\text{SL}(2, \mathbb{C})$. We summarize this in the diagram below:

$$
\begin{array}{ccc}
\mathfrak{T}(S) & \xrightarrow{\ \mu\ } & V_{\text{top}}(S) \\
\big\uparrow \Pi \circ \pi & & \big\uparrow \pi \\
\pi^{-1}(\mu(\mathfrak{T}(S))) & \subset & F_{\text{top}}(S) \\
\big\downarrow RH & & \big\downarrow RH \\
Y_{ne}(\pi_1(S), \text{PSL}(2, \mathbb{C})) & \subset & Y(\pi_1(S), \text{PSL}(2, \mathbb{C}))
\end{array} \cdot
$$

In view of [Gu1], the image $RH(\pi^{-1}(\mu(\mathfrak{T}(S))))$ consists of (projective) monodromy representations of complex projective structures on the surface $S$. On the other hand, each holomorphic bundle in $M_R^*$ is maximally unstable. Let $V_\rho$ be a maximally unstable bundle associated with a representation $\rho : \pi_1(R) \to \text{SL}(2, \mathbb{C})$. Thus, for all characters $\chi : \pi_1(R) \to \{\pm 1\}$, the bundles $\chi \otimes V_\rho = V_{\chi \cdot \rho}$ are also maximally unstable. The inverse image of the subvariety $\pi^{-1}(M_R)$ in $Y(\pi_1(R), \text{SL}(2, \mathbb{C}))$ has $2^{2g}$ components. Each component consists of holomorphically isomorphic vector bundles over $R$, but members of distinct components are not holomorphically isomorphic to each other.

Therefore, by applying Theorem 1.1.1, we obtain,

THEOREM 11.3.3. *The map $RH$ sends $\pi^{-1}(\mu(\mathfrak{T}(S)))$ onto $Y_{ne}(\pi_1(S), \text{PSL}(2, \mathbb{C}))$. In other words, let $\rho \in Y(\pi_1(S), \text{SL}(2, \mathbb{C}))$ be a nonelementary representation. It is the monodromy of a holomorphic flat connection on a maximally unstable holomorphic vector bundle over a Riemann surface $R$; $R$ is diffeomorphic to $S$ via an orientation-preserving diffeomorphism.*

11.4.  *Proper embeddings in the representation variety.* In this section we will give a detailed proof of the "divergence" theorem. It was first suggested by Hejhal in [He1] that such theorem could be true. This theorem shows that on a fixed Riemann surface, if any sequence of quadratic differentials diverge, so must the conjugacy classes of corresponding monodromy representations. A brief outline of the proof was given in [Ka, §7.2].[2]

As before, $R$ denotes a closed Riemann surface of genus exceeding one and $Q(R)$ its space of holomorphic quadratic differentials. Let hol denote the map that sends each $\phi \in Q(R)$ to the monodromy homomorphism determined by the corresponding Schwarzian equation $S(f) = \phi$. By Theorem 1.1.1, the image lies in the component of the representation variety $V_g$ containing the identity (cf., §1.5). That is,

$$\mathrm{hol} : Q(R) \to Y(\pi_1(R), \mathrm{PSL}(2, \mathbb{C})).$$

THEOREM 11.4.1 (Divergence Theorem).  *The map* hol *is proper.*

*Proof.* Let $\tilde{Z} \subset \mathrm{Hom}(\pi_1(R), \mathrm{SL}(2, \mathbb{C}))$ denote the preimage of $Z$, where $Z$ is the image of hol.

Our first goal is to show that $\tilde{Z}$ is a properly embedded complex analytic subvariety in $\mathrm{Hom}(\pi_1(R), \mathrm{SL}(2, \mathbb{C}))$. Indeed, if $\rho : \pi_1(R) \to \mathrm{SL}(2, \mathbb{C})$ is any representation, the associated vector bundle $V_\rho \to R$ is maximally unstable if and only if $\rho \in \tilde{Z}$. Equivalently,

$$\rho \in \tilde{Z} \iff H^0(R, L^* \otimes V_\rho) \neq 0 \ \text{ for some } \ L \in \sqrt{K}.$$

The set $\sqrt{K}$ is finite. Thus, by the upper semicontinuity theorem for cohomology (see [B-S]), the subset $\tilde{Z}$ is closed and is equal to a finite union of disjoint complex analytic subvarieties $X_L$ properly embedded in $\mathrm{Hom}(\pi_1(R), \mathrm{SL}(2, \mathbb{C}))$ (these subvarieties are indexed by $L \in \sqrt{K}$).

Recall that $\tilde{Z}$ is contained in the open subset $\mathrm{Hom}_{\mathrm{ne}}(\pi_1(R), \mathrm{SL}(2, \mathbb{C}))$ of nonelementary representations, *i.e.* those whose projectivizations are nonelementary. The group $\mathrm{SL}(2, \mathbb{C})$ acts on $\mathrm{Hom}_{\mathrm{ne}}(\pi_1(R), \mathrm{SL}(2, \mathbb{C}))$ by conjugation and the quotient is $Y_{\mathrm{ne}}(\pi_1(R), \mathrm{SL}(2, \mathbb{C}))$. Hence the projection

$$\mathrm{Hom}_{\mathrm{ne}}(\pi_1(R), \mathrm{SL}(2, \mathbb{C})) \to Y_{\mathrm{ne}}(\pi_1(R), \mathrm{SL}(2, \mathbb{C}))$$

is a principal $\mathrm{SL}(2, \mathbb{C})$-bundle. Since $\tilde{Z}$ is invariant under this action, the projection $Z^*$ of $\tilde{Z}$ to $Y_{\mathrm{ne}}(\pi_1(R), \mathrm{SL}(2, \mathbb{C}))$ is again a closed properly embedded complex analytic subvariety. It consists of $2^{2g}$ components indexed by elements of $\sqrt{K}$.

The restriction of the projection

$$p : Y_{\mathrm{ne}}(\pi_1(R), \mathrm{SL}(2, \mathbb{C})) \to Y_{\mathrm{ne}}(\pi_1(R), \mathrm{PSL}(2, \mathbb{C}))$$

---

[2]Note that the discussion in [Ka, §7.2] does not distinguish linear and projective monodromy representations.

to each component of $Z^*$ is a bijection onto hol$(Q(R))$. Now $p(Z^*) = Z$ is closed, since $p$ is a finite covering. It is disjoint from the collection of conjugacy classes of elementary representations because all elementary representations correspond to semistable bundles over $R$. Consequently we can restrict our study to the smooth (Hausdorff) manifold $Y_{ne}(\pi_1(R), \mathrm{PSL}(2, \mathbb{C}))$.

According to [Gu2], the partition of $Y_{ne}(\pi_1(R), \mathrm{SL}(2, \mathbb{C}))$ into holomorphic equivalence classes is a smooth foliation. The components of $Z^*$ are leaves of this foliation; hence they are complex submanifolds in $Y_{ne}(\pi_1(R), \mathrm{SL}(2, \mathbb{C}))$. This implies that $Z \subset Y_{ne}(\pi_1(R), \mathrm{PSL}(2, \mathbb{C}))$ is a properly embedded complex submanifold. On the other hand, the mapping hol : $Q(R) \to Z$ is a continuous bijection, hence a homeomorphism. Therefore hol : $Q(R) \to Z \subset Y_{ne}(\pi_1(R), \mathrm{PSL}(2, \mathbb{C}))$ is proper. Hence the composition of

$$\mathrm{hol} : Q(R) \to Y_{ne}(\pi_1(R), \mathrm{PSL}(2, \mathbb{C}))$$

with the inclusion

$$Y_{ne}(\pi_1(R), \mathrm{PSL}(2, \mathbb{C})) \hookrightarrow Y(\pi_1(R), \mathrm{PSL}(2, \mathbb{C}))$$

is a proper map $Q(R) \to Y(\pi_1(R), \mathrm{PSL}(2, \mathbb{C}))$.                    □

*Remark* 11.4.2.   The above proof shows that elementary representations cannot be limits of sequences from hol$(Q(R))$. It was proven [Ka] only that the mapping hol : $Q(R) \to Y_{ne}(\pi_1(R), \mathrm{PSL}(2, \mathbb{C}))$ is proper. Tanigawa [Tani] recently gave a nice geometric proof of this statement in contrast to algebro-geometric proof presented here and in [Ka]. However Tanigawa's arguments do not seem to prove that $Z = \mathrm{hol}(Q(R))$ is closed in $Y(\pi_1(R), \mathrm{PSL}(2, \mathbb{C}))$, only in the submanifold corresponding to nonelementary representations. See also §12.4.

11.5. *An analogue of Hejhal's holonomy theorem for branched projective structures.* The nonelementary representation variety $V'_g$ has two components [Go2]. These correspond to the representations that lift to $\mathrm{SL}(2, \mathbb{C})$, and those that do not. Each of these has dimension $6g-6$. By a *singly branched projective structure* we mean one that has exactly one branch point and that is of order two. In the next section we will show that the monodromy of each singly branched projective structure is a nonelementary representation but we will use this fact in this section.

Let $R$ be a closed Riemann surface of genus $g \geq 2$ and $p \in R$ a given point. We will first parameterize singly branched structures on $R$ with branch point at the designated point $p$. Let $D$ be the divisor of $p$ and $Q_D(R)$ the space of holomorphic quadratic differentials on $R$ which have at most a simple pole at $p$.

Recall from §1.4, equation (6), that the meromorphic quadratic differential $\phi_0$ generates a singly branched complex projective structure if its Laurent expansion at the chosen branch point $p$ has the form

$$\phi_0 = -3/z^2 + b/z + \alpha_0 + \alpha_1 z + \cdots, \qquad \text{where} \ \ b^2 + 2\alpha_0 = 0,$$

(here and below we choose local coordinates so that $p$ is identified with zero). The side condition comes from the requirement that the solution of the Schwarzian equation has no logarithmic term.

We note that there exists such a differential $\phi_0$. First of all the Riemann-Roch theorem implies there is a quadratic differential with a double pole at any point $p$. Secondly it also implies that there is an abelian differential $\omega$ which does not vanish at $p$. The holomorphic differential $\omega^2$ can be employed to insure that the side condition is satisfied ([Man1]). Fix one such quadratic differential $\phi_0$.

There is a meromorphic quadratic differential with a single pole at $p$ with the Laurent expansion

$$\psi_0 = 1/z + a_0 + d_1 z + \cdots.$$

Adding $\omega^2$, which does not vanish at $p$, to $\psi_0$ if necessary, we may assume that $a_0 + b \neq 0$.

Let $\psi_i$, $1 \leq i \leq 3g - 3$ be a basis of the holomorphic quadratic differentials on $R$. Then $\psi_i$, $0 \leq i \leq 3g - 3$ is a basis of the space $Q_D(R)$.

Let $a_i$ be the constant term in the Laurent expansion of $\psi_i$ at $p$. Not all $a_i$ can vanish.

The vector space $Q_D(R)$ consists of the differentials $\psi = \sum_{i=0}^{3g-3} c_i \psi_i$. When is $\phi_0 + \psi$ an admissible quadratic differential in the sense of §1.4, equation (6)? The answer is when $u^2 + 2bu + 2v = 0$, where $v$ is the constant term in in the Laurent expansion of $\psi$ at $p$, and $u$ is its residue.

The constant term in $\psi$ is $v = \sum_{i=0}^{3g-3} c_i a_i$. The residue term is just $c_0$. Hence the condition reads

$$(11) \qquad c_0^2 + 2bc_0 + 2 \sum_{i=0}^{3g-3} c_i a_i = 0.$$

Recall that $a_0 + b \neq 0$, thus the implicit function theorem implies that the collection of vectors $\vec{c} = (c_0, c_1, \ldots, c_{3g-3})$ satisfying the above equation is a complex manifold of dimension $3g - 3$ provided that the norm $|\vec{c}|$ is sufficiently small. (Actually it suffices to require that only $|c_0|$ is sufficiently small.)

Consequently we can choose a small neighborhood $U$ of $\phi_0$ in the affine space of meromorphic quadratic differentials $\phi_0 + Q_D(R)$ with the following property.

The collection of differentials $\phi_0 + \sum_{i=0}^{3g-3} c_i \psi_i \in U$ satisfying (11) forms a $(3g - 3)$-dimensional complex manifold $\Delta$ containing $\phi_0$.

Let $B_g$ denote the holomorphic variety which consists of singly branched complex projective structures on closed Riemann surfaces $S$ of genus $g \geq 2$.

Let $\mathfrak{T}(S - \{q\})$ denote the Teichmüller space of surfaces $S$ with one marked point. There is a holomorphic mapping $\nu : B_g \to \mathfrak{T}(S - \{q\})$ whose fiber over a marked Riemann surface $R$ with a marked point $p$ is the space $I(R, D)$ of singly branched complex projective structures with the underlying complex structure $R$ and branching at $D = p$. It follows from the above discussion that $B_g$ is a holomorphic variety of generic dimension $6g - 5$: the Teichmüller space of once punctured surfaces $\mathfrak{T}(S - \{q\})$ has complex dimension $3g - 2$ and the fiber of $\nu$ has complex dimension $3g - 3$. There is an open and dense subset of $B_g$ which is a complex manifold of dimension $6g - 5$. We will use the notation $(S, q, \varphi)$ for elements of $B_g$, where $S$ denotes the marked Riemann surface, $q$ the branch point and $\varphi$ the meromorphic quadratic differential.

We will need the following explicit description of the space $B_g$. Choose a point $R$ as the "origin" in $\mathfrak{T}(S)$ and write $R = \mathbb{H}^2/G$ where $\mathbb{H}^2$ is the unit disk $\{|z| < 1\}$ and $G$ is a fuchsian group acting on $\mathbb{H}^2$. In the "Bers' slice" model, Teichmüller space $\mathfrak{T}(S)$ is identified with that subset of the space $Q(R)$ of holomorphic quadratic differentials on $R$, lifted to $\mathbb{H}^2$, such that the corresponding developing map $h_\tau : \mathbb{H}^2 \to \mathbb{S}^2$, $\tau \in Q(R)$, is a univalent holomorphic mapping with homeomorphic extension to $\{|z| = 1\}$. Thus $h = h_\tau$ solves the Schwarzian equation for $\tau$; we will normalize it by the requirement that $h(0) = 0, h'(0) = 1$. Let $\rho_\tau : G \to G_\tau$ denote the corresponding monodromy representation. As $\tau \to 0$, $G_\tau$ converges algebraically back to $G$.

The image $\rho_\tau(G) = G_\tau$ is a quasifuchsian group. Its set of discontinuity has two components. One is $\underline{\Omega}_\tau = h_\tau(\mathbb{H}^2)$. The other $\Omega_\tau$ represents the marked Riemann surface $R_\tau := \Omega_\tau/G_\tau \in \mathfrak{T}(S)$. The homotopy marking of this point in $\mathfrak{T}(S)$ is given by the isomorphism $\rho_\tau : G = \pi_1(R) \to G_\tau = \pi_1(R_\tau)$. If we mark a point $p \in R_\tau$ we get an element of $\mathfrak{T}(S - \{q\})$.

Any given compact subset of $\Omega_0$ belongs to $\Omega_\tau$ for $\tau$ sufficiently close to 0; likewise any neighborhood of the closure of $\Omega_0$ contains $\Omega_\tau$ for $\tau$ sufficiently close to 0. Here $\Omega_0 = \{z : |z| > 1\} \cup \infty$.

LEMMA 11.5.1. *There is a locally defined holomorphic map $P : B_g \to \mathbb{S}^2$ that "records" the position of the branch points.*

*Proof.* We construct $P$ in a small neighborhood $\Delta$ of a given point $(R, q, \varphi) \in B_g$, where $\varphi$ is a quadratic differential on $R$, the surface $R$ is the "origin" in the Teichmüller space and $q \in R$ is the branch point. We will denote points $\sigma \in \Delta$ by $(R_\tau, q_\sigma, \varphi_\sigma)$ where $\tau = \tau(\sigma) \in \mathfrak{T}(S)$ and $q_\sigma \in R_\tau$ is the branch-point. If $\tau = 0$, then $\sigma$ represents a change of branch point from $q$ to $q_\sigma$ on $R$ itself. The point $\sigma = 0$ is $(R, q, \varphi)$. Let $f_0 : \Omega_0 \to \mathbb{S}^2$ denote the (as yet unnormalized) developing map of $(R, q, \varphi)$ and $\theta : G = \pi_1(R) \to \mathrm{PSL}(2, \mathbb{C})$ the associated nonelementary monodromy representation (here we are applying Theorem 11.6.1 that will be proven in the next section).

Let $f_\sigma : \Omega_\tau \to \mathbb{S}^2$ be the associated holomorphic developing map. We will show in the next paragraph how to fix a consistent normalization for $f_\sigma$ given $f_0$ so that the restrictions of $f_\sigma$ to compact domains in $\Omega_0$ depend holomorphically on $\sigma$. Each developing mapping $f_\sigma$ corresponds to the monodromy representation

$$\theta'_\sigma := \theta_\sigma \circ \rho_\tau : G \to G_\tau \to H_\sigma \subset \mathrm{PSL}(2, \mathbb{C}).$$

At the origin, $\theta'_0 = \theta_0$.

Consider the projection $\mathrm{Hom}(\pi_1(R), \mathrm{PSL}(2,\mathbb{C})) \to V_g$. We will construct a local cross section $\widetilde{V}_g$ near $\theta'_0$ as follows. We know from Part A that we can find in $H_0$ three loxodromic elements $h_1 = \theta'_0(g_1), h_2 = \theta'_0(g_2), h_3 = \theta'_0(g_3)$ with distinct attracting fixed points $a_1, a_2, a_3$, where $g_1, g_2, g_3 \in G$. Normalize each developing mapping $f_\sigma$ so that the attracting fixed point of $\theta'_\sigma(g_j)$ remains $a_j$, $j = 1, 2, 3$. This can be done for all $\sigma \in \Delta$ if $\Delta$ is sufficiently small, i.e., if the attracting fixed points remain distinct and the elements $\theta'_\sigma(g_j)$ remain loxodromic. Thus in $\Delta$ we have a holomorphic lift

$$\widetilde{\mathrm{hol}} : B_g \to \widetilde{V}_g \subset \mathrm{Hom}(G, \mathrm{PSL}(2, \mathbb{C})).$$

Now given a lift $q^* \in \Omega_0$ of $q \in R$, there is, in the set of lifts of $q_\sigma$ to $\Omega_\sigma$, a closest (in the spherical metric) point $q^*_\sigma$ to $q^*$. Define

$$P : (R_\tau, q_\sigma, \varphi_\sigma) \mapsto f_\sigma(q^*_\sigma) \in \mathbb{S}^2.$$

It is clear that the mapping $P$ is holomorphic provided that $\Delta$ is so small that the point $q^*_\sigma$ is unique.  □

Thus, by the previous lemma we have a locally defined holomorphic map

$$\mu \quad = \quad (P, hol) : B_g \to \mathbb{S}^2 \times V_g$$

and its lift

$$\tilde{\mu} \quad = \quad (P, \widetilde{\mathrm{hol}}) : B_g \to \mathbb{S}^2 \times \widetilde{V}_g .$$

We are now ready to state our theorem.

THEOREM 11.5.2.  *The holonomy map* $\mathrm{hol} : B_g \to V_g$ *is locally a topological fiber bundle with fiber of complex dimension one.*

*Remark* 11.5.3.  The fibers reflect the choice of branch point. This result should generalize to the space of $D$-branched projective structures where $D$ is a fixed (topological) branching divisor, provided we consider structures with nonelementary monodromy.

*Proof.* In Lemmas 11.5.5 and 11.5.4 below we will prove that $\mu$ is injective and an open map. Hence $\mu$ is a local homeomorphism. Since $\mathbb{S}^2 \times V'_g$ is a complex manifold of dimension $(6g-5)$ we can therefore use $\mu$ to locally identify $B_g$ with the product $\mathbb{S}^2 \times V_g$ so that hol is identified with the projection to the second factor. $\square$

LEMMA 11.5.4. *Let $X$ be a holomorphic variety of generic complex dimension $n$ (i.e. there is an open dense subset $U \subset X$ which is a complex manifold of dimension $n$). Let $f : X \to M$ be a locally injective holomorphic mapping, where $M$ is a complex manifold of dimension $n$. Then $f$ is open.*

*Proof.* Since this is a local question it suffices to consider the germ of $X$ at a point $x \in X$ and the germ of $f$ at $x$. Since $f$ is locally injective, the germ of the mapping $f$ at $x$ is "finite" in the terminology of [Gu4, p. 56].

Suppose that the germ of $f$ at $x$ is not onto. Apply [Gu4, Corollary 9]: it follows that there exists a nonzero germ of a holomorphic function $h$ on $M$ at $m = f(x)$ such that $h \circ f = 0$. The germ at $m$ of the zero level set $\{h = 0\}$ of $h$ is a holomorphic subvariety of dimension strictly less than $n$, by the uniqueness principle of holomorphic functions. Thus the germ of the image $f(X)$ at $m$ has generic dimension less than $n$. However $f(X)$ is generically a manifold, hence $f(U)$ has dimension less than $n$, a contradiction to invariance of domain for manifolds. $\square$

LEMMA 11.5.5. *The mapping $\mu$ is locally injective.*

*Proof.* It suffices to show that two nearby branched structures with the same monodromy representation are identical provided that the images of their branch points under $P$ are the same. Our proof is analogous to that of [He1, Theorem 1]. It clearly enough to show local injectivity of the holomorphic lift

$$\tilde{\mu} = (P, \widetilde{\text{hol}}) : B_g \to \mathbb{S}^2 \times \tilde{V}_g .$$

We consider the points $\sigma = (R_\tau, q_\sigma, \varphi_\sigma)$ of a small neighborhood $\Delta$ of the point $(R, q, \varphi) \in B_g$.

Let $\mathcal{F}_\sigma \subset \Omega_\tau$ denote the (closed) Dirichlet fundamental domain for $G_\tau$ in the hyperbolic metric on $\Omega_\tau$ and with center $q_\sigma^*$; $\tau = \tau(\sigma)$. Let $\mathcal{F}_0^*$ be a small open neighborhood of $\mathcal{F}_0$, and take $\Delta$ so small that $\mathcal{F}_\sigma \subset \mathcal{F}_0^*$ for all $\sigma \in \Delta$. We may also assume that the orbit $G_\tau(q_\sigma^*)$ meets the closure of $\mathcal{F}_0^*$ only at $q_\sigma^*$.

We again use the developing mappings $f_\sigma : \Omega_\tau \to \mathbb{S}^2$. Decreasing $\Delta$ even more if necessary, we may assume that:

a) For each $\sigma \in \Delta$ there is an open neighborhood $\mathcal{F}_\sigma^*$ of $\mathcal{F}_\sigma$ such that for any pair of points $\sigma, \delta \in \Delta$ we have $\mathcal{F}_\delta \subset \mathcal{F}_\sigma^*$, and

b) Given small $\varepsilon > 0$, there is a disk $V \subset \mathcal{F}_0$ about $q^*$ of radius $2\varepsilon$ with the following property. To any $z \in \mathcal{F}_0^* \setminus V$, and to any pair of points $\sigma, \delta \in \Delta$, corresponds a unique point $z_{\sigma,\delta} \in \mathcal{F}_0^*$ such that

$$f_\delta(z_{\sigma,\delta}) = f_\sigma(z), \quad \text{and} \quad d(z_{\sigma,\delta}, z) < \varepsilon.$$

Here $d(\cdot, \cdot)$ is the spherical metric.

The "membranes" $\{f_\sigma(\mathcal{F}_\sigma)\}$ over $\mathbb{S}^2$ serve as "fundamental domains" for the image groups $H_\sigma = \theta'_\sigma(G)$.

Suppose $\tilde{\mu}$ is not injective in any neighborhood $\Delta$. Then for arbitrarily small $\Delta$ there exist $\sigma \neq \delta \in \Delta$ so that $P(\sigma) = P(\delta)$ and the normalized monodromy representations are identical, i.e.,

$$(\theta'_\sigma : G \to H_\sigma) \equiv (\theta'_\delta : G \to H_\delta).$$

We claim that there is a branch $F$ of $f_\delta^{-1} \circ f_\sigma$ which is a conformal homeomorphism of the fundamental domain $\mathcal{F}_\sigma$ onto a new fundamental domain $\mathcal{F}'_\delta$ for $G_{\tau(\delta)}$ in $\Omega_{\tau(\delta)}$. Such a map $F$ would necessarily be equivariant in the sense that if $z, g(z) \in \mathcal{F}_\sigma$, $g \in G_{\tau(\sigma)}$, then $F(z), F(g(z)) \in \mathcal{F}'_\delta$ and $F(g(z)) = \rho_{\tau(\delta)} \circ \rho_{\tau(\sigma)}^{-1}(g)F(z)$. Here $\xi := \rho_{\tau(\delta)} \circ \rho_{\tau(\sigma)}^{-1} : G_{\tau(\sigma)} \to G_{\tau(\delta)}$ is the isomorphism which factors through $G$.

Indeed, for $z \notin V$ define $F(z) := z_{\sigma,\delta}$. It is clear $F$ is a univalent holomorphic mapping. Furthermore $F|\mathcal{F}_\sigma \setminus V$ extends over $V$ to a conformal mapping because both $f_\sigma$ and $f_\delta$ are 2-fold branched coverings near $q^* \in \Omega_0$ with the same critical value

$$f_\sigma(q_\sigma^*) = P(\sigma) = P(\delta) = f_\delta(q_\delta^*).$$

The mapping $F$ projects to a conformal map of $R_{\tau(\sigma)} = \Omega_{\tau(\sigma)}/G_{\tau(\sigma)}$ onto $R_{\tau(\delta)} = \Omega_{\tau(\delta)}/G_{\tau(\delta)}$. Correspondingly $F$ extends to a conformal mapping $F : \Omega_{\tau(\sigma)} \to \Omega_{\tau(\delta)}$ that induces the isomorphism $\xi : G_{\tau(\sigma)} \to G_{\tau(\delta)}$.

The map $h_{\tau(\delta)} \circ h_{\tau(\sigma)}^{-1}$ is a conformal map of $\underline{\Omega}_{\tau(\sigma)}$ onto $\underline{\Omega}_{\tau(\delta)}$ which also induces the isomorphism $\xi$. The two conformal mappings have continuous extensions to the limit set which are necessarily identical. Since the limit set is a quasicircle they are the restrictions of a Möbius transformation. In particular $F$ is a Möbius transformation and $\sigma = \delta$, a contradiction. $\qquad\square$

The following is a direct consequence of Theorem 11.5.2.

COROLLARY 11.5.6. *Let* $\sigma = (R, p, \varphi)$ *be a singly branched projective structure. Let* $\Delta \subset B_g$ *be a sufficiently small neighborhood of* $\sigma$ *in the space of singly branched structures* $\delta$ *on* $R$ *"with the same image of the branch point"* $P(\delta)$ *as* $\sigma$. *Suppose the sequence of normalized representations* $\theta_i : \pi_1(R) \to \mathrm{PSL}(2, \mathbb{C})$ *converges algebraically to the normalized monodromy representation* $\theta$ *associated with* $\sigma$. *Then for all large* $i$, $\theta_i$ *is associated with a unique* $\sigma_i \in \Delta$.

11.6. *Monodromy of singly branched projective structures.* In this section we will prove facts that have been announced in §1.6.

Theorem 11.6.1.   *Suppose that $R$ is a closed Riemann surface of genus $g \geq 2$, $\theta : \pi_1(R) \to \mathrm{PSL}(2,\mathbb{C})$ is the monodromy representation of a singly branched complex projective structure $\tau$ on $R$. Then $\Gamma = \theta(\pi_1(R))$ is a nonelementary subgroup of $\mathrm{PSL}(2,\mathbb{C})$.*

*Proof.* Since $\tau$ has exactly one branch point and the order of this branch point is 2, the representation $\theta$ has nonzero 2nd Stiefel-Whitney class. In particular, $\theta$ cannot be lifted to a representation $\pi_1(R) \to \mathrm{SL}(2,\mathbb{C})$. Suppose that the group $\Gamma = \theta(\pi_1(R))$ is elementary. There are three cases:

(a) The group $\Gamma$ has a fixed point $z \in \mathbb{S}^2$. Without loss of generality we can assume that $z = \infty$, thus $\Gamma$ is contained in the group $A$ of complex affine transformations of $\mathbb{C}$. The inclusion $A \hookrightarrow \mathrm{PSL}(2,\mathbb{C})$ admits a 1-1 lift $A \hookrightarrow \mathrm{SL}(2,\mathbb{C})$

$$a^2 z + b \mapsto \begin{pmatrix} a & ba^{-1} \\ 0 & a^{-1} \end{pmatrix}.$$

Therefore $\theta$ lifts to a representation $\theta^* : \pi_1(R) \to \mathrm{SL}(2,\mathbb{C})$, which contradicts the assumption that $\theta$ has nonzero 2nd Stiefel-Whitney class.

(b) Suppose that $\Gamma$ is conjugate into the subgroup $PU(2) \subset \mathrm{PSL}(2,\mathbb{C})$. Let $\tilde{R} \to R$ be a 2-fold covering over $R$. Thus $2(g-1) = \tilde{g}-1$, where $\tilde{g}$ denotes the genus of $\tilde{R}$. The complex projective structure $\tau$ on $R$ defines a complex projective structure $\tilde{\tau}$ on $\tilde{R}$ with two branch points of order two. Suppose that $\Gamma \subset PU(2)$; then $\theta(\pi_1(\tilde{R})) \subset PU(2)$ as well. The representation $\theta|\pi_1(\tilde{R})$ lifts to a linear representation

$$\theta^* : \pi_1(\tilde{R}) \to SU(2) \subset \mathrm{SL}(2,\mathbb{C}).$$

Consider the flat vector bundle $V$ of the rank 2 over the surface $\tilde{R}$ associated with the action $\theta^*$ of $\pi_1(\tilde{R})$ on $\mathbb{C}^2$. Clearly $\det(V) = 1$. The developing map of the branched complex projective structure $\tilde{\tau}$ defines a section

$$\sigma : \tilde{R} \to P(V).$$

According to Proposition 11.2.2, the self-intersection number $\sigma^2$ of the surface $\sigma(\tilde{R})$ in $P(V)$ equals $(2-2\tilde{g})+2$, since the structure $\tilde{\tau}$ has exactly two branch points of the order 2.

It follows from Lemma 11.1.1 that the section $\sigma$ gives rise to a line sub-bundle $L \subset V$ such that

$$\deg(L) = (\tilde{g}-1) - 1 = 2g - 3 > 0.$$

We conclude that $u(V) > 0$ and the bundle $V$ is unstable. On the other hand, every flat bundle over $\tilde{R}$ with unitary monodromy group is semistable (see for instance [N-S]). This contradiction shows that $\Gamma$ cannot be contained in $PU(2)$.

(c) Consider the case that the group $\theta(\pi_1(R))$ has an invariant pair of points in $\mathbb{S}^2$. (This does not imply that $\theta$ can be lifted to $\mathrm{SL}(2,\mathbb{C})$.) We argue as in Case (b). There is a 2-fold covering $\tilde{R} \to R$ such that the group $\theta(\pi_1(\tilde{R}))$ has a pair of fixed points in $\mathbb{S}^2$. Therefore the induced complex projective structure on $\tilde{R}$ has two branch points and the monodromy group $\theta(\pi_1(\tilde{R}))$ has a lift $\theta^*(\pi_1(\tilde{R}))$ to a subgroup of $\mathrm{SL}(2,\mathbb{C})$ conjugate to the group of diagonal matrices. Let $V$ denote the holomorphic vector bundle associated with the representation $\theta^* : \pi_1(\tilde{R}) \to \mathrm{SL}(2,\mathbb{C})$. The representation $\theta^*$ splits as the direct sum of representations. Hence the bundle $V$ is decomposable (into the direct sum of two line bundles of degree zero), which implies that $u(V) = 0$. On the other hand, the developing map of the branched complex projective structure on $\tilde{R}$ defines a section $\sigma : R \to P(V)$ with the self-intersection number

$$(2 - 2\tilde{g}) + 2 < 0,$$

where $\tilde{g}$ denotes the genus of $\tilde{R}$. Hence $u(V) > 0$ which contradicts $u(V) = 0$. □

Suppose that $\tau$ is a branched complex projective structure on the closed Riemann surface $R$ of genus at least two. We identify the universal cover of $R$ with the hyperbolic plane $\mathbb{H}^2$. Let $f : \mathbb{H}^2 \to \mathbb{S}^2$ be the developing map of $\tau$ and $\Gamma = \theta(\pi_1(R))$ the holonomy group. We say that $\tau$ is a *branched hyperbolic structure* if $\tau$ has at least one branch point and the image of $f$ is a round disk in $\mathbb{S}^2$. This definition is motivated by the fact that in such case $\Gamma$ preserves the hyperbolic metric $ds^2$ in $f(\mathbb{H}^2)$. The pull back of $ds^2$ from $f(\mathbb{H}^2)$ to $R$ is a hyperbolic metric on $R$ which has singular points at the branch points $z_j$ of $\tau$; the total angle around $z_j$ is $2\pi k_j$, where $k_j$ is the order of $z_j$.

Later we will show by example why the following result is false if we do not exclude branched hyperbolic structures. This too has been announced in §1.6.

COROLLARY 11.6.2. *Suppose that either the complex projective structure $(f, \theta)$ is unbranched, or is singly branched but is not a branched hyperbolic structure (i.e. $f(\mathbb{H}^2)$ is not a round disk). Then the following statements are equivalent*:

(i)  $f(\mathbb{H}^2) \neq \mathbb{S}^2$;

(ii)  $\mathbb{H}^2 \to f(\mathbb{H}^2)$ *is a (possibly branched) cover*;

(iii)  $\Gamma$ *acts discontinuously on $f(\mathbb{H}^2)$.*

*Proof.* The unbranched case is classical (see §1.6). Consider then the branched case. By Theorem 11.6.1, $\Gamma = \theta(\pi_1(R))$ is nonelementary. The limit set $\Lambda(\Gamma)$ is the smallest $\Gamma$-invariant closed nonempty subset of $\mathbb{S}^2$. Since $\Gamma$ is nonelementary, $\Lambda(\Gamma)$ is the closure of the set of fixed points of loxodromic elements of $\Gamma$. It follows that the $\Gamma$-orbit of any open set containing a limit point is $\mathbb{S}^2$. Suppose that $\Gamma \subset \mathrm{PSL}(2, \mathbb{C})$ is nondiscrete. Let $\bar{\Gamma}$ be the closure of $\Gamma$ in $\mathrm{PSL}(2, \mathbb{C})$. Since $\Gamma$ is nonelementary it follows that $\bar{\Gamma}$ is either $\mathrm{PSL}(2, \mathbb{C})$ or it preserves a round circle $C \subset \mathbb{S}^2$ and $\Lambda(\Gamma) = C$ [Gr]. If the latter case occurred, $f(\mathbb{H}^2)$ would be one of the two round disks in $\mathbb{S}^2$ bounded by $C$. It would follow that $\tau$ is a branched hyperbolic structure in contradiction to our assumption. If $\bar{\Gamma} = \mathrm{PSL}(2, \mathbb{C})$ then $f(\mathbb{H}^2)$ is contained in $\Lambda(\Gamma) = \mathbb{S}^2$ which implies that $f(\mathbb{H}^2) = \mathbb{S}^2$.

We conclude that if (i) holds then $\Gamma$ is a discrete subgroup of $\mathrm{PSL}(2, \mathbb{C})$ and $f(\mathbb{H}^2)$ is contained in the discontinuity domain $\Omega(\Gamma) = \mathbb{S}^2 \setminus \Lambda(\Gamma)$. Hence (i) $\Rightarrow$ (iii). Clearly, (iii) $\Rightarrow$ (i).

The implication (ii) $\Rightarrow$ (i) is immediate. Conversely if (iii) holds, $f(\mathbb{H}^2)$ must be contained in a component $\Delta$ of the domain of discontinuity of $\Gamma$. Since $f(\mathbb{H}^2)$ is connected and $\Gamma$-invariant it follows that $\Delta$ is also $\Gamma$-invariant. Hence $f$ projects to a holomorphic map $\hat{f} : R \to \hat{f}(R) \subset \Sigma = \Delta/\Gamma$. Since $\hat{f}(R)$ is a compact subsurface without boundary in $\Sigma$ we conclude that $\hat{f}(R) = \Sigma$ and $\Sigma$ is a closed surface. Any nonconstant holomorphic surjective mapping between closed Riemann surfaces is necessarily a covering, possibly branched. Consequently $f$ itself is a possibly branched covering map. $\qquad \square$

We will now construct an example of a singly branched hyperbolic structure on a surface $R$ of genus two which has nondiscrete holonomy in $\mathrm{PSL}(2, \mathbb{R})$.

Start with a regular hyperbolic octagon $X \subset \mathbb{H}^2$ with vertex angles $\pi/2$ (cf. [Tan]). Label the edges $b_1^{-1}, a_1, b_1, a_1^{-1}, \ldots a_2^{-1}$ in positive order around $X$. Identify the edges by corresponding isometries $A_1, B_1, A_2, B_2$ to obtain a Riemann surface of genus two such that $\mathbb{H}^2$ is a two sheeted cover branched over one point on $R$. Let $\sigma$ denote the line segment from the left end point of $b_1^{-1}$ to the right end point of $a_1^{-1}$. Then

$$A_1 B_1 A_1^{-1} B_1^{-1} = E = A_2 B_2 A_2^{-1} B_2^{-1}$$

where $E$ is a elliptic transformation of order two fixing the midpoint of $\sigma$. Let $\gamma$ denote the branched projective structure on $R$ with the holonomy group $\Gamma = \langle A_1, B_1, A_2, B_2 \rangle$. The quotient orbifold $\mathbb{H}^2/\Gamma$ is a torus with one cone point of order two. Clearly the holonomy $\theta : \pi_1(R) \to \Gamma$ is not injective (cf. [Go1]). According to Theorem 1.1.1, $\theta$ does not lift to $\mathrm{SL}(2, \mathbb{R})$.

Next, we will show there exists a hyperbolic structure with exactly one branch point of order two and a *nondiscrete* holonomy group. Take the example above of a branched structure $\gamma$. The representation variety

$$\mathrm{Hom}(\Gamma, \mathrm{PSL}(2, \mathbb{R}))/\mathrm{PSL}(2, \mathbb{R})$$

is 2-dimensional and the representation variety

$$\mathrm{Hom}(\pi_1(R), \mathrm{PSL}(2, \mathbb{R}))/\mathrm{PSL}(2, \mathbb{R})$$

is 6-dimensional. Therefore we can find a real-analytic curve of nonelementary representations $\theta_t : \pi_1(R) \to \mathrm{PSL}(2, \mathbb{R})$, $\theta_0 = \theta, t \in [0,1]$, which do not factor through $\theta : \pi_1(R) \to \Gamma$. The fact that $\theta_t$ is a real-analytic curve implies that there is a dense subset $J \subset [0,1]$ so that $K = \ker(\theta_t) = \ker(\theta_s), s, t \in J$. Let $\Gamma' := \pi_1(R)/K$. We claim that there cannot be a sequence of $t \in J$ which converge to $t = 0$ such that each $\Gamma_t := \theta_t(\pi_1(R))$ is discrete. For otherwise a sequence of discrete nonelementary representations $\rho_t : \Gamma' \to \Gamma_t, t \in J$ would converge to $\rho : \Gamma' \to \Gamma$ as $t \to 0$. The limit $\rho$ of such sequence has to be a faithful representation as well, as a consequence of [J-K]. This contradicts the fact that $\ker(\rho) = \ker(\theta)/K \neq \{1\}$. Thus there is an infinite sequence of nondiscrete representations $\theta_t : \pi_1(R) \to \Gamma_t$ which converges to $\theta$. In addition $\Gamma_t$ necessarily preserves the upper halfplane for $t$ close to 0. By Corollary 11.5.6, $\theta_t$ is the monodromy of a branched complex projective structure $\gamma_t$ on R with branch point likewise at $z = 0$.

## 12. Open questions about complex projective structures

In this chapter we list some unsolved problems. Some are well known in the field, others arise from the specific analysis of this paper.

There are two general issues: the monodromy representation *per se*, and the Riemann surfaces of specified type where it is induced by a particular projective structure.

We recall from §1.5 that $Q_g$ denotes the vector bundle of quadratic differentials over Teichmüller space $\mathfrak{T}_g$ and $V'_g$ is the subset of nonelementary representations in the representation variety $V_g$, modulo conjugation by $\mathrm{PSL}(2, \mathbb{C})$.

12.1. *Existence and nonuniqueness of points in $Q_g$ with given monodromy.* Our proof exhibits two sources of nonuniqueness:

- The nonuniqueness of the pants decomposition on which the monodromy is Schottky.

- The nonuniqueness of the pants configuration over $\mathbb{S}^2$ obtained from a pants decomposition: one can use $N$-sheeted branched covers for arbitrarily large $N$.

Our Theorem 1.1.1 provides a Riemann surface for every nonelementary representation $\theta$. On the other hand, if we fix attention on a particular oriented surface $R$, we don't know whether all projective structures on $R$ itself can be obtained from the pants decomposition method. For example, can there be a complex projective structure $\sigma$ on $R$ so that for each simple loop $\gamma \subset R$ with loxodromic monodromy, *no* element of its homotopy class is sent by the developing map to a simple arc in $\mathbb{S}^2$?

For the case of representations into $\mathrm{PSL}(2, \mathbb{R})$ all projective structures can be obtained by the pants decomposition method, see [F], [Go1], [Ga2]. However in all three papers the proofs that the developing map is a covering over the upper and the lower half-planes have the same gap: In general the pull-back of a complete Riemannian metric on a manifold via a local diffeomorphism can be incomplete. For complete proofs of the assertion about covering see [Kui, pp. 485–486], [Kul-Pin], or [Cho-L]).

For those projective structures on $R$ which do arise from pants decompositions, are there optimal choices for the decompositions? For example, does the developing mapping send each pants of some decomposition directly into the domain of discontinuity of the corresponding Schottky group?

PROBLEM 12.1.1. *Characterize and classify the nonuniqueness of projective structures with given monodromy.*

In particular is it possible to get one projective structure on $R$ from another by a specific series of "moves"?

One might ask to do this through a sequence of graftings. Yet, at least in the case of a once-punctured torus $R$, a connection solely by means of a grafting sequence is known to be impossible in general. The reason has to do with the fact that in the Bers slice, the result of pinching $R$ along a simple nondividing loop $\gamma$ is a $B$-group $\Gamma$ representing the punctured torus on one side, and the triply punctured sphere on the other. Specifically, construct two complex-projective structures on $R$ with the monodromy $G \to \Gamma$ as follows. Consider simple nondividing loops $\alpha$ and $\beta$ on the surface $R$ so that all the loops $\alpha, \beta, \gamma$ are mutually non-homologous. Let $\sigma_t, t \in [0, 1)$ denote the family of complex-projective structures on $R$ which is being pinched along $\gamma$ as $t \to 1$. Let $gr_\alpha(\sigma_t), gr_\beta(\sigma_t)$ be the complex-projective structures on $R$ obtained from $\sigma_t$ via grafting along $\alpha$ and $\beta$. One can show that $gr_\alpha(\sigma_t), gr_\beta(\sigma_t)$ are convergent to complex-projective structures $\sigma_1', \sigma_1''$ on $R$ as $t$ approaches 1. There results two structures $\sigma_1', \sigma_1''$ with the same orientation and the same monodromy $G \to \Gamma$. However these complex projective structures are not related by grafting. The underlying reason is that the "complex of simple loops" on the once punctured torus $R$ is totally disconnected.

For branched structures, there is another way of changing projective structures without changing the monodromy. This is the method of "bubbling."

Suppose that $R$ is a Riemann surface with a (branched) projective structure $\sigma$. Let $\alpha \subset R$ be a compact simple arc, disjoint from the singular points of $\sigma$, which the developing map sends to simple arcs in $\mathbb{S}^2$. Let $a$ be one of these arcs in $\mathbb{S}^2$. Then split $R$ open along $\alpha$, split $\mathbb{S}^2$ open along $a$, take $N$ copies of the Riemann surface $\mathbb{S}^2 - a$ and glue them to $R - \alpha$ with appropriate identification of boundary edges. The net result is a projective structure on the new

"bubble-on" Riemann surface $R_N$ with the same monodromy. The projective structure on $R_N$ has two additional branch points (at the end points of $\alpha$), both of order $N$.

"Bubble-off" is the inverse operation on $R$.

PROBLEM 12.1.2.   *Suppose that $\sigma_1, \sigma_2$ are complex-projective structures on a surface $R$ with the same monodromy representation. Can one pass from $\sigma_1$ to $\sigma_2$ using the following elementary moves*: *"grafting," inverse to "grafting," "bubble-on," "bubble-off"*?

12.2.  *Surfaces with punctures.*  What about surfaces with punctures where the corresponding quadratic differentials have at most double poles? As with compact surfaces, the dimension of the vector bundle $Q_{(g,n)}$ of quadratic differentials over the Teichmüller space $\mathfrak{T}_{(g,n)}$ agrees with that of the representation variety, if one allows arbitrary monodromy at the punctures (for an analysis of the derivative of the monodromy map for this case see [Luo]).  One can search again for pants decompositions, provided the monodromy is not elliptic of infinite order at a punctures.  With discrete monodromy at the punctures, one can look for representations of fundamental groups of pants to extended forms of Schottky groups (i.e. Klein combinations of pairs of discrete cyclic subgroups of $\mathrm{PSL}(2,\mathbb{C})$).

Suppose the genus of $R$ is positive.  We believe that our technique in Part A will yield a pants decomposition of $R$ in which the restrictions of the monodromy are onto Schottky-like groups, provided the representation around each puncture is a discrete (cyclic) group.

PROBLEM 12.2.1.   *Prove and/or explore the existence and nonuniqueness of complex projective structures with given nonelementary monodromy in the case of punctures, most importantly and most classically, punctured spheres.*

12.3.  *Linear monodromy representations.*  Throughout the paper we considered Schwarzian differential equations on Riemann surfaces.  Their monodromy representations are *projective representations* $\theta : \pi_1(R) \to \mathrm{PSL}(2,\mathbb{C})$.

One can also consider the more general case of representations into $\mathrm{GL}(2,\mathbb{C})$.  In the classical case of punctured spheres $R$, the dimension of the representation variety, modulo conjugations, is identical to the dimension of the vector bundle over $\mathfrak{T}_{(0,n)}$ of linear equations

$$u'' + pu' + q = 0,$$

where $p$ has at most simple poles and $q$ double poles at the punctures. Note that we have to restrict to the representations $\theta^*$ into $\mathrm{GL}(2,\mathbb{C})$ which map the peripheral loops of $R$ to unipotent elements.

PROBLEM 12.3.1.    *Is there an analogue of Theorem 1.1.1 for punctured spheres if one seeks a differential equation that induces a given linear representation $\theta^*$?*

12.4.    *Divergence of monodromy representations.*  Fix a closed Riemann surface $R$ of the genus $g > 1$ and let $\phi_n = \varphi_n(z)dz^2$ be a sequence of quadratic differentials on $R$ so that $||\phi_n|| \to \infty$. Let $[\rho_n]$ be the sequence of conjugacy classes of monodromy representations of $\phi_n$. We know from Theorem 11.4.1 that the sequence $[\rho_n]$ cannot subconverge to the to the conjugacy class of *any* representation.

PROBLEM 12.4.1.    *Characterize the "limit points" of divergent sequences of representations in the representation variety. Prove the Divergence Theorem 11.4.1 for complex projective structures on $R$ which have a single branch point of order* 2.

One way that the representation variety $V_g$ can be compactified is by (projective classes of) actions of the group $G = \pi_1(R)$ on metric trees. Which actions of $G$ on trees can appear as limits of the sequences $[\rho_n]$? For instance, is it true that for each sequence of quadratic differentials $\phi_n = n\phi, \phi \neq 0$, there is a limit $\rho$ of the sequence $\rho_n$ with the following property: $\rho$ is an action of $G$ on a tree that is dual to the singular foliation on $R$ determined by $\phi$?

12.5.    *Path lifting properties of monodromy mappings.*  In [He1], Hejhal proved that the natural mapping

$$P_g : Q_g \to V_g$$

is a local homeomorphism which fails to be a covering mapping.

PROBLEM 12.5.1.    *Let $\gamma : [0,1] \to V_g$ be a continuous path, $\tilde{\gamma} : [0,1) \to Q_g$ a partial lift which can not be extended to the end-point* 1. *Describe the asymptotic behavior of the path $\tilde{\gamma}$.*

For instance is it true that $\tilde{\gamma}$ has a well-defined limit

$$\lim_{t \to 1} \tilde{\gamma}(t)$$

in a natural (e.g. closed ball) compactification of $Q_g$?

12.6.    *Branched projective structures.*  As the degree of a positive divisor $D$ increases, it becomes easier to construct a complex projective structure with the branching divisor $D$. Thus, one should be able to eliminate the assumption that the representation $\theta$ is nonelementary for sufficiently large values of $\deg(D)$. For instance, if $\theta$ is the trivial representation, then branched structures with the monodromy $\theta$ are just $m$-fold ramified coverings $f : R \to \mathbb{S}^2$. Thus $\chi(R) = m\chi(\mathbb{S}^2) - \deg(D) = 2m - \deg(D)$. The number $m$ is at least 2, hence $\deg(D) \geq 4 - \chi(R) = 2g + 2$. The minimal degree is realized by a hyperelliptic ramified covering $f$, for which we have: $\deg(D) = 2g + 2$.

PROBLEM 12.6.1. *Make precise and optimize the connection between branching divisors and monodromy. Namely, compute the function*

$$d : \mathrm{Hom}(G, \mathrm{PSL}(2, \mathbb{C})) \to \mathbb{Z},$$

*where $d(\theta)$ is the smallest integer for which there exists a branched complex projective structure with branching divisor of degree $d$ and monodromy $\theta$.*[3]

We proved that $d(\theta) = 0$ for all liftable nonelementary representations $\theta$ and $d(\theta) = 1$ for all nonliftable nonelementary representations $\theta$. Is it true that $d(\theta) = 2g$ for all liftable representations $\theta : G \to \mathrm{SO}(3) \subset \mathrm{PSL}(2, \mathbb{C})$ and $d(\theta) = 2g - 1$ for all nonliftable elementary representations $\theta : G \to \mathrm{SO}(3) \subset \mathrm{PSL}(2, \mathbb{C})$, provided that the monodromy group $\theta(G)$ is dense in $\mathrm{SO}(3)$? Is it true that $d(\theta) \leq 2g + 2$ for any $\theta : G \to \mathrm{PSL}(2, \mathbb{C})$?

*Remark* 12.6.2. For the flat holomorphic bundles of rank 2 over $R$ the corresponding question is the following: given a representation $\theta^* : G \to \mathrm{SU}(2)$ with dense image, find a complex structure on $R$ so that the associated flat $\mathbb{C}^2$-bundle $V$ over $R$ has the degree of instability $u(V) = -1$.

St. John's University, Staten Island, NY
*E-mail address*: gallod@stjohns.edu

University of Utah, Salt Lake City, UT
*E-mail address*: kapovich@math.utah.edu

University of Minnesota, Minneapolis, MN
*E-mail address*: am@math.umn.edu

---

[3]This was recently solved by M. Kapovich.